%% file: AMain.tex
\documentclass[12pt]{amsart}

\input{Preamble}

\usepackage{standalone} 

\title{Higher Sheaf Theory I: Correspondences}

\author{G. Stefanich}

\date{}

\begin{document}


\begin{abstract}
We prove a universal property for the $(\infty,n)$-category of correspondences, generalizing and providing a new proof for the case $n = 2$ from \cite{GR}. We also provide conditions under which a functor out of a higher category of correspondences of $\Ccal$ can be extended to a higher category of correspondences of the free cocompletion of $\Ccal$. These results will be used in the sequels to this paper to construct $(\infty,n)$-categorical versions of the theories of quasicoherent and ind-coherent sheaves in derived algebraic geometry.
\end{abstract}

\maketitle

\tableofcontents

\input{Introduction}

\input{"Two_sided"}

\input{"Two_corr"}

\input{"Higher_correspondences"}

\bibliographystyle{amsalpha}
\bibliography{References}
 
\end{document}

%% file: Preamble.tex
\usepackage{amssymb, amsmath, amsthm}
\usepackage{calligra, mathrsfs}
\usepackage{graphicx}
\usepackage{url}
\usepackage{mathtools}
\usepackage{enumerate}
\usepackage{verbatim}
\usepackage[retainorgcmds]{IEEEtrantools}
\usepackage{tikz-cd}
\usetikzlibrary{positioning}

\usepackage[margin=1in,marginparwidth=0.8in, marginparsep=0.1in]{geometry}

\newcommand{\Acal}{{\mathcal A}}
\newcommand{\Bcal}{{\mathcal B}}
\newcommand{\Ccal}{{\mathcal C}}
\newcommand{\Dcal}{{\mathcal D}}
\newcommand{\Ecal}{{\mathcal E}}
\newcommand{\Fcal}{{\mathcal F}}

\newcommand{\Hcal}{{\mathcal H}}
\newcommand{\Ical}{{\mathcal I}}
\newcommand{\Jcal}{{\mathcal J}}

\newcommand{\Mcal}{{\mathcal M}}

\newcommand{\Ocal}{{\mathcal O}}
\newcommand{\Pcal}{{\mathcal P}}

\newcommand{\Scal}{{\mathcal S}}
\newcommand{\Tcal}{{\mathcal T}}
\newcommand{\Ucal}{{\mathcal U}}

\newcommand{\NN}{{\mathbb N}}

\usepackage{scalerel,stackengine}
\stackMath
\newcommand\reallywidehat[1]{%
\savestack{\tmpbox}{\stretchto{%
  \scaleto{%
    \scalerel*[\widthof{\ensuremath{#1}}]{\kern-.6pt\bigwedge\kern-.6pt}%
    {\rule[-\textheight/2]{1ex}{\textheight}}
  }{\textheight}%
}{0.5ex}}%
\stackon[1pt]{#1}{\tmpbox}%
}

\makeatletter
\newcommand*\bigcdot{{\mathpalette\bigcdot@{.5}}}
\newcommand*\bigcdot@[2]{\mathbin{\vcenter{\hbox{\scalebox{#2}{$\m@th#1\bullet$}}}}}
\makeatother

\newcommand{\dsh}{{\text{\normalfont-}}}

\makeatletter
\g@addto@macro\bfseries{\boldmath}
\makeatother

\tikzset{shorten <>/.style={shorten >=#1,shorten <=#1}}

\newcommand{\kr}{\kern -2pt}

\DeclareMathOperator{\id}{id}
\DeclareMathOperator{\Hom}{Hom}

\DeclareMathOperator{\End}{End}
\DeclareMathOperator{\Funct}{Funct}
\DeclareMathOperator{\ev}{ev}

\newcommand{\iso}{{\text{\normalfont iso}}}
\newcommand{\fin}{{\text{\normalfont fin}}}

\newcommand{\univ}{{\text{\normalfont univ}}}

\newcommand{\red}{{\text{\normalfont red}}}

\newcommand{\op}{{\text{\normalfont op}}}

\newcommand{\el}{{\text{\normalfont el}}}
\DeclareMathOperator{\spn}{sp}
\DeclareMathOperator{\codiscr}{codisc}

\DeclareMathOperator{\Spc}{Spc}

\DeclareMathOperator{\Cat}{Cat}
\newcommand{\twoCat}{{2\kr\Cat}}
\newcommand{\nCat}{{n\kr\Cat}}
\DeclareMathOperator{\omegaCat}{\omega\kr\Cat}

\DeclareMathOperator{\colim}{colim}

\newcommand{\pb}{{\text{\normalfont pb}}}
\newcommand{\lex}{{\text{\normalfont lex}}}

\newcommand{\cart}{{\text{\normalfont cart}}}
\newcommand{\cocart}{{\text{\normalfont cocart}}}
\newcommand{\bicart}{{\text{\normalfont two-sided}}}
\newcommand{\bivar}{{\text{\normalfont bivar}}}

\DeclareMathOperator{\Tw}{Tw}

\DeclareMathOperator{\Cocyl}{Cocyl}

\DeclareMathOperator{\Adj}{Adj}

\let\Pr\relax
\DeclareMathOperator{\Pr}{Pr}

\DeclareMathOperator{\Sp}{Sp}
\newcommand{\St}{{\text{\normalfont St}}}

\DeclareMathOperator{\Fin}{Fin}
\DeclareMathOperator{\Assos}{Assos}
\DeclareMathOperator{\BM}{BM}

\DeclareMathOperator{\Alg}{{Alg}}
\DeclareMathOperator{\CAlg}{CAlg}

\DeclareMathOperator{\Algbrd}{Algbrd}
\newcommand{\algbrd}{{\text{\normalfont algbrd}}}
\newcommand{\Catenr}[1]{\Cat^{#1}}

\DeclareMathOperator{\Corr}{Corr}
\newcommand{\twoCorr}{{2\kr \Corr}}
\newcommand{\nCorr}{{n\kr\Corr}}

\newcommand*{\shom}{{\mathscr{H}\kern -2pt om}}
\newcommand*{\Catscr}{{\mathscr{C}\kern -2pt at}}
\newcommand*{\Prscr}{{\mathscr{P}\kern -2pt r}}
\newcommand{\shext}{{\mathscr{E} \kern -2pt xt}}

\newtheorem{definition}[subsubsection]{Definition}
\newtheorem{proposition}[subsubsection]{Proposition}
\newtheorem{lemma}[subsubsection]{Lemma}
\newtheorem{theorem}[subsubsection]{Theorem}
\newtheorem{corollary}[subsubsection]{Corollary}

\newtheorem*{theorem*}{Theorem}
\newtheorem*{conjecture*}{Conjecture}

\newtheoremstyle{claim}{\topsep}{\topsep}{}{}{\it}{.}{.5em}{}
\theoremstyle{claim}
\newtheorem*{claim*}{Claim}

\newtheoremstyle{note}
{8.0pt plus 2.0pt minus 4.0pt}{8.0pt plus 2.0pt minus 4.0pt}{}{}{\bf}{.}{.5em}{} 
\theoremstyle{note}
\newtheorem{example}[subsubsection]{Example}
\newtheorem{remark}[subsubsection]{Remark}
\newtheorem{notation}[subsubsection]{Notation}
\newtheorem{construction}[subsubsection]{Construction}

\newtheorem*{example*}{Example}

\AtBeginDocument{\def\MR#1{}}

%% file: Introduction.tex

\tableofcontents

\section{Introduction}

This paper is the first of a series on the foundations of the theory of sheaves of $(\infty,n)$-categories in derived algebraic geometry. The goal of this series is to construct and study the theories of quasicoherent and ind-coherent sheaves of $(\infty,n)$-categories on prestacks, and use them to produce interesting examples of fully extended topological field theories. 

 We begin in this article by studying a general framework for higher categorical sheaf theories in the language of categories of correspondences. As emphasized in \cite{GR}, much of the functorial behaviour of a sheaf theory is encoded in its realization as a functor out of an $(\infty,2)$-category of correspondences. The existence of this realization has direct consequences in terms of dualizability properties of the sheaf theory, and is a helpful tool in the computation of various field-theoretic invariants \cite{BZN2}.
 
  We therefore take the point of view that a higher sheaf theory on a category of geometric objects $\Ccal$ is a (possibly symmetric monoidal) functor out of a higher category of correspondences on $\Ccal$, into a target $(\infty,n)$-category $\Dcal$ whose objects are to be thought of as being $(\infty, n-1)$-categories of some sort. The main goal of this paper is to provide tools for the construction of such functors.
  
 Below we provide a more detailed description of the contents of this paper. We will use the convention where all objects are $\infty$-categorical by default, and suppress this from our notation from now on.
  
  \subsection{Sheaf theories and the $2$-category of correspondences}
  
  We begin by reviewing the case of ordinary, $1$-categorical sheaf theory. Let $\Ccal$ be a category admitting finite limits, whose objects we think about as geometric spaces of some sort (for instance, $\Ccal$ can be the category of affine schemes, stacks, etc). Then one can define a $2$-category $\twoCorr(\Ccal)$ with the following properties:
  \begin{itemize}
  \item The objects of $\twoCorr(\Ccal)$ are the objects of $\Ccal$.
  \item For each pair of objects $c, c'$ in $\Ccal$, the category $\Hom_{\twoCorr(\Ccal)}(c, c')$ is the category $\Ccal_{/c, c'}$ whose objects are spans $c \leftarrow s \rightarrow c'$ in $\Ccal$.
  \item The composition of two  spans $c \leftarrow s \rightarrow c'$ and  $c' \leftarrow t \rightarrow c''$ is given by $c \leftarrow s\times_{c'} t \rightarrow c''$.
\end{itemize}

In its most basic form \footnote{Many sheaf theories of interest can be formulated in this form if we allow ourselves to work with variants of $\twoCorr(\Ccal)$ where the legs of the spans are required to belong to certain predetermined classes, and change the notion of $2$-cell. For our purposes, this basic version will suffice.}, a sheaf theory on $\Ccal$ consists of a functor $F: \twoCorr(\Ccal) \rightarrow \Dcal$ into a $2$-category $\Dcal$ whose objects are to be thought of as categories of some sort. This assigns to each object $c$ in $\Ccal$ an object $F(c)$ in $\Dcal$, subject to the following functoriality:
\begin{itemize}
\item For every morphism $\alpha: c \rightarrow c'$ in $\Ccal$ a morphism $\alpha_!: F(c) \rightarrow F(c')$ associated to the following span:
\[
\begin{tikzcd}
  & c \arrow[ld, "\id"'] \arrow[rd, "\alpha"] &    \\
c &                                          & c'
\end{tikzcd}
\]
\item For every morphism $\alpha: c \rightarrow c'$ in $\Ccal$ a morphism $\alpha^!: F(c') \rightarrow F(c)$ associated to the following span:
\[
\begin{tikzcd}
  & c \arrow[rd, "\id"] \arrow[ld, "\alpha"'] &    \\
c' &                                          & c
\end{tikzcd}
\]
\item For every morphism $\alpha: c\rightarrow c'$, a $2$-cell $\alpha_!\alpha^! \rightarrow \id_{c'}$, associated to the following morphism of spans:
\[
\begin{tikzcd}
   & c \arrow[ldd, "\alpha"'] \arrow[rdd, "\alpha"] \arrow[d,"\alpha"] &    \\
   & c' \arrow[ld, "\id"] \arrow[rd, "\id"']                  &    \\
c' &                                                          & c'
\end{tikzcd}
\]
\end{itemize}

It can be shown that the $2$-cell from the third item is the counit of an adjunction between $\alpha_!$ and $\alpha^!$. Furthermore, the composition rule of $\twoCorr(\Ccal)$ implies that for every pair of maps $\alpha: c\rightarrow c' \leftarrow c'' : \beta$ the square
\[
\begin{tikzcd}
F(c \times_{c'} c'')  \arrow{d}{\beta'_!} & F(c'') \arrow{l}[swap]{\alpha'^!} \arrow{d}{\beta_!} \\ F(c) & F(c') \arrow{l}[swap]{\alpha^!}
\end{tikzcd}
\]
which is in principle only commutative up to a natural transformation, is in fact strictly commutative.

There is a natural inclusion $\iota_\Ccal: \Ccal \rightarrow \twoCorr(\Ccal)$, which is the identity on objects, and sends each arrow $\alpha$ to the span with right leg $\alpha$ and identity left leg. It was proven in \cite{GR} that the functor $F$ can be recovered just from the knowledge of its restriction to $\Ccal$. Moreover, given a functor $f: \Ccal \rightarrow \Dcal$ one can find an extension of $F$ to $\twoCorr(\Ccal)$ if and only if $f(\alpha)$ admits a right adjoint for every arrow $\alpha$ in $\Ccal$, and for every pair of maps $\alpha: c \rightarrow c' \leftarrow c'' : \beta$ in $\Ccal$ the associated lax commutative square is strictly commutative:

\begin{theorem}[\cite{GR}] \label{teo 2corr}
Precomposition with $\iota_\Ccal$ induces an equivalence between the space of functors $F: \twoCorr(\Ccal) \rightarrow \Dcal$ and the space of functors $f: \Ccal \rightarrow \Dcal$ satisfying the left Beck-Chevalley condition.
\end{theorem}

Our first goal in this paper is to give an alternative approach to this result. Our proof relies on two main ideas:
\begin{itemize}
\item Let $\twoCorr^{\text{univ}}(\Ccal)$ be the $2$-category satisfying the universal property of theorem \ref{teo 2corr}. Then the restriction of the hom functor 
\[
\Hom_{\twoCorr^{\text{univ}}(\Ccal)}: \twoCorr^{\text{univ}}(\Ccal)^{1\dsh\op} \times \twoCorr^{\text{univ}}(\Ccal) \rightarrow \Catscr
\]
to $\twoCorr^{\text{univ}}(\Ccal)^{1\dsh\op} \times \Ccal$ is shown to satisfy a universal property - roughly speaking, we show that it can be obtained by left Kan extension of  $\Hom_\Ccal: \Ccal^\op \times \Ccal \rightarrow \Spc$ along the inclusion of $\Ccal^\op \times \Ccal$ inside $\twoCorr^{\text{univ}}(\Ccal)^{1\dsh\op} \times \Ccal$. Our proof of this relies on the description of Hom functors for enriched categories in terms of diagonal bimodules from \cite{Hinich}.
\item We prove a version of the Grothendieck construction which relates functors $\Ccal^\op \times \Ccal \rightarrow \Cat$ and the so-called two-sided fibrations over $\Ccal \times \Ccal$, and identify the image of $\Funct(\twoCorr^{\text{univ}}(\Ccal)^{1\dsh\op} \times \Ccal, \Catscr)$ under this equivalence. The problem then gets reduced to establishing a universal property for the two-sided fibration whose objects are spans in $\Ccal$. Our proof of this fact ultimately relies on the description of free fibrations from \cite{GHN}. 
\end{itemize}

\subsection{Higher sheaf theories and the $n$-category of correspondences}

For each $n > 2$ one can define an $n$-category $\nCorr(\Ccal)$ with the following properties:
\begin{itemize}
\item The objects of $\nCorr(\Ccal)$ are the objects of $\Ccal$.
\item For each pair of objects $c, c'$ we have $\Hom_{\nCorr(\Ccal)}(c, c') = (n-1)\kr\Corr(\Hom_{\twoCorr(\Ccal)}(c, c'))$.
\end{itemize}

In its most basic form, a higher sheaf theory on $\Ccal$ is a functor $F: \nCorr(\Ccal) \rightarrow \Dcal$ into an $n$-category $\Dcal$ whose objects we think about as being $(n-1)$-categories of some sort. Such a functor satisfies a large list of adjointability and base change properties. As before, its restriction to $\Ccal$ satisfies the left Beck-Chevalley condition. However, as soon as $n > 2$ it turns out that it also satisfies the right Beck-Chevalley condition: in fact, pushforwards and pullbacks form part of an ambidextrous adjunction. Furthermore, when $n > 2$ we have, for each pair of objects $c, c'$ in $\Ccal$, a functor
\[
F_* : \Ccal_{/c, c'} \rightarrow \Hom_\Dcal(F(c), F(c'))
\]
which also satisfies the left Beck-Chevalley condition. When $n > 3$, these functors also satisfy the right Beck-Chevalley condition, and for every pair of spans $S = (c \leftarrow s \rightarrow c')$ and $T = (c \leftarrow t \rightarrow c')$ there is an induced functor
\[
(F_*)_*: \Ccal_{/s, t} \rightarrow \Hom_{\Hom_\Dcal(F(c), F(c'))}(F_*(S), F_*(T))
\]
which also satisfies the left Beck-Chevalley condition.

Our first main result in this paper is a generalization of theorem \ref{teo 2corr} which singles out a minimalistic list of base change properties to verify for a functor $\Ccal \rightarrow \Dcal$ to admit an extension to $\nCorr(\Ccal)$:

\begin{theorem}\label{teo univ prop ncorr}
Precomposition with the inclusion $\Ccal \rightarrow \nCorr(\Ccal)$ induces an equivalence between the space of functors $\nCorr(\Ccal) \rightarrow \Dcal$ and the space of functors $\Ccal \rightarrow \Dcal$ satisfying the left $(n-1)$-fold  Beck-Chevalley condition.
\end{theorem}

In the same way that verifying left Beck-Chevalley condition for a functor $f: \Ccal \rightarrow \Dcal$ involves checking an adjointability statement for every cartesian square in $\Ccal$, verifying the left $(n-1)$-fold Beck-Chevalley condition involves checking a series of $n-1$ different adjointability statements for every such cartesian square. We refer the reader to section \ref{subsection higher BC} for a precise definition of the left $(n-1)$-fold Beck-Chevalley condition.

The $n$-category $\nCorr(\Ccal)$ comes equipped with a symmetric monoidal structure, and the inclusion $\Ccal \rightarrow \nCorr(\Ccal)$ can be enhanced to a symmetric monoidal functor, where $\Ccal$ is given the cartesian symmetric monoidal structure. In this paper we also prove a version of theorem \ref{teo univ prop ncorr} that takes into account this structure:

\begin{corollary}\label{prop univ con sym mon}
Let $\Dcal$ be a symmetric monoidal $n$-category. Then restriction along the inclusion $\Ccal \rightarrow \nCorr(\Ccal)$ induces an equivalence between the space of symmetric monoidal functors $\nCorr(\Ccal) \rightarrow \Dcal$ and the space of symmetric monoidal functors $\Ccal \rightarrow \Dcal$ which satisfy the left $(n-1)$-fold Beck-Chevalley condition.
\end{corollary}

 It was shown in  \cite{HaugSpan} that every object of $\nCorr(\Ccal)$ is fully dualizable in the $(n-1)$-category underlying $\nCorr(\Ccal)$. Combining this fact with corollary \ref{prop univ con sym mon} we are able to conclude that if $F: \Ccal \rightarrow \Dcal$ is a symmetric monoidal functor satisfying the $(n-1)$-fold left Beck-Chevalley condition and $c$ is an object of $\Ccal$, then $F(c)$ is a fully dualizable object in the $(n-1)$-category underlying $\Dcal$. In other words, $F(c)$ can be considered as the object of boundary conditions for an $(n-1)$-dimensional topological field theory.

\subsection{Extension along the Yoneda embedding}
The second main result of this paper concerns the interaction of the higher Beck-Chevalley condition and the procedure of Kan extension along the Yoneda embedding $\Ccal \rightarrow \Pcal(\Ccal)$. It gives conditions under which one can deduce that a colimit preserving functor $\Pcal(\Ccal) \rightarrow \Dcal$  satisfies a higher Beck-Chevalley condition, from knowing that its restriction to $\Ccal$ does. 

 This is instrumental in the study of higher sheaf theories in algebraic geometry, as one usually first gives a definition for affine schemes, and then one Kan-extends the theory to more general stacks. It can often be checked without too much difficulty that a higher Beck-Chevalley condition holds on affine schemes, but this is not so easy to do by hand for arbitrary stacks. 

There are two conditions, which we introduced in \cite{Pres}, that are required on $\Dcal$ to be able to do this. The first one is called the passage to adjoints property, and is modeled on the well-known fact that colimits of right adjointable arrows in  $\Pr^L$ can be computed as limits after passage to right adjoints.  This condition is needed to ensure that the image of any arrow in $\Pcal(\Ccal)$ will admit a right adjoint, from knowing that this is true for $\Ccal$. We require that this holds not only for $\Dcal$, but also for all Hom $(n-1)$-categories of $\Dcal$, all Hom $(n-2)$-categories of those, and so on.

The second one is that $\Dcal$ be conically cocomplete. In other words, we want that colimits in the $1$-category underlying $\Dcal$ be compatible with the $n$-categorical structure. We use this to reduce base change properties for cartesian squares in $\Pcal(\Ccal)$ to those in which the final vertex belongs to $\Ccal$. Again this is required inductively on all Hom-categories as well.

The following is a simplified version of our second main result - we refer the reader to section \ref{subsection extension} for a slightly stronger statement:

\begin{theorem}\label{teo ext}
Assume that $\Dcal$ is $(n-1)$-fold conically cocomplete, and satisfies the $(n-1)$-fold passage to adjoints property. Let $F: \Pcal(\Ccal) \rightarrow \Dcal$ be a conical colimit preserving functor and assume that $F|_{\Ccal}$ and $F|_{\Ccal}^{n\dsh\op}$ satisfy the left $(n-1)$-fold  Beck-Chevalley condition. Then $F$ satisfies the left  $(n-1)$-fold Beck-Chevalley condition.
\end{theorem}

We remark that if $\Dcal$ is an $(n+1)$-category and $F|_{\Ccal}$ satisfies the left $n$-fold Beck-Chevalley condition, then $F|_{\Ccal}$ and $F|_{\Ccal}^{n\dsh\op}$ satisfy the left $(n-1)$-fold Beck-Chevalley condition. In other words, under the stated assumptions on $\Dcal$, any functor $(n+1)\kr\Corr(\Ccal) \rightarrow \Dcal$ induces a functor $n\kr\Corr(\Pcal(\Ccal)) \rightarrow \Dcal$.

In our planned applications, $\Dcal$ will be taken to be the $n$-category $\psi_n(D)$ underlying a presentable $n$-category $D$, as defined in \cite{Pres}. We will be particularly interested in the case when $\Dcal$ is the $n$-category $(n-1)\Prscr^L$ of presentable $(n-1)$-categories, or its stable version $(n-1)\Prscr^L_{\St}$. The main results from \cite{Pres} guarantee that this satisfies the hypothesis of theorem \ref{teo ext}. Specializing theorems \ref{teo univ prop ncorr} and \ref{teo ext} to this context yields the following result:

\begin{theorem}\label{teo extension presentables}
Let $\Ccal$ be a category admitting pullbacks and let $\Dcal$ be the $n$-category underlying a presentable $n$-category. Let $F:\Ccal \rightarrow \Dcal$ be a functor such that $F|_{\Ccal}$ and $F|_{\Ccal}^{n\dsh\op}$ satisfy the left $(n-1)$-fold Beck-Chevalley condition. Then there exists a unique extension of $F$ to a functor $\nCorr(\Pcal(\Ccal)) \rightarrow \Dcal$ whose restriction to $\Pcal(\Ccal)$ preserves conical colimits.
\end{theorem}

In the sequels to this paper we will let $\Ccal$ be the category of affine schemes, and use theorem \ref{teo extension presentables} to construct higher sheaf theories on prestacks.

\subsection{Organization} We now describe the contents of the paper in more detail. We refer the reader also to the introduction of each section for an expanded outline of its contents.

In section 2 we study the theory of two-sided fibrations. We introduce a two-sided variant of the Grothendieck construction, and study the classes of two-sided fibrations which are classified by functors with various adjointability properties. We discuss two main examples: the arrow fibration and the span fibration, and prove universal properties  for them.

Section 3 deals with the $2$-category of correspondences. We recall its definition, and collect a series of basic functoriality, adjointness and duality results, for later reference. We finish this section by applying the theory of two-sided fibrations to the proof of theorem \ref{teo 2corr}

In section 4 we study the $n$-category of correspondences. We begin by giving a definition adapted to our purposes using the language of enriched category theory, and review its main adjointness and duality properties. We then introduce the higher Beck-Chevalley conditions, and prove theorem \ref{teo univ prop ncorr}. We finish this section with a proof of theorem \ref{teo ext}.

\subsection{Conventions and notation}

We use the language of higher category theory and higher algebra as developed in \cite{HTT} and \cite{HA}. All of our notions will be assumed to be homotopical or $\infty$-categorical, and we suppress this from our notation - for instance, we use the word $n$-category to mean $(\infty,n)$-category. 

We work with a nested sequence of universes. Objects belonging to the first universe are called small, objects in the second universe are called large, and objects in the third universe are called very large.

We denote by $\Spc$ and $\Cat$ the categories of (small) spaces and categories. For each $n \geq 2$ we denote by $\nCat$ the category of small $n$-categories.  We denote by $\Catscr$ the $2$-categorical enhancement of $\Cat$, and in general by $n\Catscr$ the $(n+1)$-categorical enhancement of $\nCat$. If $X$ is one of those objects (or related), we denote by $\widehat{X}$ its large variant. For instance, $\widehat{\Spc}$ denotes the category of large spaces.

We denote by $\Pr^L$ the category of presentable categories and colimit preserving functors. We usually consider this as a symmetric monoidal category, with the tensor product constructed in \cite{HA} section 4.8. By presentable (symmetric) monoidal category we mean a (commutative) algebra in $\Pr^L$. In other words, this is a presentable category equipped with a (symmetric) monoidal structure compatible with colimits.

 If $\Ccal$ is an $n$-category and $k \geq 0$, we denote by $\Ccal^{\leq k}$ the $k$-category obtained from $\Ccal$ by discarding all cells of dimension greater than $k$, and by $\prescript{\leq k}{}{\Ccal}$ the $k$-category obtained from $\Ccal$ by inverting all cells of dimension greater than $k$. In particular, for each category $\Ccal$ we denote by $\Ccal^{\leq 0}$ the space of objects of $\Ccal$.

For each category $\Ccal$ we denote by $\Hom_{\Ccal}(-, -)$ the Hom-functor of $\Ccal$. We denote by $\Pcal(\Ccal)$ its presheaf category. We usually identify $\Ccal$ with its image under the Yoneda embedding.

Given a pair of $n$-categories $\Ccal , \Dcal$ we denote by $\Funct(\Ccal, \Dcal)$ the $n$-category of functors between $\Ccal$ and $\Dcal$. When we wish to only consider the space of functors between them we will use $\Hom_{\nCat}(\Ccal, \Dcal)$ instead.

\subsection{Acknowledgments}

I would like to thank David Nadler for his advice and support these last few years, and in particular for pointing me towards the subject of categorified sheaf theory. I am also grateful to Nick Rozenblyum for a helpful conversation where he shared his thoughts on how various parts of this paper could be approached.


%% file: Two_sided.tex

\tableofcontents
\section{Two-sided fibrations}

Let $\Ccal$ and $\Dcal$ be categories. A two-sided fibration from $\Ccal$ to $\Dcal$ is a functor $p: \Ecal \rightarrow \Ccal \times \Dcal$ which behaves as a cocartesian fibration along the $\Ccal$-directions, and as a cartesian fibration along the $\Dcal$-directions. In \ref{section groth} we review the notion of two-sided fibration, and we present a two-sided analog of the straightening-unstraightening equivalence, which identifies the category of two-sided fibrations over $\Ccal \times \Dcal$ with the category of functors $\Ccal \times \Dcal^\op \rightarrow \Cat$.

In \ref{subsection arrow} we study a fundamental example of a two-sided fibration: the target-source projection $\Funct([1], \Ccal) \rightarrow \Ccal \times \Ccal$. We show that this is classified by the Hom functor $\Hom_\Ccal : \Ccal \times \Ccal^\op \rightarrow \Spc$. This is a two-sided counterpart to the fact that the twisted arrow category of $\Ccal$ is the pairing of categories classified by the Hom functor of $\Ccal$. As an application,  we provide a concrete description of the so-called representable bifibrations, which are classified by functors $\Ccal \times \Dcal^\op \rightarrow \Spc$ of the form $\Hom_\Dcal(F(-), -)$ for some functor $F: \Ccal \rightarrow \Dcal$.  We also establish a universal property for the arrow bifibration - in its most basic form it states that the arrow category is the free bifibration on the diagonal functor $\Delta: \Ccal \rightarrow \Ccal \times \Ccal$.

In \ref{subsection bivar} we introduce the class of bivariant fibrations. These are functors $p: \Ecal \rightarrow \Ccal \times \Dcal$ which are simultaneously cocartesian fibrations, cartesian fibrations, and two-sided fibrations over both $\Ccal \times \Dcal$ and $\Dcal \times \Ccal$. In the same way that a two-sided fibration over $\Ccal \times \Dcal$ is classified by a functor $\Ccal \times \Dcal^\op \rightarrow \Cat$, we show that a bivariant fibration is equivalent to the data of four functors
\begin{align*}
\Ccal \times \Dcal \rightarrow \Cat & \hspace{2cm} \Ccal \times \Dcal^\op \rightarrow \Cat \\\Ccal^\op \times \Dcal \rightarrow \Cat & \hspace{2cm} \Ccal^\op \times \Dcal^\op \rightarrow \Cat
\end{align*}
each of which determines the rest by passage to right or left adjoints in one or both coordinates.

In \ref{subsection Beck Chev} we specialize to the case when $\Ccal$ and $\Dcal$ admit pullbacks, and study the class of bivariant fibrations satisfying the so-called Beck-Chevalley condition. We show that these are classified by functors which satisfy familiar base change properties in each variable. We finish with a fundamental example of a bivariant fibration satisfying the Beck-Chevalley condition: the source-target projection $\Funct(\Lambda^2_0, \Ccal) \rightarrow \Ccal \times \Ccal$, where $\Funct(\Lambda^2_0, \Ccal)$ denotes the category whose objects are spans in $\Ccal$. We show that this enjoys a universal property: it is the free cocartesian and two-sided fibration which satisfies the Beck-Chevalley condition in the first coordinate, on the arrow two-sided fibration.

\subsection{The two-sided Grothendieck construction}\label{section groth}

We begin by reviewing the notion of two-sided fibration.

\begin{definition}
Let $\Ccal, \Dcal, \Ecal$ be categories, and let $p = (p_1, p_2): \Ecal \rightarrow \Ccal \times \Dcal$ be a functor. We say that $p$ is a lax two-sided fibration from $\Ccal$ to $\Dcal$ if the following conditions hold:
\begin{itemize}
\item Let $e$ be an object of $\Ecal$ and write $p(e) = (c, d)$. Then for every arrow $\alpha: c \rightarrow c'$ in $\Ccal$ there is a $p$-cocartesian lift of $(\alpha, \id_d)$ with source $e$.
\item Let $e$ be an object of $\Ecal$ and write $p(e) = (c, d)$. Then for every arrow $\beta: d \rightarrow d'$ in $\Ccal$ there is a $p$-cartesian lift of $(\id_c, \beta)$ with target $e$.
\end{itemize}
\end{definition}

\begin{remark}
The condition that the functor $p$ be a lax two-sided fibration depends on the decomposition of $\Ccal \times \Dcal$ as an ordered product. It will usually be clear from context what this decomposition is. Unless otherwise stated, we will work in this section with the convention where a (lax) two-sided fibration has cocartesian lifts of arrows in the first factor, and cartesian lifts of arrows in the second factor. In other words, unless otherwise stated,  if we say that a functor $p: \Ecal \rightarrow \Ccal \times \Dcal$ is a (lax) two-sided fibration we mean that it is a (lax) two-sided fibration from $\Ccal$ to $\Dcal$.
\end{remark}

\begin{remark}\label{remark lax bicartesian arrows}
Let $\Ccal, \Dcal, \Ecal$ be categories, and let $p = (p_1, p_2): \Ecal \rightarrow \Ccal \times \Dcal$ be a functor. Let $\alpha: c \rightarrow c'$ be an arrow in $\Ccal$ and $d$ be an object in $\Dcal$. It follows from \cite{HTT} proposition 4.3.1.5 item (2) that a lift of $(\alpha, \id_d)$ to $\Ecal$  is $p$-cocartesian if and only if it is $p_1$-cocartesian. Similarly, if $\beta: d \rightarrow d'$ is a morphism in $\Dcal$ and $c$ is an object in $\Ccal$, a lift of $(\id_c, \beta)$ to $\Ecal$ is $p$-cartesian if and only if it is $p_2$-cartesian. 

In particular, we see that $p$ is a lax two-sided fibration if and only if the following two conditions hold:
\begin{itemize}
\item $p_1$ is a cocartesian fibration and $p$ is a morphism of cocartesian fibrations over $\Ccal$.
\item $p_2$ is a cartesian fibration and $p$ is a morphism of cartesian fibrations over $\Dcal$.
\end{itemize}
\end{remark}

\begin{definition}
Let $\Ccal, \Dcal$ be categories, and let $p = (p_1, p_2): \Ecal \rightarrow \Ccal \times \Dcal$ be a lax two-sided fibration from $\Ccal$ to $\Dcal$. Let $f: e \rightarrow e'$ be an arrow in $\Ecal$, lying above an arrow $(\alpha, \beta): (c, d) \rightarrow (c', d')$. Let $\overline{(\alpha, \id_d)}: e \rightarrow m$ be a $p$-cocartesian lift of $(\alpha, \id_d)$, and $\overline{(\id_c', \beta)}: m' \rightarrow e$ be a $p$-cartesian lift of $(\id_c, \beta)$. We say that $f$ is $p$-bicartesian if the induced map $m \rightarrow m'$ is an isomorphism. We say that $p$ is a two-sided fibration from $\Ccal$ to $\Dcal$ if composition of $p$-bicartesian arrows is $p$-bicartesian.
\end{definition}

\begin{remark} \label{remark bicartesian arrows}
Let $\Ccal, \Dcal$ be categories, and let $p = (p_1, p_2): \Ecal \rightarrow \Ccal \times \Dcal$ be a lax two-sided fibration. An arrow $f : e \rightarrow e'$ in $\Ecal$ is $p$-bicartesian if and only if it can be written as a composition of a $p_1$-cocartesian arrow followed by a $p_2$-cartesian arrow. The projection $p$ is a two-sided fibration if and only if any arrow of the form $f_1 \circ f_2$ where $f_2$ is $p_2$-cartesian and $f_1$ is $p_1$-cocartesian, is bicartesian.
\end{remark}

\begin{definition}
Let $\Ccal, \Dcal$ be categories, and let $p: \Ecal \rightarrow \Ccal \times \Dcal$ and $q: \Ecal' \rightarrow \Ccal \times \Dcal$ be two-sided fibrations. A functor $F: \Ecal \rightarrow \Ecal'$ equipped with an identification $qF = p$ is said to be a morphism of two-sided fibrations if it maps bicartesian arrows to bicartesian arrows.
\end{definition}

\begin{remark}
Let $\Ccal, \Dcal$ be categories, and let $p: \Ecal \rightarrow \Ccal \times \Dcal$ and $q: \Ecal' \rightarrow \Ccal \times \Dcal$ be two-sided fibrations. Let $F: \Ecal \rightarrow \Ecal'$ be a functor equipped with an identification $qF = p$. It follows from remark \ref{remark bicartesian arrows} that $F$ is a morphism of two-sided fibrations if and only if it is a morphism of cocartesian fibrations over $\Ccal$ and a morphism of cartesian fibrations over $\Dcal$.
\end{remark}

\begin{notation}
Let $\Ccal$ be a category. We denote by $\Cat^\cart_{/\Ccal}$ (resp. $\Cat^\cocart_{/\Ccal}$) the full subcategory of the overcategory $\Cat_{/\Ccal}$ on the (co)cartesian fibrations and morphisms of (co)cartesian fibrations. Given another category $\Dcal$, we denote by $\Cat^{\bicart}_{/\Ccal, \Dcal}$ the subcategory of $\Cat_{/\Ccal, \Dcal}$ on the two-sided fibrations and morphisms of two-sided fibrations.
\end{notation}

If we allow ourselves to break symmetry, we can give an alternative characterization of two-sided fibrations and morphisms of two-sided fibrations.

\begin{proposition}\label{prop bicart 0}
Let $\Ccal, \Dcal, \Ecal$ be categories, and let $p = (p_1, p_2): \Ecal \rightarrow \Ccal \times \Dcal$ be a  functor. Then
\begin{enumerate}[\normalfont (i)]
\item The map $p$ is a lax two-sided fibration if and only if $p_1$ is a cocartesian fibration, the functor $p$ is a morphism of cocartesian fibrations over $\Ccal$, and for every object $c$ in $\Ccal$ the projection $p_1^{-1}(c) \rightarrow \Dcal$ is a cartesian fibration.
\item The map $p$ is a two-sided fibration if and only if it is a lax two-sided fibration and the induced functor $\Ccal \rightarrow \Cat_{/\Dcal}$ factors through $\Cat^{\cart}_{/\Dcal}$. 
\item Assume that $p$ is a two-sided fibration and let $q: \Ecal' \rightarrow \Ccal \times \Dcal$ be another two-sided fibration. Then a functor $F: \Ecal \rightarrow \Ecal'$ equipped with an identification $qF = p$ is a morphism of two-sided fibrations if and only if it is a morphism of cocartesian fibrations over $\Ccal$, and for every object $c$ in $\Ccal$ the induced functor $p^{-1}(c) \rightarrow q^{-1}(c)$ is a morphism of cartesian fibrations over $\Dcal$.
\end{enumerate}
\end{proposition}
\begin{proof}
This follows from a combination of remarks \ref{remark lax bicartesian arrows} and \ref{remark bicartesian arrows}, together with (the dual of) \cite{HTT} corollary 4.3.1.15.
\end{proof}

Our next goal is to present a two-sided analog of the Grothendieck construction. We first recall the notion of weighted colimits of categories.

\begin{definition}
Let $\Bcal$ be a category and denote by $\Tw(\Bcal)$ the twisted arrow category of $\Bcal$.
Let $H: \Bcal \rightarrow \Cat$  and $W: \Bcal^\op \rightarrow \Cat$ be functors. The colimit of $H$ weighted by $W$ is defined to be the colimit of the composite functor
\[
\Tw(\Bcal) \rightarrow \Bcal^\op \times \Bcal \xrightarrow{W(-) \times H(-)} \Cat.
\]
The lax colimit of $H$ is the colimit of $H$ weighted by the functor $\Bcal_{-/}: \Bcal^\op \rightarrow \Cat$ obtained by straightening of the source fibration $\Funct([1], \Bcal) \rightarrow \Bcal$. The oplax colimit of $H$ is the colimit of $H$ weighted by the functor $(\Bcal_{-/})^\op = (\Bcal^\op)_{/-}$.
\end{definition}

It is proven in \cite{GHN} that the lax colimit of a functor $H: \Bcal \rightarrow \Cat$ recovers the total category of the cocartesian fibration associated to $H$. Similarly, the oplax colimit of $H$ recovers the total category of the cartesian fibration associated to $H$. We now define our two-sided Grothendieck construction to be a mixed lax-oplax colimit.

\begin{definition}
Let $\Ccal, \Dcal$ be categories and $H: \Ccal \times \Dcal^\op \rightarrow \Cat$ be a functor. We define the bilax colimit of $H$ to be the colimit of $H$ weighted by the functor $\Ccal_{-/} \times \Dcal_{/-} : \Ccal^\op \times \Dcal \rightarrow \Cat$. We denote it by $\int_{\Ccal \times \Dcal^\op} H$.
\end{definition}

\begin{example}
Let $\Ccal$ be a category and $H: \Ccal = \Ccal \times [0]^\op \rightarrow \Cat$ be a functor. Then we have
\[
\int_{\Ccal \times [0]^\op} \Ccal = \colim_{\Tw(\Ccal \times [0]^\op)} \Ccal_{-/} \times [0]_{/-} \times H = \colim_{\Tw(\Ccal)} \Ccal_{-/} \times H
\]
Using \cite{GHN} theorem 7.4 we see that $\int_{\Ccal \times [0]^\op} H$ is the total category of the cocartesian fibration associated to $H$.
\end{example}
\begin{example}
Let $\Dcal$ be a category and $H: \Dcal^\op = [0] \times \Dcal^\op \rightarrow \Cat$ be a functor. Then by \cite{GHN} corollary 7.6 we have that $\int_{[0] \times \Dcal^\op} H$ is the total category of the cartesian fibration associated to $H$.
\end{example}

\begin{example}\label{example integral terminal}
Let $\Ccal, \Dcal$ be categories and $H: \Ccal \times \Dcal^\op \rightarrow \Cat$ be the functor which is constant $[0]$. We have
\[
\int_{\Ccal \times \Dcal^\op} H = \colim_{\Tw(\Ccal \times \Dcal^\op)} \Ccal_{-/} \times \Dcal_{/-} = \colim_{\Tw(\Ccal)}  \Ccal_{-/}  \times  \colim_{\Tw(\Dcal^\op)}  \Dcal_{/-}
\]
Using \cite{GHN} corollary 7.5 we conclude that $\int_{\Ccal \times \Dcal^\op} H = \Ccal \times \Dcal^\op$.
\end{example}

\begin{remark}
Let $\Ccal, \Dcal$ be categories. The assignment $H \mapsto \int_{\Ccal \times \Dcal^\op} H$ can be enhanced to a functor $\Funct(\Ccal \times \Dcal^\op, \Cat) \rightarrow \Cat$. The functor constant $[0]$ is the terminal object in $\Funct(\Ccal \times \Dcal^\op, \Cat)$. It follows from example \ref{example integral terminal} that the bilax colimit functor can be enhanced to a functor
\[
\int_{\Ccal \times \Dcal^\op} : \Funct(\Ccal \times \Dcal^\op, \Cat) \rightarrow \Cat_{/\Ccal \times \Dcal^\op}.
\]
\end{remark}

\begin{proposition}\label{prop bicart} Let $\Ccal, \Dcal$ be categories. 
The functor $\int_{\Ccal \times \Dcal^\op}$ factors through $\Cat^\bicart_{/\Ccal, \Dcal}$, and induces an equivalence 
\[\Funct(\Ccal \times \Dcal^\op, \Cat) = \Cat^\bicart_{/\Ccal, \Dcal}.
\]
\end{proposition}
\begin{proof}
Let $H: \Ccal \times \Dcal^\op \rightarrow \Cat$ be a functor. We have
\begin{align*}
\int_{\Ccal \times \Dcal^\op} H & = \colim_{\Tw(\Ccal \times \Dcal^\op)} \Ccal_{-/} \times \Dcal_{/-} \times H \\ & = \colim_{\Tw(\Ccal)} (\colim_{\Tw(\Dcal)} \Dcal_{/-} \times H) \times \Ccal_{-/} 
\\ & = \int_{\Ccal \times [0]^\op} \left( \int_{[0] \times\Dcal^\op} H \right).
\end{align*}
It follows that the functor $\int_{\Ccal \times \Dcal^\op} : \Funct(\Ccal \times \Dcal^\op, \Cat) \rightarrow \Cat_{/\Ccal \times \Dcal^\op}$ is equivalent to the composite functor
\[
\Funct(\Ccal \times \Dcal^\op, \Cat) = \Funct(\Ccal, \Funct(\Dcal^\op, \Cat)) \rightarrow \Funct(\Ccal, \Cat_{/\Dcal}) \rightarrow \Cat_{/\Ccal \times \Dcal}
\]
where the middle arrow is the cartesian Grothendieck construction, and the last arrow is the cocartesian Grothendieck construction. This is an embedding, and it follows from proposition \ref{prop bicart 0} that its image coincides with $\Cat_{/\Ccal \times \Dcal}^\bicart$.
\end{proof}

\begin{remark}\label{remark groth constr}
The proof of proposition \ref{prop bicart} shows that the two-sided Grothendieck construction can be computed in two steps, by first doing cartesian unstraightening along $\Dcal$, and then cocartesian unstraightening along $\Ccal$. Alternatively, we can also compute it by first doing cocartesian unstraightening along $\Ccal$ and then cartesian straightening along $\Dcal$.
\end{remark}

\subsection{The arrow bifibration} \label{subsection arrow}

The notion of two-sided fibration specializes in the case when the fibers are groupoids to the notion of bifibration introduced in \cite{HTT}.

\begin{definition}
Let $\Ccal, \Dcal$ be categories. A two-sided fibration $p = (p_1, p_2): \Ecal \rightarrow \Ccal \times \Dcal$ is said to be a bifibration if  for every pair $(c, d)$ in $\Ccal \times \Dcal$ the fiber $p^{-1}((c, d))$ is a groupoid.
\end{definition}
 
\begin{remark}\label{remark lax bicart is bifibr}
If $p: \Ecal \rightarrow \Ccal \times \Dcal$ is a functor whose fibers are groupoids, then $p$ is a two-sided fibration if and only if it is a lax two-sided fibration. Moreover, if $p$ is a bifibration and $q : \Ecal' \rightarrow \Ccal \times \Dcal$ is another bifibration, then any functor $F: \Ecal \rightarrow \Ecal'$ equipped with an identification $qF = p$ is automatically a morphism of two-sided fibrations.
\end{remark}

\begin{remark}
Let $\Ccal, \Dcal$ be categories. 
Under the equivalence of proposition \ref{prop bicart}, the full subcategory of $\Cat^\bicart_{/\Ccal \times \Dcal}$ on the bifibrations gets identified with $\Funct(\Ccal \times \Dcal^\op, \Spc)$.
\end{remark}

For each category $\Ccal$, the arrow category of $\Ccal$ equipped with its target-source projection turns out to be a bifibration. As the  following proposition shows, it in fact plays a similar role in the theory of bifibrations as the twisted arrow category does in the theory of pairings of categories (see \cite{HA} section 5.2.1).

\begin{proposition}\label{prop arrow cat}
Let $\Ccal$ be a category. Then the bifibration associated to the functor $\Hom_\Ccal: \Ccal \times \Ccal^\op \rightarrow \Spc$ is equivalent to the projection $p = (\ev_1, \ev_0) : \Funct([1], \Ccal) \rightarrow \Ccal$.
\end{proposition}
\begin{proof}
The fact that $p$ is a bifibration follows directly from the criteria of remark  \ref{remark lax bicartesian arrows} together with remark \ref{remark lax bicart is bifibr}. Let $\ev_1^\vee: \Funct([1], \Ccal)^\vee \rightarrow \Ccal^\op$ be the cartesian fibration classified by the same functor as the cocartesian fibration $\ev_1$. This comes equipped with a projection $p^\vee = (\ev_1^\vee, \ev_0^\vee): \Funct([1], \Ccal)^\vee \rightarrow \Ccal^\op \times \Ccal$, which is a right fibration classified by the same functor that classifies the bifibration $p$. According to \cite{HA} proposition 5.2.1.11, the right fibration associated to $\Hom_\Ccal: \Ccal \times \Ccal^\op \rightarrow \Spc$ is equivalent to the canonical projection $\lambda = (t, s): \Tw(\Ccal) \rightarrow \Ccal^\op \times \Ccal$. Hence it suffices to show that $p^\vee$ is equivalent to $\lambda$.

We will use the description of dual fibrations from \cite{BGN} (translated to the language of simplicial spaces rather than simplicial sets). The space of $n$-simplices of the category $\Funct([1], \Ccal)^\vee$ is the  space of functors $\sigma: \Tw([n])^\op \rightarrow \Funct([1], \Ccal)$ such that 
for each $0 \leq i < j \leq n$ the cospan
\[
\begin{tikzcd}
& \sigma(i \rightarrow j) &  \\ \sigma(i \rightarrow j-1)\arrow{ur} & & \arrow{ul} \sigma(i+1 \rightarrow j)
\end{tikzcd}
\]
has a cocartesian right leg, and the image of its left leg under $\ev_1$ is an isomorphism. Under the isomorphism $\Hom_{\Cat}(\Tw([n])^\op,\Funct([1], \Ccal)) = \Hom_{\Cat}(\Tw([n])^\op \times [1], \Ccal)$ we have that the space of $n$-simplices in $\Funct([1], \Ccal)^\vee$ is the space of maps $\tau: \Tw([n])^\op \times [1] \rightarrow \Ccal$ such that for every $0 \leq i < j \leq n$ the morphisms $\tau(i \rightarrow j-1, 1) \rightarrow \tau(i \rightarrow j, 1)$ and $\tau(i+1 \rightarrow j, 0) \rightarrow \tau(i \rightarrow j, 0)$ are isomorphisms. In other words, this is the space of maps $\tau: \Tw([n])^\op \times [1] \rightarrow \Ccal$ which factor through the localization of $\Tw([n])^\op \times [1]$ at the  collection $S_n$ of arrows of the form $(i \rightarrow j-1, 1) \rightarrow (i \rightarrow j, 1)$ and $(i+1 \rightarrow j,0) \rightarrow (i \rightarrow j, 0)$.

Consider for each $n \geq 0$ the functor $\Tw([n])^\op \times [1] \rightarrow [2n+1]$ that maps $(i \rightarrow j, 0)$ to $j$ and $(i \rightarrow j, 1)$ to $2n+1 - i$. This induces an isomorphism $S_n^{-1}(\Tw([n])^\op \times [1]) = [2n+1]$. This isomorphism is natural in $n$ - namely, we have an isomorphism of simplicial categories $[n] \mapsto S_n^{-1}(\Tw([n]^\op \times [1])$  and $n \mapsto [n]\star [n]^\op = [2n+1]$. The latter corepresents the twisted arrow category construction, hence we see that $\Funct([1], \Ccal)^\vee$ is equivalent to $\Tw(\Ccal)$.

The projection $\ev_1^\vee : \Funct([1], \Ccal)^\vee \rightarrow  \Ccal^\op$ sends an $n$-simplex $\sigma: \Tw([n])^\op \rightarrow \Funct([1], \Ccal)$ to the composition of $\ev_1 \sigma$ with the functor $[n]^\op \rightarrow \Tw([n])^\op$ given by the formula $i \mapsto (i \rightarrow n)$. The projection $\ev_0^\vee: \Funct([1], \Ccal)^\vee \rightarrow \Ccal$ sends an $n$-simplex $\sigma: \Tw([n])^\op \rightarrow \Funct([1],\Ccal)$ to the composition of $\ev_0\sigma$ with the functor $[n] \rightarrow \Tw([n])^\op$ given by the formula $j \mapsto (0 \rightarrow j)$. Under the isomorphism $\Funct([1],\Ccal)^\vee([n]) = \Ccal([2n+1]) = \Ccal([n] \star [n]^\op)$ these assignments amount to precomposing with the natural inclusions $[n]^\op \rightarrow [n]\star[n]^\op$ and $[n] \rightarrow [n]\star [n]^\op$. We conclude that under the isomorphism $\Funct([1], \Ccal)^\vee = \Tw(\Ccal)$, the projection $p^\vee = (\ev_1^\vee, \ev_0^\vee)$ agrees with the projection $\lambda$, as desired.
\end{proof}

The following proposition shows that the arrow category enjoys a universal property in the category of two-sided fibrations. Specializing it to the category of bifibrations, we are able to conclude that the arrow  bifibration is the free bifibration on the diagonal functor  $\Delta: \Ccal \rightarrow \Ccal \times \Ccal$.

\begin{proposition}\label{prop univer arrow}
Let $\Ccal$ be a category and let $p= (\ev_1, \ev_0) : \Funct([1], \Ccal) \rightarrow \Ccal \times \Ccal$. Denote by $\psi: \Ccal = \Funct([0] , \Ccal) \rightarrow \Funct([1], \Ccal)$ the functor given by precomposition with the projection $[1] \rightarrow [0]$ and by $\Delta = p\psi: \Ccal \rightarrow \Ccal \times \Ccal$ the diagonal map. Let $r = (r_1, r_2): \Ecal \rightarrow \Ccal \times \Ccal$ be a two-sided fibration. Then precomposition with $\psi$ induces an embedding
\[
\Hom_{\Cat^\bicart_{/\Ccal \times \Ccal}}(p, r) \rightarrow \Hom_{\Cat_{/\Ccal \times \Ccal}}(\Delta, r)
\]
which identifies $\Hom_{\Cat^\bicart_{/\Ccal \times \Ccal}}(p, r)$ with the space of maps $\Delta \rightarrow r$ which send arrows in $\Ccal$ to bicartesian arrows in $\Ecal$.
\end{proposition}
\begin{proof}
Recall from \cite{GHN} section 4 that the functor $\psi$ presents $\ev_1$ as the free cocartesian fibration on $\id_\Ccal$. It follows that precomposition with $\psi$ induces an equivalence 
\[
\Hom_{(\Cat^\cocart_{/\Ccal})_{/\Ccal \times \Ccal}}(p, r) =  \Hom_{\Cat_{/\Ccal \times \Ccal}}(\Delta, r).
\]
The space $\Hom_{\Cat^\bicart_{/\Ccal \times \Ccal}}(p, r)$ is the subspace of $\Hom_{(\Cat^\cocart_{/\Ccal})_{/\Ccal \times \Ccal}}(p, r)$ containing those maps $F: p \rightarrow r$ which map $\ev_0$-cartesian arrows to $r_2$-cartesian arrows.

Let $F$ be an object in $\Hom_{(\Cat^\cocart_{/\Ccal})_{/\Ccal \times \Ccal}}(p, r) $.  We have to show  that $F$ maps $\ev_0$-cartesian arrows to $r_2$-cartesian arrows if and only if $F\psi$ maps arrows in $\Ccal$ to bicartesian arrows in $\Ecal$.

Assume first that $F\psi$ maps arrows in $\Ccal$ to bicartesian arrows in $r$. Let $f: \sigma \rightarrow \sigma'$ be an $\ev_0$-cartesian arrow in $\Funct([1],\Ccal)$. In other words, we have $\ev_1 f$ invertible. We have a commutative diagram
\[
\begin{tikzcd}
\psi \sigma(0) \arrow[d] \arrow[r] & \sigma \arrow[d] \\
\psi \sigma'(0) \arrow[r]          & \sigma'         
\end{tikzcd}
\]
where the horizontal arrows are $\ev_1$-cocartesian. Applying the functor $F$ we obtain a commutative diagram
\[
\begin{tikzcd}
F\psi \sigma(0) \arrow[d] \arrow[r] & F\sigma \arrow[d] \\
F\psi \sigma'(0) \arrow[r]          & F\sigma'.         
\end{tikzcd}
\]
The horizontal arrows are $r_1$-cocartesian, and the left vertical arrow is bicartesian. Since the image of the right vertical arrow under $r_1$ is an isomorphism and $r$ is a two-sided fibration we conclude that the right vertical arrow is $r_2$-cartesian, as desired.

Assume now that $F$ maps  $\ev_0$-cartesian arrows to $r_2$-cartesian arrows. Since $p$ is a bifibration and $F$ also maps $\ev_1$-cocartesian arrows to $r_1$-cocartesian arrows, we conclude that $F$ maps all arrows in $\Funct([1], \Ccal)$ to bicartesian arrows. It follows that the same is true for $F\psi$, as desired.
\end{proof}
\begin{corollary}
Let $\Ccal$ be a category, and let $p, \psi$ be as in the statement of proposition \ref{prop univer arrow}. Let $r: \Ecal \rightarrow \Ccal \times \Ccal$ be a bifibration. Then precomposition with $\psi$ induces an equivalence
\[
\Hom_{\Cat^\bicart_{/\Ccal \times \Ccal}}(p, r) = \Hom_{\Cat_{/\Ccal \times \Ccal}}(\Delta, r).
\]
\end{corollary}

We now study the notion of representable bifibration, which is a two-sided analog of the notion of representable pairing of categories from \cite{HA}.
\begin{definition}
Let $\Ccal, \Dcal$ be categories.  We say that  a bifibration $p = (p_1, p_2): \Ecal \rightarrow \Ccal \times \Dcal$ is representable if for every object $c$ in $\Ccal$ the fiber $p_1^{-1}(c)$ has a final object.
\end{definition}

\begin{remark}
Let $\Ccal, \Dcal$ be categories. Recall that the full subcategory of $\Cat^\bicart_{/\Ccal \times \Dcal}$ on the bifibrations corresponds under the equivalence of proposition \ref{prop bicart} to the category $\Funct(\Ccal, \Pcal(\Dcal)) = \Funct(\Ccal \times \Dcal^\op, \Spc) \subseteq \Funct(\Ccal \times \Dcal^\op, \Cat)$. The full subcategory of $ \Cat^\bicart_{/\Ccal \times \Dcal}$ on the representable bifibrations agrees under this correspondence with the subcategory $\Funct(\Ccal , \Dcal) \subseteq \Funct(\Ccal, \Pcal(\Dcal))$.
\end{remark} 
 
\begin{notation}
Let $\Ccal, \Dcal$ be categories and $F: \Ccal \rightarrow \Dcal$ be a functor. We define the mapping cocylinder of $F$ as the following pullback:
\[
\begin{tikzcd}
\Cocyl(F) \arrow{r}{\overline{F}} \arrow{d}{\overline{\ev_1}} & \Funct([1], \Dcal) \arrow{d}{\ev_1} \\ \Ccal \arrow{r}{F} & \Dcal.
\end{tikzcd}
\]
Note that $\Cocyl(F)$ comes equipped with a projection $p^F= (p^F_1, p^F_2) : \Cocyl(F) \rightarrow \Ccal \times \Dcal$ given by $p^F_1 = \overline{\ev_1}$ and $p^F_2 = \ev_0 \overline{F}$.
\end{notation}

\begin{proposition}\label{prop repr bifibr}
Let $\Ccal, \Dcal$ be categories and $F: \Ccal \rightarrow \Dcal$ be a functor. Then the projection $p^F: \Cocyl(F) \rightarrow \Ccal \times \Dcal$ is a representable bifibration, and the induced functor $\Ccal \rightarrow \Dcal$ is equivalent to $F$.
\end{proposition}
\begin{proof}
The fact that $p^F$ is a bifibration follows directly from the fact that it is the base change of the bifibration $p: \Funct([1], \Dcal) \rightarrow \Dcal \times \Dcal$ along the functor $(F, \id_\Dcal): \Ccal \times \Dcal \rightarrow \Dcal \times \Dcal$. The functor $\Ccal \rightarrow \Cat^\cart_{/\Dcal}$ classifying $p^F_1$ is the composition
\[
\Ccal \xrightarrow{F} \Dcal \rightarrow \Cat^\cart_{/\Dcal}
\]
where the second arrow is the functor classifying the projection $\ev_1: \Funct([1], \Dcal) \rightarrow \Dcal$. By proposition \ref{prop arrow cat}, this corresponds under the equivalence $\Cat^\cart_{\Dcal} = \Funct(\Dcal^\op, \Cat)$  to the composition of the Yoneda embedding $\Dcal \rightarrow \Pcal(\Dcal)$ and the inclusion $\Pcal(\Dcal) \rightarrow \Funct(\Dcal, \Cat)$. Hence we see that the functor $\Ccal \rightarrow \Pcal(\Dcal)$ associated to $p^F$ is the composition of $F$ with the Yoneda embedding of $\Dcal$. This means that $p^F$ is representable and the associated functor $\Ccal \rightarrow \Dcal$ is equivalent to $F$, as desired.
\end{proof}

\subsection{Bivariant fibrations}\label{subsection bivar}
We now introduce the class of bivariant fibrations.

\begin{definition}
Let $\Ccal, \Dcal, \Ecal$ be categories, and $p: \Ecal \rightarrow \Ccal \times \Dcal$ be a functor. We say that $p$ is a lax bivariant fibration if it is both a cocartesian and a cartesian fibration.
\end{definition}

\begin{remark}
Let $\Ccal, \Dcal$ be categories, and $p = (p_1, p_2) : \Ecal \rightarrow \Ccal \times \Dcal$ be a lax bivariant fibration. Then $p$ is a lax two-sided fibration from $\Ccal$ to $\Dcal$, and a lax two-sided fibration from $\Dcal$ to $\Ccal$.
\end{remark}

\begin{definition}
Let $\Ccal, \Dcal$ be categories, and $p = (p_1, p_2): \Ecal \rightarrow \Ccal \times \Dcal$ be a lax bivariant fibration. We say that $p$ is a bivariant fibration if $p$ is a two-sided fibration from $\Ccal$ to $\Dcal$ and also from $\Dcal$ to $\Ccal$.
\end{definition}

\begin{definition}
Let $\Ccal, \Dcal$ be categories and $p: \Ecal \rightarrow \Ccal \times \Dcal$ and $p': \Ecal' \rightarrow \Ccal \times \Dcal$ be bivariant fibrations. A functor $F: \Ecal \rightarrow \Ecal'$ equipped with an identification $p'F = p$ is said to be a morphism of bivariant fibrations if it is a morphism of cartesian and cocartesian fibrations over $\Ccal \times \Dcal$. We denote by $\Cat^\bivar_{/\Ccal, \Dcal}$ the subcategory of $\Cat_{/\Ccal \times \Dcal}$ on the bivariant fibrations and morphisms of bivariant fibrations.
\end{definition}

Our next goal is to identify the image of $\Cat^\bivar_{/\Ccal, \Dcal}$ across the different versions of the straightening equivalence.  We first review the notion of adjointable diagram of categories.

\begin{definition}
We say that a commutative diagram of categories
\[
\begin{tikzcd}
\Ecal_{00} \arrow{d}{\beta'} \arrow{r}{\alpha'} & \Ecal_{10} \arrow{d}{\beta} \\ \Ecal_{01} \arrow{r}{\alpha} & \Ecal_{11}
\end{tikzcd}
\]
is vertically right adjointable if the following conditions hold:
\begin{itemize}
\item The functors $\beta$ and $\beta'$ have right adjoints $\beta^R$ and $\beta'^R$.
\item The natural transformation
\[
\alpha' (\beta')^R \rightarrow \beta^R \beta \alpha' (\beta')^R = \beta^R \alpha \beta' (\beta')^R \rightarrow \beta^R \alpha
\]
built from the unit $\id_{\Ecal_{10}} \rightarrow \beta^R \beta$ and the counit $\beta'(\beta')^R \rightarrow \id_{\Ecal_{01}} $, is a natural isomorphism.
\end{itemize}

We say that the above diagram is horizontally right adjointable if its transpose is vertically right adjointable. We say that it is right adjointable if it is both horizontally and vertically right adjointable. We say that it is (vertically / horizontally) left adjointable if the diagram obtained by taking opposites of all the categories involved is (vertically / horizontally) right adjointable.
\end{definition}

\begin{remark}\label{remark passage diagrams to adjoints}
Let
\[
\begin{tikzcd}
\Ecal_{00} \arrow{d}{\beta'} \arrow{r}{\alpha'} & \Ecal_{10} \arrow{d}{\beta} \\ \Ecal_{01} \arrow{r}{\alpha} & \Ecal_{11}
\end{tikzcd}
\]
be a vertically right adjointable commutative diagram of categories. Then the natural isomorphism $\alpha'(\beta')^R = \beta^R \alpha$ exhibits the following diagram as commutative:
\[
\begin{tikzcd}
\Ecal_{00}  \arrow{r}{\alpha'} & \Ecal_{10}  \\ \Ecal_{01} \arrow{u}{(\beta')^R} \arrow{r}{\alpha} & \Ecal_{11} \arrow{u}{\beta^R}.
\end{tikzcd}
\]
We say that this diagram arises from the first one by passage to right adjoints of the vertical arrows. We can similarly talk about diagrams obtained by passage to right adjoints of horizontal arrows, or left adjoints of vertical / horizontal arrows, by requiring that the original diagram satisfy the appropriate adjointability condition.
\end{remark}

The following proposition provides a link between the notion of adjointability of commutative squares of categories and the theory of two-sided fibrations.

\begin{proposition}\label{prop adjointability y bicart}
Let $\Ccal, \Dcal$ be categories and $H: \Ccal \times \Dcal \rightarrow \Cat$ be a functor. The following conditions are equivalent:
\begin{enumerate}[\normalfont (i)]
\item For every pair of arrows $\alpha: c \rightarrow c'$ in $\Ccal$ and $\beta: d \rightarrow d'$ in $\Dcal$ the commutative diagram of categories
\[
\begin{tikzcd}[column sep = large]
H(c, d) \arrow{r}{H(\alpha, \id_d)} \arrow{d}{H(\id_c, \beta)} & H(c', d) \arrow{d}{H(\id_{c'},\beta)} \\ H(c, d') \arrow{r}{H(\alpha, \id_{d'})} & H(c', d')
\end{tikzcd}
\]
is vertically right adjointable.
\item The cocartesian fibration $p^\cocart: \Ecal^\cocart \rightarrow \Ccal \times \Dcal$ classified by $H$ is a two-sided fibration.
\item The two-sided fibration $p^\bicart: \Ecal^\bicart \rightarrow \Dcal \times \Ccal^\op$ classified by $H$ is a cartesian fibration.
\item Let $G: \Ccal \rightarrow \Cat_{/\Dcal}$ be the image of $H$ under the composite map 
\[
\Funct(\Ccal \times \Dcal, \Cat) = \Funct(\Ccal, \Cat^\cocart_{/\Dcal}) \rightarrow \Funct(\Ccal, \Cat_{/\Dcal}).
\]
Then $G$ factors through the subcategory $\Cat^\cart_{/\Dcal}$.
\end{enumerate}
\end{proposition}
\begin{proof}
The equivalence between conditions (ii) and (iv) is a direct consequence of proposition \ref{prop bicart 0}. To show that conditions (iii) and (iv) are equivalent, note that the two-sided fibration $p^\bicart$ is obtained by applying the cartesian unstraightening construction to the functor $G$. The equivalence then follows from a combination of \cite{HTT} propositions 2.4.2.8 and 2.4.2.11.

We now establish the equivalence between conditions (i) and (ii). Let $c$ be an object in $\Ccal$, and $\beta: d \rightarrow d'$ be an arrow in $\Dcal$. The functor $H(\id_c, \beta)$ admits a right adjoint if and only if the arrow $(\id_c, \beta)$ in $\Ccal \times \Dcal$ admits a locally $p^\cocart$-cartesian lift. By (the dual of) \cite{HTT} corollary 4.3.1.15, this  happens if and only if $(\id_c, \beta)$ admits a $p^\cocart$-cartesian lift. This holds for all pairs $c, \beta$ if and only if $p^\cocart$ is a lax two-sided fibration.

Assume now that $p^\cocart$ is indeed a lax two-sided fibration. Let $\alpha: c \rightarrow c'$ and $\beta: d \rightarrow d'$ be a pair of arrows in $\Ccal$ and $\Dcal$ and let $e$ be an object in $(p^\cocart)^{-1}(c, d') = H(c, d')$. Consider the commutative diagram
\[
\begin{tikzcd}[column sep = large, row sep= large]
e' \arrow{r}{ \overline{(\alpha, \id_d)}} \arrow{d}{\widehat{(\id_c, \beta)}} & f \arrow{r}{\zeta} & f' \arrow{d}{\widehat{(\id_{c'},\beta)}} \\ e \arrow{rr}{\overline{(\alpha, \id_{d'})}} &  & e''
\end{tikzcd}
\]
where $\widehat{A}$ denotes a $p^\cocart$-cocartesian lift of the arrow $A$, and $\overline{A}$ denotes a $p^\cocart$-cartesian lift of $A$ . We have $f = H(\alpha, \id_d) H(\id_c,\beta)^R e$ and $f' = H(\id_{c'},\beta)^R H(\alpha, \id_{d'})e$. The projection $p^\cocart$ is a two-sided fibration if and only if for every choice of $\alpha, \beta$ and $e$, the resulting map $\zeta$ is an isomorphism.

We now enlarge the above diagram as follows:
\[
\begin{tikzcd}[column sep = large, row sep= large]
e' \arrow{dd}{\widehat{(\id_c, \beta)}} \arrow{rd}{\overline{(\id_c, \beta)}}\arrow{rr}{\overline{(\alpha, \id_d)}} &                         & f \arrow{rd}[swap]{\overline{(\id_{c'}, \beta)}} \arrow{r}{\zeta_1} & g \arrow{r}{\zeta_2}       \arrow{d}{\widehat{(\id_{c'},\beta)}}    & f'  \arrow{dd}{\widehat{(\id_{c'},\beta)}}            \\
                                    & x \arrow{ld}{\eta} \arrow{rr}[swap]{\overline{(\alpha, \id_{d'})}} &                        & y \arrow{rd}{\eta'}  &                \\
e \arrow{rrrr}{\overline{(\alpha, \id_{d'})}}                      &                         &                        &                        & e''. 
\end{tikzcd}
\]
The map $\zeta$ from the previous diagram is now decomposed as the composite map $f\xrightarrow{\zeta_1} g \xrightarrow{\zeta_2} f'$. We have equivalences
\[
H(\id_{c'},\beta)^R H(\id_{c'},\beta) H(\alpha, \id_d) H^R(\id_c,\beta)e = g = H(\id_{c'},\beta)^R H(\alpha, \id_{d'}) H(\id_c,\beta) H^R(\id_c,\beta)e.
\]
Under these equivalences, the map $\zeta_1$ is the unit of the adjunction $H(\id_{c'},\beta)\dashv H(\id_{c'},\beta)^R$ evaluated at $f$, and the map $\zeta_2$ is the image under $H(\id_{c'},\beta)^R H(\alpha, \id_{d'})$ of the counit of the adjunction $H(\id_c, \beta) \dashv H(\id_c, \beta)^R$ applied to $e$. We conclude that the map $\zeta$ is a component of the natural transformation witnessing the lax commutativity of the diagram obtained by passing to right adjoints of the vertical arrows of the square in the statement. It follows that the square in the statement is vertically right adjointable if and only  if for every choice of $\alpha, \beta$ and $e$ the resulting map $\zeta$ is an isomorphism, which we already observed is equivalent to $p^\cocart$ being a two-sided fibration.
\end{proof}

\begin{corollary}\label{coro equiv left adj}
Let $\Ccal, \Dcal$ be categories and $H: \Ccal \times \Dcal \rightarrow \Cat$ be a functor. The following conditions are equivalent:
\begin{enumerate}[\normalfont (i)]
\item For every pair of arrows $\alpha: c \rightarrow c'$ in $\Ccal$ and $\beta: d \rightarrow d'$ in $\Dcal$ the commutative diagram of categories
\[
\begin{tikzcd}[column sep = large]
H(c, d) \arrow{r}{H(\alpha, \id_d)} \arrow{d}{H(\id_c, \beta)} & H(c', d) \arrow{d}{H(\id_{c'},\beta)} \\ H(c, d') \arrow{r}{H(\alpha, \id_{d'})} & H(c', d')
\end{tikzcd}
\]
is vertically left adjointable.
\item The cartesian fibration $p^\cart: \Ecal^\cart \rightarrow \Dcal^\op \times \Ccal^\op$ classified by $H$ is a two-sided fibration.
\item The two-sided fibration $p^\bicart: \Ecal^\bicart \rightarrow \Ccal \times \Dcal^\op$ classified by $H$ is a cocartesian fibration.
\item Let $G: \Ccal \rightarrow \Cat_{/\Dcal^\op}$ be the image of $H$ under the composite map 
\[
\Funct(\Ccal \times \Dcal, \Cat) = \Funct(\Ccal, \Cat^\cart_{/\Dcal^\op}) \rightarrow \Funct(\Ccal, \Cat_{/\Dcal^\op}).
\]
Then $G$ factors through the subcategory $\Cat^\cocart_{/\Dcal^\op}$.
\end{enumerate}
\end{corollary}
\begin{proof}
This follows from proposition \ref{prop adjointability y bicart} applied to the functor $H^\op: \Ccal \times \Dcal \rightarrow \Cat$.
\end{proof}

\begin{definition}
Let $\Ccal, \Dcal$ be categories and $H: \Ccal \times \Dcal \rightarrow \Cat$ be a functor. We say that $H$ is right adjointable in the $\Dcal$ coordinate if the equivalent conditions of proposition \ref{prop adjointability y bicart} are satisfied. We say that $H$ is left adjointable in the $\Dcal$ coordinate if the equivalent conditions of corollary \ref{coro equiv left adj} are satisfied. By switching the role of $\Ccal$ and $\Dcal$ we can similarly talk about right/left adjointability in the $\Ccal$ coordinate.
\end{definition}

\begin{notation}
Let $\Ccal, \Dcal$ be categories and $H: \Ccal \times \Dcal \rightarrow \Cat$ be a functor. If $H$ is right adjointable in the $\Dcal$ coordinate then condition (iv) in proposition \ref{prop adjointability y bicart} provides a functor $\Ccal \rightarrow \Cat^\cart_{/\Dcal} = \Funct(\Dcal^\op, \Cat)$. We denote by $H^{R_\Dcal}: \Ccal \times \Dcal^\op \rightarrow \Cat$ the induced functor. We say that this is obtained from $H$ by passage to right adjoints in the $\Dcal$ coordinate. Similarly, if $H$ is left adjointable then out of equivalence (iv) in corollary \ref{coro equiv left adj} we obtain a functor $H^{L_\Dcal}: \Ccal \times \Dcal^\op \rightarrow \Cat$ which is said to arise from $H$ by passage to left adjoints in the $\Dcal$ coordinate.
\end{notation}

\begin{notation}
Let $\Ccal$ be a category and $H: \Ccal \rightarrow \Cat$ be a functor. Identifying $\Ccal$ with $[0] \times \Ccal$, it makes sense to talk about right and left adjointability of $H$ in the coordinate $\Ccal$. In this situation, we simply say that $H$ is right (resp. left) adjointable if it is right (resp. left) adjointable in the coordinate $\Ccal$. We will use the notation $H^R$ (resp. $H^L$) instead of $H^{R_\Ccal}$ (resp. $H^{L_\Ccal})$.
\end{notation}

\begin{remark}
Let $\Ccal, \Dcal$ be categories. Then passage to right and left adjoints in the $\Dcal$ coordinate define inverse equivalences between the space of functors $\Ccal \times \Dcal \rightarrow \Cat$ which are right adjointable in the $\Dcal$ coordinate, and the space of functors $\Ccal \times \Dcal^\op \rightarrow \Cat$ which are left adjointable in the $\Dcal$ coordinate.
\end{remark}

\begin{remark}\label{prop equiv passage right}
Let $\Ccal, \Dcal$ be categories and $H: \Ccal \times \Dcal \rightarrow \Cat$ be a functor.  Let $p^\cocart: \Ecal^\cocart \rightarrow \Ccal \times \Dcal$ be the cocartesian fibration classified by $H$ and $p^\bicart: \Ecal^\bicart \rightarrow \Dcal \times \Ccal^\op$ be the two-sided fibration classified by $H$. Assume that $H$ is right adjointable in the $\Dcal$ coordinate. Then the two-sided fibration $p^\cocart$ and the cartesian fibration $p^\bicart$ are both classified by the functor $H^{R_\Dcal}$. Moreover it follows from the proof of proposition \ref{prop adjointability y bicart} that for every pair of arrows $\alpha: c \rightarrow c'$ and $\beta: d \rightarrow d'$ in $\Ccal$ and $\Dcal$ respectively, the commutative diagram
\[
\begin{tikzcd}[column sep = huge]
H^{R_\Dcal}(c, d) \arrow{r}{H^{R_\Dcal}(\alpha, \id_d)} & H^{R_\Dcal}(c', d) \\ H^{R_\Dcal}(c, d') \arrow{u}{H^{R_\Dcal}(\id_c, \beta)}  \arrow{r}{H^{R_\Dcal}(\alpha, \id_{d'})} & H^{R_\Dcal}(c', d')  \arrow{u}{H^{R_\Dcal}(\id_{c'},\beta)}
\end{tikzcd}
\]
is equivalent to the diagram obtained from
\[
\begin{tikzcd}[column sep = large]
H(c, d) \arrow{r}{H(\alpha, \id_d)} \arrow{d}{H(\id_c, \beta)} & H(c', d) \arrow{d}{H(\id_{c'},\beta)} \\ H(c, d') \arrow{r}{H(\alpha, \id_{d'})} & H(c', d')
\end{tikzcd}
\]
by passage to right adjoints of vertical arrows (see remark \ref{remark passage diagrams to adjoints}).
\end{remark}

\begin{remark}\label{coro equiv passage left}
Let $\Ccal, \Dcal$ be categories and $H: \Ccal \times \Dcal \rightarrow \Cat$ be a functor.  Let $p^\cart: \Ecal^\cart \rightarrow \Dcal^\op \times \Ccal^\op$ be the cartesian fibration classified by $H$ and $p^\bicart: \Ecal^\bicart \rightarrow \Ccal \times \Dcal^\op$ be the two-sided fibration classified by $H$. Assume that $H$ is left adjointable in the $\Dcal$ coordinate. Then the two-sided fibration $p^\cart$ and the cocartesian fibration $p^\bicart$ are both classified by the functor $H^{L_\Dcal}$. Moreover, for every pair of arrows $\alpha: c \rightarrow c'$ and $\beta: d \rightarrow d'$ in $\Ccal$ and $\Dcal$ respectively, the commutative diagram
\[
\begin{tikzcd}[column sep = huge]
H^{L_\Dcal}(c, d) \arrow{r}{H^{L_\Dcal}(\alpha, \id_d)} & H^{L_\Dcal}(c', d) \\ H^{L_\Dcal}(c, d') \arrow{u}{H^{L_\Dcal}(\id_c, \beta)}  \arrow{r}{H^{L_\Dcal}(\alpha, \id_{d'})} & H^{L_\Dcal}(c', d')  \arrow{u}{H^{L_\Dcal}(\id_{c'},\beta)}
\end{tikzcd}
\]
is equivalent to the diagram obtained from
\[
\begin{tikzcd}[column sep = large]
H(c, d) \arrow{r}{H(\alpha, \id_d)} \arrow{d}{H(\id_c, \beta)} & H(c', d) \arrow{d}{H(\id_{c'},\beta)} \\ H(c, d') \arrow{r}{H(\alpha, \id_{d'})} & H(c', d')
\end{tikzcd}
\]
by passage to left adjoints of vertical arrows.
\end{remark}

The notion of bivariant fibration arises naturally when studying functors of two variables which have adjointability properties with respect to both of them. We consider first the case of functors which are either right or left adjointable with respect to both variables.

\begin{proposition}\label{prop equiv being bivariant cartesian}
Let $\Ccal, \Dcal$ be categories and $p : \Ecal \rightarrow \Ccal \times \Dcal$ be a cocartesian fibration classified by a functor $H: \Ccal \times \Dcal \rightarrow \Cat$. Then $p$ is a bivariant fibration if and only if $H$ is right adjointable in both the $\Ccal$ coordinate and the $\Dcal$ coordinate.
\end{proposition}
\begin{proof}
According to proposition \ref{prop adjointability y bicart}, the functor $H$ is right adjointable in both the $\Ccal$ and $\Dcal$ coordinates if and only if $p$ is a two-sided fibration from $\Ccal$ to $\Dcal$ and from $\Dcal$ to $\Ccal$. In this case, for every arrow $(\alpha, \beta): (c, d ) \rightarrow (c', d')$ in $\Ccal \times \Dcal$, the functor $H(\alpha, \beta)$ has a right adjoint, since it can be written as the composition of the right adjointable functors $H(\alpha, \id_d)$ and $H(\id_{c'}, \beta)$. We conclude that in this case $p$ is also a cartesian fibration, and therefore a bivariant fibration, as desired.
\end{proof}

\begin{corollary}\label{coro equiv being bivariant cocartesian}
Let $\Ccal, \Dcal$ be categories and $p : \Ecal \rightarrow \Ccal \times \Dcal$ be a cartesian fibration classified by a functor $H: \Ccal^\op \times \Dcal^\op \rightarrow \Cat$. Then $p$ is a bivariant fibration if and only if $H$ is left adjointable in both the $\Ccal$ coordinate and the $\Dcal$ coordinate.
\end{corollary}
\begin{proof}
Apply proposition \ref{prop equiv being bivariant cartesian} to the functor $H^\op: \Ccal \times \Dcal \rightarrow \Cat$.
\end{proof}

We now deal with the case of functors which enjoy mixed adjointability properties.

\begin{proposition}\label{prop equiv being bivariant bicart}
Let $\Ccal, \Dcal$ be categories and $p : \Ecal \rightarrow \Ccal \times \Dcal$ be a two-sided fibration classified by a functor $H: \Ccal \times \Dcal^\op \rightarrow \Cat$. Then the following are equivalent
\begin{enumerate}[\normalfont (i)]
\item The map $p$ is a bivariant fibration.
\item The functor $H$ is left adjointable in the $\Dcal^\op$ coordinate and the functor $H^{L_\Dcal}$ is right adjointable in the $\Ccal$ coordinate.
\item The functor $H$ is right adjointable in the $\Ccal$ coordinate and the functor $H^{R_\Ccal}$ is left adjointable in the $\Dcal^\op$ coordinate.
\end{enumerate}
\end{proposition}
\begin{proof}
By corollary \ref{coro equiv left adj}, $H$ is left adjointable in the $\Dcal^\op$ coordinate if and only if $p$ is a cocartesian fibration. In this case, $p$ is the cocartesian fibration classified by $H^{L_\Dcal}$. The equivalence between conditions (i) and (ii) is now a direct consequence of proposition \ref{prop equiv being bivariant cartesian}. Similarly, proposition \ref{prop adjointability y bicart} shows that $H$ is right adjointable in the $\Ccal$ coordinate if and only if $p$ is a cartesian fibration, and moreover in this case $p$ is the cartesian fibration classified by $H^{R_\Ccal}$. The equivalence between conditions (i) and (iii) now follows from corollary \ref{coro equiv being bivariant cocartesian}.
\end{proof}

\begin{remark}
Let $\Ccal, \Dcal$ be categories and $p : \Ecal \rightarrow \Ccal \times \Dcal$ be a bivariant fibration. Then the different versions of the straightening equivalence yield four functors
\begin{align*}
H_{\cocart}: \Ccal \times \Dcal \rightarrow \Cat & \hspace{2cm} H_{\cocart,\cart}:\Ccal \times \Dcal^\op \rightarrow \Cat \\ H_{\cart, \cocart}: \Ccal^\op \times \Dcal \rightarrow \Cat & \hspace{2cm} H_{\cart}: \Ccal^\op \times \Dcal^\op \rightarrow \Cat.
\end{align*}
It follows from remarks \ref{prop equiv passage right} and \ref{coro equiv passage left} that each of these four functors determines the rest, by passage to right or left adjoints in the appropriate coordinates. In particular, the operations of passing to adjoints in different coordinates commute.
\end{remark}

We now examine the image of the class of morphisms of bivariant fibrations under the different versions of the straightening equivalence.

\begin{remark}\label{remark adjointability triple product}
Let $\Ical, \Ccal, \Dcal$ be categories, and $H: \Ical \times \Ccal \times \Dcal \rightarrow \Cat$ be a functor.  Identifying $\Ical \times \Ccal \times \Dcal$ with $(\Ical \times \Ccal) \times \Dcal$, we may talk about right adjointability of $H$ with respect to the $\Dcal$ coordinate. Since any map in $\Ical \times \Ccal$ is a composition of maps which are constant in one coordinate, the functor $H$ is right adjointable in the $\Dcal$ coordinate if and only if the following two conditions hold:
\begin{itemize}
\item For every object $i$ in $\Ical$ the induced functor $H(i, -, -): \Ccal \times \Dcal \rightarrow \Cat$ is right adjointable in the $\Dcal$ coordinate.
\item For every object $c$ in $\Ccal$ the induced functor $H(-, c, -) : \Ical \times \Dcal \rightarrow \Cat$ is right adjointable in the $\Dcal$ coordinate.
\end{itemize}
The above can be adapted in a straightforward way to yield characterizations for left adjointability in the $\Dcal$ coordinate, or right/left adjointability in the $\Ccal$ coordinate.
\end{remark}

\begin{proposition}\label{prop characterization morph of bivariant}
Let $\Ccal, \Dcal$ be categories. Let $p : \Ecal \rightarrow \Ccal \times \Dcal$ and $p': \Ecal \rightarrow \Ccal \times \Dcal$ be bivariant fibrations, and let $F: p \rightarrow p'$ be a morphism in $\Cat_{/\Ccal \times \Dcal}$. The following are equivalent: 
\begin{enumerate}[\normalfont (i)]
\item The map $F$ is a morphism of bivariant fibrations.
\item The map $F$ is a morphism of cocartesian fibrations, and the induced functor $[1] \times \Ccal \times \Dcal \rightarrow \Cat$ is right adjointable in the $\Ccal$ and $\Dcal$ coordinates.
\item The map $F$ is a morphism of two-sided fibrations from $\Ccal$ to $\Dcal$, and the induced functor $[1] \times \Ccal \times \Dcal^\op \rightarrow \Cat$ is right adjointable in the $\Ccal$ coordinate, and left adjointable in the $\Dcal^\op$ coordinate.
\item The map $F$ is a morphism of cartesian fibrations, and the induced functor $[1] \times \Ccal^\op \times \Dcal^\op \rightarrow \Cat$ is left adjointable in the $\Ccal$ and $\Dcal$ coordinates.
\item The map $F$ is a morphism of two-sided fibrations from $\Dcal$ to $\Ccal$, and the induced functor $[1] \times \Dcal \times \Ccal^\op \rightarrow \Cat$ is right adjointable in the $\Dcal$ coordinate, and left adjointable in the $\Ccal^\op$ coordinate.
\end{enumerate}
\end{proposition}
\begin{proof}
We show that conditions (i) and (ii) are equivalent - the equivalence between (i) and each of the items (iii)-(v) follows along similar lines. Assume that $F$ is a morphism of cocartesian fibrations and denote by $G: [1] \times \Ccal \times \Dcal \rightarrow \Cat$ the induced functor. Since $p$ is a bivariant fibration, the functors $\Ccal \times \Dcal \rightarrow \Cat$ classifying $p$ and $p'$  are right adjointable in the $\Ccal$ and $\Dcal$  coordinates. By virtue of remark \ref{remark adjointability triple product}, the functor $G$ is right adjointable in the $\Ccal$ and $\Dcal$ coordinates if and only if the following conditions are satisfied:
\begin{itemize}
\item For every object $c$ in $\Ccal$ the induced functor $G(-, c, -): [1] \times \Dcal \rightarrow \Cat$ is right adjointable in the $\Dcal$ coordinate.
\item For every object $d$ in $\Dcal$ the induced functor $G(-, - , d) : [1] \times \Ccal \rightarrow \Cat$ is right adjointable in the $\Ccal$ coordinate.
\end{itemize}
The result now follows from the characterization of adjointability given by item (iv) in proposition \ref{prop adjointability y bicart}.
\end{proof}

\subsection{The Beck-Chevalley condition}\label{subsection Beck Chev}
We now specialize to the class of bivariant fibrations satisfying an extra base change property.

\begin{definition}
Let $\Bcal$ be a category admitting pullbacks, and $p: \Ecal \rightarrow \Bcal$ be a functor which is both a cocartesian and a cartesian fibration. We say that  $p$ satisfies the Beck-Chevalley condition if for every cartesian square $C: [1] \times [1] \rightarrow \Bcal$ the base change of $p$ along $C$ is a bivariant fibration.
\end{definition}

\begin{remark}
Let $\Ccal, \Dcal$ be categories admitting pullbacks. Let $\alpha: c \rightarrow c'$ and $\beta: d \rightarrow d'$ be arrows in $\Ccal$ and $\Dcal$ respectively. Then the commutative diagram
\[
\begin{tikzcd}
(c, d) \arrow{r}{(\alpha, \id_d)} \arrow{d}{(\id_c, \beta)} & (c', d) \arrow{d}{(\id_{c'}, \beta)} \\ (c, d') \arrow{r}{(\alpha, \id_{d'})} & (c', d')
\end{tikzcd}
\]
is a cartesian square in $\Ccal \times \Dcal$. It follows that if $p: \Ecal \rightarrow \Ccal \times \Dcal$ is a lax bivariant fibration satisfying the Beck-Chevalley condition then $p$ is a bivariant fibration.
\end{remark}

\begin{proposition}\label{prop beck chev chequear algunas solas}
Let $\Ccal, \Dcal$ be categories admitting pullbacks, and $p: \Ecal \rightarrow \Ccal \times \Dcal$ be a bivariant fibration. Then $p$ satisfies the Beck-Chevalley condition if and only if the following conditions are satisfied:
\begin{itemize}
\item For every cartesian square $C: [1] \times [1] \rightarrow \Ccal$ and every object $d$ in $\Dcal$ the base change of $p$ along $C \times d$ is a bivariant fibration.
\item For every cartesian square $C: [1] \times [1] \rightarrow \Dcal$ and every object $c$ in $\Ccal$ the base change of $p$ along $c \times C$ is a bivariant fibration.
\end{itemize}
\end{proposition}
\begin{proof}
The pair of conditions given in the statement are evidently implied by the Beck-Chevalley condition. Assume now that $p$ satisfies the two conditions in the statement, and consider  a general cartesian square $C: [1] \times [1] \rightarrow \Ccal \times \Dcal$ depicted as follows:
\[
\begin{tikzcd}
(c' \times_c c'', d' \times_d d'') \arrow{r}{} \arrow{d}{} & (c'', d'') \arrow{d}{} \\ (c', d') \arrow{r}{} & (c,d).
\end{tikzcd}
\]
We can see $C$ as the outer boundary of a commutative diagram $G: [2]\times [2] \rightarrow \Ccal \times \Dcal$, as follows:
\[
\begin{tikzcd}
(c' \times_c c'', d' \times_d d'') \arrow{r}{} \arrow{d}{} & (c'', d' \times_d d'') \arrow{d}{} \arrow{r}{} & (c'', d'') \arrow{d}{} \\ (c', d' \times_d d'') \arrow{d}{} \arrow{r}{} & (c,d' \times_d d'') \arrow{r}{} \arrow{d}{} & (c, d'') \arrow{d}{} \\ (c', d') \arrow{r}{}& (c, d') \arrow{r}{} &  (c,d).
\end{tikzcd}
\]
To show that the base change $C^*p$ is a  bivariant fibration it suffices to show that $G^*p$ is a bivariant fibration. To see this we must show that the base change of $p$ along each of the four small commutative diagrams inside $G$ is a bivariant fibration. The fact that this holds for the lower left and the upper right squares is a direct consequence of the fact that $p$ itself is a bivariant fibration. The fact that this holds for the upper left and lower right squares follows from our assumptions on $p$.
\end{proof}

\begin{corollary}\label{coro characterization bc}
Let $\Ccal, \Dcal$ be categories admitting pullbacks, and $p: \Ecal \rightarrow \Ccal \times \Dcal$ be a bivariant fibration classified by a functor $H: \Ccal \times \Dcal^\op \rightarrow \Cat$. Then $p$ satisfies the Beck-Chevalley condition if and only if the following conditions are satisfied:
\begin{itemize}
\item For every pair of maps $c' \rightarrow c \leftarrow c''$ in $\Ccal$ and every object $d$ in $\Dcal$
the commutative diagram of categories
\[
\begin{tikzcd}
H(c' \times_c c'', d) \arrow{d}{} \arrow{r}{} & H(c'',d) \arrow{d}{} \\
H(c',d) \arrow{r}{} & H(c,d)
\end{tikzcd}
\]
is right adjointable.
\item For every object $c$ in $\Ccal$ and every pair of arrows $d' \rightarrow d \leftarrow d''$ in $\Dcal$ the commutative diagram of categories
\[
\begin{tikzcd}
H(c, d) \arrow{d}{} \arrow{r}{} & H(c,d') \arrow{d}{} \\
H(c,d'') \arrow{r}{} & H(c,d' \times_d d'')
\end{tikzcd}
\]
is left adjointable.
\end{itemize}
\end{corollary}
\begin{proof}
Combine propositions \ref{prop beck chev chequear algunas solas} and \ref{prop adjointability y bicart}.
\end{proof}

We now study a fundamental example of a bivariant fibration satisfying the Beck-Chevalley condition, which plays a role analogous to that of the arrow category in the theory of two-sided fibrations.

\begin{notation}
Denote by $\Lambda^2_0$ the category with objects $0, 1, 2$ and nontrivial arrows $0 \rightarrow 1$ and $0 \rightarrow 2$. For each category $\Ccal$ and $i$ in $\Lambda^2_0$ we let $\ev_i: \Funct(\Lambda^2_0, \Ccal) \rightarrow \Ccal$ be the functor of evaluation at $i$. 
\end{notation}

\begin{proposition}\label{proof q is bivariant}
Let $\Ccal$ be a category admitting pullbacks. Then 
\[
q = (\ev_1, \ev_2): \Funct(\Lambda^2_0, \Ccal) \rightarrow \Ccal \times \Ccal
\]
 is a bivariant fibration which satisfies the Beck-Chevalley condition.
\end{proposition}
\begin{proof}
Let $\rho: \mu \rightarrow \mu'$ be an arrow in $\Funct(\Lambda^2_0, \Ccal)$. Then
\begin{itemize}
\item $\rho$ is $q$-cocartesian if and only if the induced map $\mu(0) \rightarrow \mu'(0)$ is an isomorphism.
\item $\rho$ is $q$-cartesian if and only if the diagram
\[
\begin{tikzcd}
                 & \mu(0) \arrow[ld] \arrow[rd] \arrow[d] &                  \\
\mu(1) \arrow[d] & \mu'(0) \arrow[ld] \arrow[rd]          & \mu(2) \arrow[d] \\
\mu'(1)          &                                        & \mu'(2)         
\end{tikzcd}
\]
is a limit diagram.
\end{itemize}
We thus see that $q$ has all cocartesian lifts, and that the existence of cartesian lifts is guaranteed by the fact that $\Ccal$ has pullbacks. Therefore $q$ is a lax bivariant fibration. 

It follows from the above characterization of cartesian and cocartesian arrows that a morphism $\mu \rightarrow \mu'$ in $\Funct(\Lambda^2_0, \Ccal)$ is bicartesian (from $\Ccal$ to $\Dcal$) if and only if the square
\[
\begin{tikzcd}
\mu(0) \arrow[d] \arrow[r] & \mu'(0) \arrow[d] \\
\mu(2) \arrow[r]           & \mu'(2)          
\end{tikzcd}
\]
is cartesian. We thus see that the class of bicartesian arrows is stable under composition, and therefore $q$ is a two-sided fibration from $\Ccal$ to $\Dcal$. Switching the roles of $\Ccal$ and $\Dcal$ we conclude that $q$ is also a two-sided fibration from $\Dcal$ to $\Ccal$, and therefore $q$ is a bivariant fibration.

It remains to show that $q$ satisfies the Beck-Chevalley condition. Consider a commutative diagram $C : [1] \times [1] \rightarrow \Funct(\Lambda^2_0, \Ccal)$ depicted as follows
\[
\begin{tikzcd}
\nu \arrow{r}{} \arrow{d}{} & \mu'' \arrow{d}{} \\ \mu' \arrow{r}{} & \mu
\end{tikzcd}
\]
and such that $qC$ is cartesian. Assume that the horizontal arrows are $q$-cocartesian and that the left vertical arrow is $q$-cartesian. We have to show that the right vertical arrow is $q$-cartesian. By virtue of proposition \ref{prop beck chev chequear algunas solas}, and since the bivariant fibrations $q$ and $(\ev_2, \ev_1)$ are equivalent, it suffices to consider the case when the image of the above square under $\ev_2$ is constant. 

Consider the commutative diagram
\[
\begin{tikzcd}
\nu(0) \arrow{d}{} \arrow{r}{} & \nu(1) \arrow{d}{} \arrow{r}{} & \mu''(1) \arrow{d}{} \\
\mu'(0) \arrow{r}{} & \mu'(1) \arrow{r}{} & \mu(1).
\end{tikzcd}
\]
The right inner square is $\ev_1(C)$ which is cartesian since $qC$ is cartesian. The left inner square is cartesian since the morphism $\nu \rightarrow \mu'$ is $\ev_1$-cartesian. Hence the outmost square is cartesian. Since the maps $\nu \rightarrow \mu''$ and $\mu' \rightarrow \mu$ are cocartesian, the outermost square is equivalent to
\[
\begin{tikzcd}
\mu''(0) \arrow{r}{} \arrow{d}{} & \mu''(1) \arrow{d}{} \\ \mu(0) \arrow{r}{} & \mu(1).
\end{tikzcd}
\]
Since the map $\mu'' \rightarrow \mu$ is constant under $\ev_2$ this shows that it is in fact $\ev_1$-cartesian, as desired.
\end{proof}

We finish by proving a universal property for $\Funct(\Lambda^2_0, \Ccal)$. Although it is possible to formulate a universal property for it in the category of bivariant fibrations, for our purposes we will need a version which treats it as a cocartesian and two-sided fibration.

\begin{notation}
Let $\Ccal, \Dcal$ be categories. We denote by $\Cat^{\cocart,\bicart}_{/\Ccal \times \Dcal}$ the intersection of $\Cat^{\cocart}_{/\Ccal \times \Dcal}$ and $\Cat^{\bicart}_{/\Ccal \times \Dcal}$ inside $\Cat_{/\Ccal \times \Dcal}$.
\end{notation}

\begin{proposition}\label{prop univer span}
Let $\Ccal$ be a category admitting pullbacks.  Let $p = (\ev_1, \ev_0): \Funct([1], \Ccal) \rightarrow \Ccal \times \Ccal$ and $q = (\ev_1, \ev_2) : \Funct(\Lambda^2_0, \Ccal) \rightarrow \Ccal \times \Ccal$. Let $\phi: p \rightarrow q$ be the map given by precomposition with the functor $\Lambda^2_0 \rightarrow [1]$ sending $0$ to $0$, $1$ to $1$, and $2$ to $0$.  Then 
\begin{enumerate}[\normalfont (i) ] 
\item The map $\phi$ is a morphism of two-sided fibrations. 
\item Let $r = (r_1, r_2): \Ecal \rightarrow \Ccal \times \Ccal$ be a cocartesian and two-sided fibration such that for every $c$ in $\Ccal$ the base change of $\Ecal$ along $c \times \id_\Ccal: \Ccal \rightarrow \Ccal \times \Ccal$ satisfies the Beck-Chevalley condition. Then precomposition with $\phi$ induces an equivalence
\[
\Hom_{\Cat^{\cocart,\bicart}_{/\Ccal \times \Ccal}}(q, r) = \Hom_{\Cat^{\bicart}_{/\Ccal \times \Ccal}}(p, r).
\]
\end{enumerate}
\end{proposition}

\begin{proof}
Let $g: \sigma \rightarrow \sigma'$ be an arrow in $\Funct([1], \Ccal)$. Then
\begin{itemize}
\item The arrow $g$ is $\ev_1$-cocartesian if and only if the induced map $\sigma(0) \rightarrow \sigma'(0)$ is an isomorphism.
\item The arrow $g$ is $\ev_0$-cartesian if and only if the induced map $\sigma(1) \rightarrow \sigma'(1)$ is an isomorphism.
\end{itemize}
It follows from this together with the description of $q$-cocartesian and $q$-cartesian arrows from the proof of proposition \ref{proof q is bivariant} that the map $\phi$ is a morphism of two-sided fibrations. 

It remains to check the universality of $\phi$. Let $\Delta: \Ccal \rightarrow \Ccal \times \Ccal$ be the diagonal map. Let $\psi: \Delta \rightarrow p$ be the map induced by precomposition with the projection $[1] \rightarrow [0]$.
It follows from \cite{GHN} section 4 that the map $\psi$ presents $\ev_1$ as the free cocartesian fibration on $\id_\Ccal$, and that $\phi\psi$ presents $q$ as the free cocartesian fibration on $\Delta$.  We therefore have equivalences
\[
\Hom_{\Cat^{\cocart}_{/\Ccal \times \Ccal}}(q, r) = \Hom_{\Cat_{/\Ccal \times \Ccal}}(\Delta, r) = \Hom_{(\Cat^\cocart_{/\Ccal})_{/\Ccal \times \Ccal}}(p, r).
\]
Under the above equivalence, the space $\Hom_{\Cat^\bicart_{/\Ccal \times \Ccal}}(p,r)$ becomes identified with the space of morphisms of cocartesian fibrations $F: q \rightarrow r$ whose composition with $\phi$ is a morphism of two-sided fibrations. Since $\phi$ is a morphism of two-sided fibrations, we see that precomposition with $\phi$ induces an inclusion 
\[
\Hom_{\Cat^{\cocart, \bicart}_{/\Ccal \times \Ccal}}(q, r)  \subseteq \Hom_{\Cat^{\bicart}_{/\Ccal \times \Ccal}}(p, r). 
\]

It remains to show that the above inclusion is an equivalence. Let $F: q \rightarrow r$ be a morphism of cocartesian fibrations and assume that $F\phi$ is a morphism of two-sided fibrations. We have to show that $F$ is a morphism of two-sided fibrations. Let $\rho: \mu \rightarrow \mu'$ be an $\ev_2$-cartesian arrow in $\Funct(\Lambda^2_0, \Ccal)$. Consider the commutative diagram $C$ in $\Ecal$ given as follows:
\[
\begin{tikzcd}
F\phi (\mu(1) \leftarrow \mu(0)) \arrow{r}{} \arrow{d}{} & F\mu \arrow{d}{} \\
F\phi (\mu'(1) \leftarrow \mu'(0)) \arrow{r}{} & F\mu'.
\end{tikzcd}
\]
Since the maps $\phi(\mu(1) \leftarrow \mu(0)) \rightarrow \mu$ and $\phi(\mu'(1) \leftarrow \mu'(0)) \rightarrow \mu'$ are $q$-cocartesian, we see that the horizontal arrows in $C$ are $r$-cocartesian. Moreover, since $\rho$ is $\ev_2$-cartesian the map $(\mu(1) \leftarrow \mu(0)) \rightarrow (\mu'(1) \leftarrow \mu'(0))$ is $\ev_0$-cartesian. The fact that $F\phi$ is a morphism of two-sided fibrations then implies that the left vertical arrow in $C$ is $r_2$-cartesian. The image of $C$ under $r$ is the commutative square
\[
\begin{tikzcd}
(\mu(1), \mu(0)) \arrow{r}{} \arrow{d}{} & (\mu(1), \mu(2)) \arrow{d}{} \\ (\mu'(1), \mu'(0)) \arrow{r}{} & (\mu'(1), \mu'(2))
\end{tikzcd}
\]
which is cartesian and has constant first coordinate since $\rho$ was taken to be $\ev_2$-cartesian. Since $r$ satisfies the Beck-Chevalley condition in the second coordinate we conclude that the right vertical arrow in $C$ is $r$-cartesian, as desired.
\end{proof}
 

%% file: Two_corr.tex

\tableofcontents

\section{The $2$-category of correspondences} \label{section twocorr}

Let $\Ccal$ be a category admitting pullbacks. We can attach to $\Ccal$ a $2$-category $\twoCorr(\Ccal)$ called the $2$-category of correspondences of $\Ccal$. Its space of objects coincides with the space of objects of $\Ccal$, and for each pair of objects $c, c'$ in $\Ccal$, the hom category $\Hom_{\twoCorr(\Ccal)}(c, c')$ is the category of diagrams in $\Ccal$ of the form
\[
\begin{tikzcd}
  & s \arrow[ld] \arrow[rd] &    \\
c &                         & c'.
\end{tikzcd}
\]
 Our goal in this section is to review the definition and main properties of the $2$-category of correspondences\footnote{Many of the basic properties of $\twoCorr(\Ccal)$ that we discuss in \ref{subsection 2Corr} - \ref{subsection adj}  (namely, its construction, symmetric monoidal structure, adjunctions, duals) can be found in some way in \cite{HaugSpan} or \cite{GR}. We chose to include statements and proofs of these facts for completenes and for ease of reference, as our notation differs from that of previous sources. Some of our proofs of these facts contain some level of novelty - for instance we prove the adjointness and duality properties of $\twoCorr(\Ccal)$ by appeal to the functoriality of $\twoCorr$, therefore reducing to checking that they hold in the universal examples, which is often manageable.}, and to provide a new proof of its universal property.

We begin in \ref{subsection 2Corr}  by recalling the definition and basic properties of $\twoCorr(\Ccal)$. We define $\twoCorr(\Ccal)$ first as a simplicial category, and show that it is in fact a complete Segal object in $\Cat$, so it defines a $2$-category. We provide here a description of the degeneracies and composition maps for $\twoCorr(\Ccal)$.

In \ref{subsection functoriality twocorr} we study the functoriality of the assignment $\Ccal \mapsto \twoCorr(\Ccal)$. We construct $\twoCorr$ as a functor on the category of categories with pullbacks and pullback preserving morphisms. We show that this is in fact a limit preserving functor - in particular, if $\Ccal$ comes equipped with a symmetric monoidal structure which is compatible with pullbacks, we obtain an induced symmetric monoidal structure on $\twoCorr(\Ccal)$.

The $2$-category $\twoCorr(\Ccal)$ comes equipped with inclusions $\iota_\Ccal: \Ccal \rightarrow \twoCorr(\Ccal)$ and $\iota_\Ccal^R : \Ccal^\op \rightarrow \twoCorr(\Ccal)$. We show that these also depend functorially on $\Ccal$. In particular, in the presence of a symmetric monoidal structure on $\Ccal$ compatible with pullbacks, the inclusions $\iota_\Ccal$ and $\iota_\Ccal^R$ inherit  canonical symmetric monoidal structures.

In \ref{subsection adj} we record two basic dualizability and adjointness properties of the $2$-category of correspondences. We show that if $\alpha$ is an arrow in $\Ccal$ then the maps $\iota_\Ccal(\alpha)$ and $\iota_\Ccal^R(\alpha)$ are adjoint to each other. In the case when $\Ccal$ admits finite limits, we show that every object of $\twoCorr(\Ccal)$ is dualizable under the symmetric monoidal structure inherited from the cartesian symmetric monoidal structure on $\Ccal$. The proof of these results appeals in a fundamental way to the functoriality properties of $\twoCorr$: rather than showing that a candidate (co)unit indeed defines a duality or adjunction on $\twoCorr(\Ccal)$ for an arbitrary $\Ccal$, one proves it in the universal example, in which case the verification of the triangle identities becomes simple.

 In \ref{subsection BC} we review the so-called Beck-Chevalley condition, and provide a new proof of the fact that the inclusion $\iota_\Ccal: \Ccal \rightarrow \twoCorr(\Ccal)$ is the universal embedding of $\Ccal$ into a $2$-category satisfying the left Beck-Chevalley condition.   This provides a concrete way of constructing functors out of $\twoCorr(\Ccal)$: given a $2$-category $\Dcal$, a functor $\twoCorr(\Ccal) \rightarrow \Dcal$ is the same data as a functor $\Ccal \rightarrow \Dcal$ satisfying familiar base change properties.

\subsection{Construction and basic properties} \label{subsection 2Corr}
We begin by reviewing the construction of the $2$-category of correspondences.

\begin{notation}
Let $n \geq 0$ and let $\Tw([n])$ be the twisted arrow category of $[n]$. We identify the objects in $\Tw([n])$ with pairs $(i, j)$ in $[n]\times [n]$ such that $i \leq j$, so that there is a unique arrow $(i, j) \rightarrow (i', j')$ whenever $i \leq i' \leq j' \leq j$.  We denote by $\Tw([n])_\el$ the full subcategory of $\Tw([n])$ on the objects of the form $(i, i+1)$ for $0 \leq i < n$.
\end{notation}

\begin{definition}
Let $\Ccal$ be a category admitting pullbacks and let $n \geq 0$. We say that a functor $\Tw([n]) \rightarrow \Ccal$ is cartesian if it is the right Kan extension of its restriction to $\Tw([n])_\el$.
\end{definition}

\begin{proposition}\label{prop equiv being cartesian}
Let $\Ccal$ be a category admitting pullbacks and let $n \geq 0$. Let $S: \Tw([n]) \rightarrow \Ccal$ be a functor. Then the following conditions are equivalent:
\begin{enumerate}[\normalfont (i)]
\item The functor $S$ is cartesian.
\item For every object $(i, j)$ in $\Tw([n])$ such that $j \geq i+2$ the commutative square
\[
\begin{tikzcd}
S(i, j) \arrow{d}{} \arrow{r}{} & S(i, j-1) \arrow{d}{} \\
S(i+1, j) \arrow{r}{} & S(i+1, j-1)
\end{tikzcd}
\]
is cartesian.
\end{enumerate}
\end{proposition}
\begin{proof}
 For each integer $1 \leq k \leq n$ let $\Tw([n])_{\leq k}$ be the full subcategory of $\Tw([n])$ on the objects of the form $(i, j)$ with $j - i \leq k$. Note that if $2 \leq k \leq n$ and $(i,j)$ is such that $j-i = k$ then the undercategory $((\Tw([n])_{\leq k-1})_{(i,j)/}$ contains the diagram 
 \[(i+1,j) \rightarrow (i+1, j-1) \leftarrow (i, j-1).
 \]
 This diagram is in fact final in $((\Tw([n])_{\leq k-1})_{(i,j)/}$, and therefore we have that the restriction of $S$ to $(\Tw([n]))_{\leq k}$ is the right Kan extension of its restriction to $(\Tw([n]))_{\leq k-1}$ if and only if the diagram in the statement is cartesian for every $(i, j)$ such that $j-i = k$. The result now follows by induction on $k$.
\end{proof}

\begin{notation}
Let $\Ccal$ be a category admitting pullbacks. Let $\overline{\twoCorr}(\Ccal): \Delta^\op \rightarrow \Cat$ be the simplicial category given by the formula $\overline{\twoCorr}(\Ccal)([n]) = \Funct(\Tw([n]), \Ccal)$. By virtue of proposition \ref{prop equiv being cartesian}, for each map $[n] \rightarrow [n']$ in $\Delta$, the induced map $\overline{\twoCorr}(\Ccal)([n']) \rightarrow \overline{\twoCorr}(\Ccal)([n]))$ sends cartesian objects to cartesian objects. We denote by $\twoCorr'(\Ccal)$ the sub-simplicial category of $\overline{\twoCorr}(\Ccal)$ such that for every $[n]$ in $\Delta$ the category $\twoCorr'(\Ccal)([n])$ is the full subcategory of $\overline{\twoCorr}(\Ccal)([n])$ on the cartesian objects.
\end{notation}

\begin{proposition}\label{prop twocorrprime is segal cat}
Let $\Ccal$ be a category admitting pullbacks. Then $\twoCorr'(\Ccal)$ is a Segal category. 
\end{proposition}
\begin{proof}
Let $n \geq 0$ and denote by $\spn([n])$ the spine of $[n]$, that is, the union inside $\Pcal(\Delta)$ of all the edges of $[n]$ of the form $i \rightarrow i+1$. Any simplicial category $\Scal$ determines by right Kan extension a functor $\Pcal(\Delta)^\op \rightarrow \Cat$, which we will also denote by $\Scal$.

We have a commutative diagram
\[
\begin{tikzcd}
\overline{\twoCorr}(\Ccal)([n]) \arrow{r}{} & \overline{\twoCorr}(\Ccal)(\spn[n])  \\ \arrow{u}{}
\twoCorr'(\Ccal)([n]) \arrow{r}{} & \arrow{u}{} \twoCorr'(\Ccal)(\spn[n]) .
\end{tikzcd}
\]
Since $\twoCorr'(\Ccal)$ and $\overline{\twoCorr}(\Ccal)$ agree on simplices of dimension at most $1$, the right vertical arrow is an isomorphism. Observe that the left Kan extension along $\Delta \rightarrow \Pcal(\Delta)$ of the cosimplicial category $\Tw|_\Delta$ maps the inclusion $\spn([n]) \rightarrow [n]$ to the inclusion $\Tw([n])_\el \rightarrow \Tw([n])$. It follows that the top horizontal arrow in the above diagram is equivalent to the  restriction map
\[
\Funct(\Tw([n]), \Ccal) \rightarrow \Funct(\Tw([n])_{\el}, \Ccal).
\]
It follows from the definition of $\twoCorr'(\Ccal)([n])$ that restriction along the inclusion $\Tw([n])_\el \rightarrow \Tw([n])$ provides an equivalence $\twoCorr'(\Ccal)([n]) = \Funct(\Tw([n])_\el, \Ccal)$. Therefore the bottom horizontal arrow in the diagram is an equivalence, which means that $\twoCorr'(\Ccal)$ satisfies the Segal conditions, as desired.
\end{proof}

Note that the category $\twoCorr'(\Ccal)([0])$ is equivalent to $\Ccal$. We now consider the Segal category obtained from $\twoCorr'(\Ccal)$ by discarding noninvertible arrows in $\twoCorr'(\Ccal)([0])$.

\begin{notation}\label{notation red}
Denote by $\codiscr: \Cat = \Funct([0], \Cat) \rightarrow \Funct(\Delta^\op, \Cat)$ the functor of right Kan extension along the inclusion $\lbrace [0] \rbrace \rightarrow \Delta^\op$.  We have a  diagram of endofunctors of $\Funct(\Delta^\op, \Cat)$ as follows
\[
\begin{tikzcd}
 & \id_{\Funct(\Delta^\op, \Cat)} \arrow{d}{} \\
 \codiscr \ev_{[0]}^{\leq 0} \arrow{r}{} & \codiscr \ev_{[0]}
\end{tikzcd}
\]
where the vertical arrow is the unit of the adjunction $\ev_{[0]} \dashv \codiscr$ and the horizontal arrow is the canonical inclusion. We let $(-)_\red: \Funct(\Delta^\op, \Cat) \rightarrow \Funct(\Delta^\op, \Cat)$ be the fiber product of the above diagram. 
\end{notation}

\begin{remark}
The functor $\codiscr$ can alternatively described as the functor induced from the composite map
\[
\Delta^\op \times \Cat \xrightarrow{(-)^{\leq 0} \times \id_{\Cat}} \Cat^\op \times \Cat \xrightarrow{\Funct(-,-)} \Cat.
\]
In other words, for each category $\Ccal$ the simplicial category $\codiscr(\Ccal)$ has $\Ccal^n$ as its category of $[n]$-simplices. Note that this is in fact a Segal category.
\end{remark}

\begin{remark}
The functor $(-)_\red$ sends a simplicial category $\Scal$ to a simplicial subcategory $\Scal_{\red}$ of $\Scal$ such that for each $n \geq 0$ we have that $\Scal_\red([n])$ is the subcategory of $\Scal([n])$ containing all objects, and only those arrows whose images under the $n+1$ functors $\Scal([n]) \rightarrow \Scal([0])$ are invertible. It follows from its presentation as the fiber product $\codiscr(\Scal([0])^{\leq 0}) \times_{\codiscr(\Scal([0]))} \Scal$ that the map $\Scal_\red \rightarrow \Scal$ is universal among the maps of simplicial categories $\Scal' \rightarrow \Scal$ such that $\Scal'([0])$ is a space. Moreover, since codiscrete simplicial categories are Segal, we see that if $\Scal$ is a Segal category then $\Scal_\red$ is also a Segal category. 
\end{remark}

\begin{notation}
Let $\Ccal$ be a category admitting pullbacks. We denote by $\twoCorr(\Ccal)$ the Segal category $\twoCorr'(\Ccal)_\red$.
\end{notation}

\begin{remark}\label{remark hom categories}
Let $\Ccal$ be a category admitting pullbacks. The space of objects of $\twoCorr(\Ccal)$ is the space of maps $\Tw([0]) = [0] \rightarrow \Ccal$, and it therefore agrees with the space of objects of $\Ccal$. Given two objects, $c, c'$ in $\Ccal$, the category of morphisms $\Hom_{\twoCorr(\Ccal)}(c, c')$ is  the overcategory 
\[
\Ccal_{/c, c'} = \Funct(\Tw([1]), \Ccal) \times_{\Ccal \times \Ccal} \lbrace (c, c') \rbrace
\]
where the projection $\Funct(\Tw([1]), \Ccal) \rightarrow \Ccal\times \Ccal$ is given by evaluation at $(0,0)$ and $(1,1)$.

 The objects of $\Hom_{\twoCorr(\Ccal)}(c, c')$ are therefore spans $c \leftarrow s \rightarrow c'$, and a morphism $(c \leftarrow s \rightarrow c') \rightarrow (c \leftarrow t \rightarrow c')$ is a commutative diagram
\[
\begin{tikzcd}
  & s \arrow[ldd] \arrow[rdd] \arrow[d] &    \\
  & t \arrow[ld] \arrow[rd]             &    \\
c &                                     & c'.
\end{tikzcd}
\]
The degeneracy map of $\twoCorr(\Ccal)$ is given by precomposition with the projection $\Tw([1]) \rightarrow \Tw([0]) = [0]$. We thus see that for every object $c$ in $\Ccal$ the identity in $\Hom_\Ccal(c, c)$ is given by the span $c \xleftarrow{\id_c} c \xrightarrow{\id_c} c$.

Consider the diagram of simplicial spaces $\spn([2]) \rightarrow [2] \leftarrow [1]$, where the right arrow is the unique active map from $[1]$ to $[2]$, and the left arrow is the inclusion of the spine $\spn([2]) =[1] \cup_{[0]} [1]$ inside $[2]$. As in the proof of proposition \ref{prop twocorrprime is segal cat}, via right Kan extension along the inclusion $\Delta^\op \rightarrow \Pcal(\Delta)^\op$, we allow ourselves to evaluate Segal categories on such a diagram. Applying this to the maps $\twoCorr'(\Ccal) \rightarrow \overline{\twoCorr}(\Ccal) \rightarrow \codiscr(\Ccal)$ we obtain a commutative diagram of categories 
\[
\begin{tikzcd}
\twoCorr'(\spn([2])) \arrow{d}{}{} & \arrow{l}{} \twoCorr'(\Ccal)([2]) \arrow{d}{} \arrow{r}{} & \twoCorr'(\Ccal)([1]) \arrow{d}{} \\
\overline{\twoCorr}(\spn([2])) \arrow{d}{}{} & \arrow{l}{} \overline{\twoCorr}(\Ccal)([2]) \arrow{d}{} \arrow{r}{} & \overline{\twoCorr}(\Ccal)([1]) \arrow{d}{}\\
\codiscr(\Ccal)(\spn([2])) & \arrow{l}{} \codiscr(\Ccal)([2]) \arrow{r}{} & \codiscr(\Ccal)([1]).
\end{tikzcd}
\]
Here the top left and bottom left horizontal arrows are isomorphisms, the top left and top right vertical arrows are isomorphisms, and the top middle vertical arrow is a monomorphism. The top and bottom left squares are horizontally right adjointable. Using the canonical identification $\Tw([2])_\el = \Tw([1]) \cup_{[0]} \Tw([1])$ we see that the bottom two rows of the diagram can be rewritten as follows:
\[
\begin{tikzcd}
\twoCorr'(\Ccal)(\spn[2]) \arrow{d}{} & \arrow{l}{} \twoCorr'(\Ccal)([2]) \arrow{d}{} \arrow{r}{} & \twoCorr'(\Ccal)([1]) \arrow{d}{} \\
\Funct(\Tw([2])_\el, \Ccal) \arrow{d}{}{} & \arrow{l}{} \Funct(\Tw([2]),\Ccal) \arrow{d}{} \arrow{r}{} & \Funct(\Tw([1]),\Ccal) \arrow{d}{}\\
\Ccal \times \Ccal \times \Ccal & \arrow{l}{} \Ccal \times \Ccal \times \Ccal \arrow{r}{} & \Ccal \times \Ccal.
\end{tikzcd}
\]
Here the bottom vertical arrows are given by evaluation at the objects of the form $(i, i)$, the bottom left vertical arrow is the identity, the bottom right arrow is the projection onto the first and third coordinate, the middle left horizontal arrow is restriction along the inclusion $\Tw([2])_\el \rightarrow \Tw([2])$, and the middle right horizontal arrow is induced from the active map $[1] \rightarrow [2]$. 

The top and bottom left commutative squares in the above diagram are horizontally right adjointable. Passing to horizontal right adjoints of these yields a commutative diagram
\[
\begin{tikzcd}
\twoCorr'(\Ccal)(\spn([2])) \arrow{d}{} \arrow{r}{} & \twoCorr'(\Ccal)([2]) \arrow{d}{} \arrow{r}{} &  \twoCorr'(\Ccal)([1]) \arrow{d}{} \\
\Funct(\Tw([2])_\el, \Ccal) \arrow{d}{}{} \arrow{r}{} & \Funct(\Tw([2]),\Ccal) \arrow{d}{} \arrow{r}{} & \Funct(\Tw([1]),\Ccal) \arrow{d}{}\\
\Ccal \times \Ccal \times \Ccal \arrow{r}{} & \Ccal \times \Ccal \times \Ccal \arrow{r}{} & \Ccal \times \Ccal.
\end{tikzcd}
\]
Here the middle left horizontal arrow is given by right Kan extension along the inclusion $\Tw([2])_\el \rightarrow \Tw([2])$. The top row recovers the composition map for the Segal category $\twoCorr'(\Ccal)$. The composition map for $\twoCorr(\Ccal)$ can be obtained from the above by base change along the inclusion of \[
\Ccal^{\leq 0} \times \Ccal^{\leq 0} \times \Ccal^{\leq 0} \longrightarrow \Ccal^{\leq 0} \times \Ccal^{\leq 0} \times \Ccal^{\leq 0} \longrightarrow \Ccal^{\leq 0} \times \Ccal^{\leq 0}
\] into the bottom row.

Let $c, c', c''$ be a triple of objects in $\Ccal$. The bottom row receives a map from the final row $[0] \rightarrow [0] \rightarrow [0]$ that picks out the triple $(c, c', c'')$ inside $\Ccal \times \Ccal \times \Ccal$. Base change along this map yields a commutative diagram
\[
\begin{tikzcd}
\Funct(\Tw([2])_\el, \Ccal)_{(c, c', c'')} \arrow{d}{} \arrow{r}{} & \twoCorr(\Ccal)([2])_{(c, c', c'')} \arrow{d}{} \arrow{r}{} & \Funct(\Tw([1]), \Ccal)_{(c, c'')}\arrow{d}{} \\
\Funct(\Tw([2])_\el, \Ccal)_{(c, c', c'')}\arrow{r}{} & \Funct(\Tw([2]),\Ccal)_{(c, c', c'')} \arrow{r}{} & \Funct(\Tw([1]),\Ccal)_{(c, c'')}
\end{tikzcd}.
\]
The composition of the two maps in the top row recovers the composition map
 \[
\Hom_{\twoCorr(\Ccal)}(c, c') \times \Hom_{\twoCorr(\Ccal)}(c', c'') \rightarrow \Hom_{\twoCorr(\Ccal)}(c, c'').
\] 
In particular, we see that the composition map for $\twoCorr(\Ccal)$ sends a pair of spans $c \leftarrow s \rightarrow c'$ and $c' \leftarrow s' \rightarrow c''$ to the span $c \leftarrow s\times_{c'} s' \rightarrow c''$. Given a pair of morphisms of spans
\[
\begin{tikzcd}
  & s \arrow[ldd] \arrow[rdd] \arrow[d] &    \\
  & t \arrow[ld] \arrow[rd]             &    \\
c &                                     & c'
\end{tikzcd} 
\begin{tikzcd}
  & s' \arrow[ldd] \arrow[rdd] \arrow[d] &    \\
  & t' \arrow[ld] \arrow[rd]             &    \\
c' &                                     & c''
\end{tikzcd}
\]
their composition is the morphism $(c \leftarrow s \times_{c'} s' \rightarrow c'') \rightarrow (c \leftarrow t \times_{c'} t'  \rightarrow c'')$ that arises from the unique commutative diagram of the form
\[
\begin{tikzcd}
  &                                     & s \times_{c'} s' \arrow[ld] \arrow[rd] \arrow[d] &                                      &     \\
  & s \arrow[ldd] \arrow[d] \arrow[rdd] & t\times_{c'} t' \arrow[ld] \arrow[rd]            & s' \arrow[ldd] \arrow[rdd] \arrow[d] &     \\
  & t \arrow[ld] \arrow[rd]             &                                                  & t' \arrow[ld] \arrow[rd]             &     \\
c &                                     & c'                                               &                                      & c''
\end{tikzcd}
\]
that extends the given morphisms of spans.
\end{remark}

\begin{proposition}\label{prop isos in corr}
Let $\Ccal$ be a category admitting pullbacks. Then a span $c \leftarrow s \rightarrow c'$ is invertible in $\twoCorr(\Ccal)$ if and only if both legs are isomorphisms.
\end{proposition}
\begin{proof}
If both legs of a span $c \leftarrow s \rightarrow c'$ are invertible, then it is equivalent to an identity span, which is invertible. Conversely, assume that the span $c \leftarrow s \rightarrow c'$ is invertible and let $c' \leftarrow t \rightarrow c$ be the inverse. Consider the composite span
\[
\begin{tikzcd}
  &                         & u \arrow[ld] \arrow[rd] &                         &     \\
  & s \arrow[rd] \arrow[ld] &                                     & t \arrow[rd] \arrow[ld] &     \\
c &                         & c'                                  &                         & c
\end{tikzcd}
\]  
where $u = s \times_{c'} t$. The maps $c \leftarrow u$ and $u \rightarrow c''$ are invertible, and therefore we see that the projections $c \leftarrow s$ and $t \rightarrow c$ admit sections. Considering the composition in the opposite order reveals that the projections $s \rightarrow c'$ and $c' \leftarrow t$ also admit sections. In particular we conclude that the projection $s \leftarrow u$ admits a section, and therefore the map $c \leftarrow s$ is an isomorphism. Similarly, the projection $u \rightarrow t$ admits a section, and thus the map $t \rightarrow c$ is also an isomorphism.
\end{proof}

\begin{corollary}\label{coro twocorr is a 2cat}
Let $\Ccal$ be a category admitting pullbacks. Then the Segal space underlying $\twoCorr(\Ccal)$ is complete.
\end{corollary}
\begin{proof}
Denote by $\twoCorr(\Ccal)([1])^\iso$ the space of invertible $1$-morphisms in $\twoCorr(\Ccal)$ and by $\Ccal([1])^\iso$ the space of invertible arrows in $\Ccal$. Proposition \ref{prop isos in corr} implies that $\twoCorr(\Ccal)([1])^\iso$ is equivalent to the fiber product $\Ccal([1])^\iso \times_{\Ccal^{\leq 0}} \Ccal([1])^\iso$, where the fiber product is taken with respect to the source projection in both coordinates. Our claim now follows from the fact that the degeneracy $\Ccal^{\leq 0} \rightarrow \Ccal([1])^\iso$ is an equivalence.
\end{proof}

 In other words, corollary \ref{coro twocorr is a 2cat} states that $\twoCorr(\Ccal)$ belongs to the image of the inclusion $\twoCat \rightarrow \Funct(\Delta^\op, \Cat)$. 
 
 \begin{definition}
 Let $\Ccal$ be a category admitting pullbacks. We call $\twoCorr(\Ccal)$ the $2$-category of correspondences of $\Ccal$.
 \end{definition}
 
 \begin{remark}\label{remark inclusiones C y Cop}
 Let $\Ccal$ be a category admitting pullbacks. For each $n \geq 0$ we have source and target projections $[n] \leftarrow \Tw([n]) \rightarrow [n]^\op$ that map the pair $(i, j)$ to $i$ and $j$ respectively. These projections are natural in $[n]$, and precomposing with them yields functors
 \[
 \Ccal \xrightarrow{\iota_\Ccal} \twoCorr(\Ccal) \xleftarrow{\iota_\Ccal^R} \Ccal^\op.
 \] 
 The above functors are the identity on objects. Moreover, for every arrow $\alpha: c \rightarrow c'$ in $\Ccal$, we have $\iota_\Ccal(\alpha) = (c \xleftarrow{\id_c} c \xrightarrow{\alpha} c')$ and $\iota^R_{\Ccal}(\alpha) =  (c' \xleftarrow{\alpha} c \xrightarrow{\id_{c}} c)$.
 \end{remark}

\subsection{Functoriality of $\twoCorr$}\label{subsection functoriality twocorr}
We now examine the functoriality of the assignment $\Ccal \mapsto \twoCorr(\Ccal)$.

\begin{construction}\label{construction twocorr}
We denote by $\Cat_\pb$ the subcategory of $\Cat$ on the categories admitting pullbacks, and functors which preserve pullbacks. Let $\overline{\twoCorr}: \Cat_\pb \rightarrow \Funct(\Delta^\op, \Cat)$ be the functor induced from the composite map
\[
\Delta^\op \times \Cat_\pb \xrightarrow{\Tw(-) \times \id_{\Cat_\pb}} \Cat^\op \times \Cat_\pb \subset \Cat^\op \times \Cat \xrightarrow{\Funct(-, -)} \Cat.
\]

It follows from proposition \ref{prop equiv being cartesian} that for every morphism $\Ccal \rightarrow \Dcal$ in $\Cat_\pb$ the induced morphism of simplicial categories $\overline{\twoCorr}(\Ccal) \rightarrow \overline{\twoCorr}(\Dcal)$ restricts to a morphism $\twoCorr'(\Ccal) \rightarrow \twoCorr'(\Dcal)$. This in turn restricts to a functor $\twoCorr(\Ccal) \rightarrow \twoCorr(\Dcal)$. We therefore have that the assignment $\Ccal \mapsto \twoCorr(\Ccal)$ extends to a functor $\twoCorr: \Cat_\pb \rightarrow \twoCat$ whose composition with the embedding $\twoCat \rightarrow \Funct(\Delta^\op, \Cat)$ is a subfunctor of $\overline{\twoCorr}$.
\end{construction}

\begin{proposition}\label{prop twocorr pres limits}
The functor $\twoCorr: \Cat_{\text{\normalfont pb}} \rightarrow \twoCat$ preserves limits.
\end{proposition}
\begin{proof}
Denote by $i: \twoCat \rightarrow \Funct(\Delta^\op, \Cat)$ the inclusion. To show that $\twoCorr$ preserves limits, it suffices to show that the functor $\ev_{[n]} i \twoCorr: \Cat_\pb \rightarrow \Cat$ preserves limits for every $n \geq 0$. We now fix  a category $\Jcal$ and a limit diagram $F: \Jcal^\lhd \rightarrow \Cat_\pb$. Denote by $\ast$ the initial object of $\Jcal^\lhd$.

 By construction, the functor $\ev_{[n]} \overline{\twoCorr}$ is given by the composition of the forgetful functor $\Cat_\pb \rightarrow \Cat$ and the functor $\Funct(\Tw([n]), -): \Cat \rightarrow \Cat$. Both of these preserve limits, so we see that $\ev_{[n]}\overline{\twoCorr} F$ is a limit diagram. The natural functor 
 \[
 \lim (\ev_{[n]} i \twoCorr F|_\Jcal) \rightarrow \ev_{[n]}\overline{\twoCorr} F(\ast)
 \]
  is a limit of monomorphisms, and is therefore a monomorphism. 
 
 An object $S$ in $\ev_{[n]}\overline{\twoCorr} F(\ast)$ belongs to $\lim (\ev_{[n]} i \twoCorr F|_\Jcal)$ if and only if its projection to $\ev_{[n]}\overline{\twoCorr} F(j)$ belongs to $\ev_{[n]}i\twoCorr F (j)$ for every $j$ in $\Jcal$. It follows from proposition \ref{prop equiv being cartesian} and the fact that the transition maps in $F$ preserve pullbacks that this happens if and only if $S$ belongs to $\ev_{[n]}i\twoCorr F(\ast)$.
 
Similarly, a morphism $g: S \rightarrow S'$ in $\ev_{[n]}\overline{\twoCorr} F(\ast)$ belongs to $\lim (\ev_{[n]} i \twoCorr F|_\Jcal)$ if and only if its projection to $\ev_{[n]}\overline{\twoCorr} F(j)$ belongs to $\ev_{[n]}i\twoCorr F (j)$ for every $j$ in $\Jcal$. This again happens if and only if $g$ belongs to $\ev_{[n]}i\twoCorr F(\ast)$.

We have thus seen that the functor $\lim (\ev_{[n]} i \twoCorr F|_\Jcal) \rightarrow \ev_{[n]}\overline{\twoCorr} F(\ast)$ is a monomorphism, and its image coincides with $\ev_{[n]}i\twoCorr F(\ast)$. This implies that $\ev_{[n]} i \twoCorr F$ is indeed a limit diagram.
\end{proof}

\begin{remark}\label{remark symmetric monoidal structure}
Equip $\Cat_\pb$ and $\twoCat$ with their cartesian symmetric monoidal structures. It follows from proposition \ref{prop twocorr pres limits} that $\twoCorr$ has a canonical symmetric monoidal structure. As a consequence, if $\Ccal$ is a symmetric monoidal category admitting pullbacks and such that the tensor product functor $\Ccal \times \Ccal \rightarrow \Ccal$ preserves pullbacks, there is an induced symmetric monoidal structure on $\twoCorr(\Ccal)$. In particular, given any finitely complete category $\Ccal$ there is a symmetric monoidal structure on $\twoCorr(\Ccal)$ inherited from the cartesian symmetric monoidal structure on $\Ccal$. This assignment is functorial in $\Ccal$ - namely, it can be enhanced to yield functors
\[
\Cat_\lex \rightarrow \CAlg(\Cat_\pb) \rightarrow \CAlg(\twoCat)
\]
where $\Cat_\lex$ denotes the category if categories with finite limits and left exact functors.
\end{remark}

The inclusions $\iota_\Ccal$ and $\iota_\Ccal^R$ from remark \ref{remark inclusiones C y Cop} turn out to be compatible with the symmetric monoidal structure of remark \ref{remark symmetric monoidal structure}. To show this, we will need to make the transformations $\iota_\Ccal$ and $\iota_\Ccal^R$ functorial in $\Ccal$.

\begin{construction}\label{constr natural transf}
Consider the commutative diagram of cosimplicial categories
\[
\begin{tikzcd}
\Tw|_\Delta \arrow{r}{} \arrow{d}{} & i_{\Delta}^\op \arrow{d}{} \\{ i_\Delta }\arrow{r}{} & {[0]}
\end{tikzcd}
\]
where $i_\Delta: \Delta \rightarrow \Cat$ is the canonical inclusion, and the left vertical and top horizontal maps are given by the source and target projections (see remark \ref{remark inclusiones C y Cop}). This induces a commutative square of functors $\Cat_\pb \rightarrow \Funct(\Delta^\op, \Cat)$ as follows
\[
\begin{tikzcd}
\operatorname{disc} \arrow{d}{} \arrow{r}{} & j(-)^\op \arrow{d}{} \\ j \arrow{r}{} & \overline{\twoCorr}
\end{tikzcd}
\]
where $\operatorname{disc}$ denotes the functor that maps each category $\Ccal$ to the constant simplicial category on $\Ccal$, the functor $j$ is given by the formula $j\Ccal([n]) = \Funct([n], \Ccal)$, and $j(-)^\op$ denotes the composition of $j$ with the functor $(-)^\op : \Cat \rightarrow \Cat$. The bottom horizontal and right vertical arrows factor through $\twoCorr'$, so we obtain a commutative square of functors $\Cat_\pb \rightarrow \Funct(\Delta^\op, \Cat)$ as follows:
\[
\begin{tikzcd}
\operatorname{disc} \arrow{d}{} \arrow{r}{} & j(-)^\op \arrow{d}{} \\ j \arrow{r}{} & \twoCorr'
\end{tikzcd}
\]
Composing with the functor $(-)_\red$ from notation \ref{notation red} we obtain a commutative square of functors $\Cat_\pb \rightarrow \twoCat$
\[
\begin{tikzcd}
i_{\Spc}(-)^{\leq 0} \arrow{d}{} \arrow{r}{} & i^\op_{\Cat_\pb} \arrow{d}{} \\ i_{\Cat_\pb} \arrow{r}{} & \twoCorr
\end{tikzcd}
\]
where $i_{\Cat_\pb}$ denotes the canonical inclusion $\Cat_\pb \rightarrow \twoCat$, the functor $i_{\Cat_\pb}^\op$ is the composition of $i_{\Cat_\pb}$ with the functor $(-)^\op: \Cat_\pb \rightarrow \Cat_\pb$, and the top left corner is the composition of the truncation functor $(-)^{\leq 0}: \Cat \rightarrow \Spc$ and the canonical inclusion $i_{\Spc}:\Spc \rightarrow \twoCat$. 

We denote by $\iota$ the bottom horizontal arrow of the above diagram, and $\iota^R$ the right vertical arrow. When evaluated at a category $\Ccal$ in $\Cat_\pb$, the above diagram recovers the commutative diagram
\[
\begin{tikzcd}
\Ccal^{\leq 0} \arrow{r}{} \arrow{d}{} & \Ccal^\op \arrow{d}{\iota^R_\Ccal} \\ \Ccal \arrow{r}{\iota_\Ccal} & \twoCorr(\Ccal)
\end{tikzcd}
\]
from remark \ref{remark inclusiones C y Cop}.
\end{construction}

\begin{remark}\label{remark sm structure 2}
We can think about the commutative diagram
 \[
\begin{tikzcd}
i_{\Spc}(-)^{\leq 0} \arrow{d}{} \arrow{r}{} & i^\op_{\Cat_\pb} \arrow{d}{\iota^R} \\ i_{\Cat_\pb} \arrow{r}{\iota} & \twoCorr
\end{tikzcd}
\]
as a functor $\Cat_\pb \rightarrow \Funct([1] \times [1], \twoCat)$. Thanks to proposition \ref{prop twocorr pres limits}, this functor is limit preserving, so it can be given a canonical symmetric monoidal structure, where we equip $\Cat_\pb$ and $\twoCat$ with their cartesian  symmetric monoidal structures. It follows that if $\Ccal$ is a symmetric monoidal category admitting pullbacks and such that the tensor product functor $\Ccal \times \Ccal \rightarrow \Ccal$ preserves pullbacks, the inclusions $\iota_\Ccal$ and $\iota_\Ccal^R$ can be given symmetric monoidal structures, and we have a commutative diagram of symmetric monoidal categories and symmetric monoidal functors 
\[
\begin{tikzcd}
\Ccal^{\leq 0} \arrow{r}{} \arrow{d}{} & \Ccal^\op \arrow{d}{\iota^R_\Ccal} \\ \Ccal \arrow{r}{\iota_\Ccal} & \twoCorr(\Ccal).
\end{tikzcd}
\]
\end{remark}

\begin{remark}\label{remark passing op}
The span of cosimplicial categories 
\[
\begin{tikzcd}
\Tw|_\Delta \arrow{r}{} \arrow{d}{} & i_{\Delta}^\op \arrow{d}{} \\{ i_\Delta }\arrow{r}{} & {[0]}
\end{tikzcd}
\]
from  construction \ref{constr natural transf} is equivalent to the transpose of the diagram obtained by composing the above with the functor $\Funct(\Delta, \Cat) \rightarrow \Funct(\Delta, \Cat)$ induced from the order reversing automorphism of $\Delta$. It follows that the commutative square of functors $\Cat_\pb \rightarrow \twoCat$
\[
\begin{tikzcd}
i_{\Spc}(-)^{\leq 0} \arrow{d}{} \arrow{r}{} & i^\op_{\Cat_\pb} \arrow{d}{\iota^R} \\ i_{\Cat_\pb} \arrow{r}{\iota} & \twoCorr
\end{tikzcd}
\]
is equivalent to the transpose of the square
\[
\begin{tikzcd}
i_{\Spc}(-)^{\leq 0} \arrow{d}{} \arrow{r}{} & i_{\Cat_\pb} \arrow{d}{\iota^{R, 1\dsh\op}} \\ i^\op_{\Cat_\pb} \arrow{r}{\iota^{1\dsh\op}} & \twoCorr^{1\dsh\op}
\end{tikzcd}
\]
where $(-)^{1\dsh\op}: \twoCat \rightarrow \twoCat$ denotes the functor that reverses the directions of $1$-arrows. In other words, for every object $\Ccal$ in $\Cat_\pb$ we have an equivalence $\twoCorr(\Ccal) = \twoCorr(\Ccal)^{1\dsh\op}$ which is the identity on objects, and exchanges $i_\Ccal(\alpha)$ and $i^R_\Ccal(\alpha)$ for each arrow $\alpha$ in $\Ccal$.
\end{remark}

\subsection{Adjointness and duality in $\twoCorr(\Ccal)$} \label{subsection adj}

We now review some basic adjointness and duality properties of morphisms and objects in the $2$-category of correspondences.

\begin{proposition}\label{prop adjunctos in twocrr}
Let $\Ccal$ be a category admitting pullbacks and let $\alpha: c \rightarrow c'$ be an arrow in $\Ccal$. Then the morphism $\eta_\alpha: \iota_\Ccal(\alpha) \iota^R_{\Ccal}(\alpha) \rightarrow \id_{\iota_\Ccal(c')}$ given by the following diagram
\[
\begin{tikzcd}
   & c \arrow{ldd}[swap]{\alpha} \arrow[rdd, "\alpha"] \arrow[d, "\alpha"] &    \\
   & c' \arrow[ld, "\id_{c'}"] \arrow{rd}[swap]{\id_{c'}}                   &    \\
c' &                                                                   & c'
\end{tikzcd}
\]
is the counit of an adjunction between $\iota_\Ccal(\alpha)$ and $\iota^R_\Ccal(\alpha)$.
\end{proposition}
\begin{proof}
Recall from (the dual version of) \cite{HTT} proposition 5.3.6.2 that the forgetful functor $\Cat_\pb \rightarrow \Cat$ has a left adjoint $\Fcal$, which maps a category $\Ical$ to the smallest full subcategory of the free completion of $\Ical$ (namely, $\Pcal(\Ical^\op)^\op$) containing $\Ical$ and closed under pullbacks. Consider the pullback preserving functor $\Fcal([1]) \rightarrow \Ccal$ induced from the map $[1] \rightarrow \Ccal$ that sends the unique arrow $a: 0 \rightarrow 1$ in $[1]$ to $\alpha$. The arrows $\iota_\Ccal(\alpha)$ and $\iota_\Ccal^R(\alpha)$ are the images of $\iota_{\Fcal([1])}(a)$ and $\iota_{\Fcal([1])}^R(a)$ under the induced functor $\twoCorr(\Fcal([1])) \rightarrow \twoCorr(\Ccal)$. Moreover, the morphism of spans in the statement is the image of the morphism of spans
\[
\begin{tikzcd}
   & 0\arrow{ldd}{} \arrow[rdd] \arrow[d] &    \\
   & 1 \arrow[ld] \arrow{rd}{}                   &    \\
1 &                                                                   & 1.
\end{tikzcd}
\]
We have thus reduced to proving the result in the case when $\Ccal = \Fcal([1])$ and $\alpha = a: 0 \rightarrow 1$. Note that in $\twoCorr(\Fcal([1]))$ there is a unique map $\iota_{\Fcal([1])}(a) \rightarrow \iota_{\Fcal([1])}(a)$ and a unique map $\iota^R_{\Fcal([1])}(a) \rightarrow \iota^R_{\Fcal([1])}(a)$. Therefore any morphism
\[
\epsilon_a: \id_{\iota_{\Fcal([1])}(0)} \rightarrow \iota_{\Fcal([1])}^R(a) \iota_{\Fcal([1])}(a)
\]
will satisfy the triangle identities with $\eta_a$. Such a morphism is unique, and given by the following diagram:
\[
\begin{tikzcd}
  & 0 \arrow[ldd] \arrow[rdd] \arrow[d] &   \\
  & 0\times_{1} 0 \arrow[ld] \arrow[rd]                 &   \\
0 &                                                      &0
\end{tikzcd}
\] 
\end{proof}

\begin{remark}
Let $\Ccal$ be a category admitting pullbacks and let $\alpha: c \rightarrow c'$ be an arrow in $\Ccal$. The composition $\iota^R_\Ccal(\alpha) \iota_\Ccal(\alpha)$ is given by the span 
$c \leftarrow c\times_d c \rightarrow c$. It follows from the proof of proposition \ref{prop adjunctos in twocrr} that the counit of the adjunction $\iota_\Ccal(\alpha) \dashv \iota_\Ccal^R(\alpha)$ is given by the diagram
\[
\begin{tikzcd}
  & c \arrow[ldd, "\id_c"'] \arrow[rdd, "id_c"] \arrow[d] &   \\
  & c\times_{d} c \arrow[ld] \arrow[rd]                 &   \\
c &                                                      & c
\end{tikzcd}
\]
where the map $c \rightarrow c\times_{d} c$ is the diagonal map.
\end{remark}

\begin{proposition}\label{prop dualizability}
Let $\Ccal$ be a category admitting finite limits and let $c$ be an object in $\Ccal$. Denote by $\Delta: c \rightarrow c \times c$ the diagonal map, and by $\pi: c \rightarrow 1_\Ccal$ the map into the final object of $\Ccal$. Equip $\twoCorr(\Ccal)$ and the inclusion $\iota_\Ccal$ with the symmetric monoidal structures from remark \ref{remark sm structure 2}. Then the morphism 
\[
\eta_c : 1_{\twoCorr(\Ccal)} = \iota_\Ccal(1_\Ccal) \rightarrow \iota_\Ccal(c \times c) = \iota_\Ccal(c) \otimes \iota_\Ccal(c) 
\]
given by the span
\[
\begin{tikzcd}
 & c \arrow{dl}[swap]{\pi} \arrow{dr}{\Delta} & \\ 1_\Ccal & & c\times c
\end{tikzcd}
\]
is the unit of a self duality for $\iota_\Ccal(c)$.
\end{proposition}
\begin{proof}
Let $\Spc_\fin$ be the category of finite spaces. This is obtained from $[0]$ by adjoining finite colimits, and therefore there is a unique left exact functor $F: \Spc_\fin^\op \rightarrow \Ccal$ that maps the point $\ast$ to $c$. As observed in remark \ref{remark sm structure 2}, the functor 
\[
\twoCorr(F): \twoCorr(\Spc_\fin^\op) \rightarrow \twoCorr(\Ccal)
\]
inherits a symmetric monoidal structure, which is compatible the transformations $\iota$ and $\iota^R$. The map $\eta_c$ is the image under $\twoCorr(F)$ of the functor
\[
\eta_\ast : 1_{\twoCorr(\Spc^\op_\fin)} = \iota_{\Spc^\op_\fin}(\emptyset) \rightarrow \iota_{\Spc^\op_\fin}(\ast \amalg \ast) = \iota_{\Spc^\op_\fin}(\ast) \otimes \iota_{\Spc^\op_\fin}(\ast) 
\]
defined by the span
\[
\begin{tikzcd}
 & \ast \arrow{dl}{} \arrow{dr}{} & \\ \emptyset & &  \mathllap{\ast} \amalg \mathrlap{\ast.}
\end{tikzcd}
\]
It therefore suffices to prove the proposition in the case $\Ccal =\Spc_\fin^\op$ and $c = \ast$. Since there is a unique map $\ast \rightarrow \ast$ in $\Spc_\fin^\op$, any map $\epsilon_\ast : \iota_{\Spc^\op_\fin}(\ast) \otimes \iota_{\Spc^\op_\fin}(\ast) \rightarrow 1_{\twoCorr(\Spc^\op_\fin)}$ satisfies the triangle identities with $\eta_\ast$. Such a map exists and is unique, and is given by the following span
\[
\begin{tikzcd}
 & \ast \arrow{dl}{} \arrow{dr}{} & \\ \mathllap{\ast} \amalg \mathrlap{\ast} & & \emptyset.
\end{tikzcd}
\]
\end{proof}

\begin{remark}\label{remark description dualizable cont}
Let $\Ccal$ be a category admitting finite limits and let $c$ be an object in $\Ccal$. Then the proof of proposition \ref{prop dualizability} shows that the morphism
\[
\epsilon_c : \iota_\Ccal(c) \otimes \iota_\Ccal(c) = \iota_\Ccal(c \times c) \rightarrow \iota_\Ccal(1_\Ccal) = 1_{\twoCorr(\Ccal)}
\]
given by the span
\[
\begin{tikzcd}
 & c \arrow{dl}[swap]{\Delta} \arrow{dr}{\pi} & \\ c \times c & &  1_\Ccal
\end{tikzcd}
\]
is the counit of the self duality of proposition \ref{prop dualizability}.
\end{remark}

\begin{proposition} \label{prop description dual map}
Let $\Ccal$ be a category admitting finite limits. Let $c, c'$ be objects of $\Ccal$, and let $\sigma: \iota_\Ccal(c) \rightarrow \iota_\Ccal(c')$ be a morphism between them, represented by a span
\[
\begin{tikzcd}
  & s \arrow[ld, "\alpha"'] \arrow[rd, "\beta"] &    \\
c &                                             & c'.
\end{tikzcd}
\]
  Then the morphism $\sigma^\vee: \iota_\Ccal(c') \rightarrow \iota_{\Ccal}(c)$ dual to $\sigma$ under the self duality of proposition \ref{prop dualizability} is given by the span
  \[
  \begin{tikzcd}
  & s \arrow[ld, "\beta"'] \arrow[rd, "\alpha"] &    \\
c' &                                             & c.
\end{tikzcd}
  \]
\end{proposition}
\begin{proof}
Note that the morphism $\sigma$ is equivalent to $\iota_\Ccal(\beta) \iota_\Ccal^R(\alpha)$, and we have $\sigma^\vee = (\iota_\Ccal^R(\alpha))^\vee \iota_\Ccal(\beta)^\vee$. It suffices therefore to show that there are equivalences $\iota_\Ccal(\beta)^\vee = \iota_\Ccal^R(\beta)$ and $\iota_\Ccal^R(\alpha)^\vee = \iota_\Ccal(\alpha)$. We may furthermore restrict to showing the first identity only - the second one follows from the first one by replacing $\beta$ with $\alpha$ and passing to  adjoints.

Let $\Fcal([1])$ be the free finitely complete category with limits on the arrow category. The morphism $\beta$ is the image of the walking arrow $0 \rightarrow 1$ under a finite limit preserving functor $\Fcal([1]) \rightarrow \Ccal$. It therefore suffices to prove our proposition in the case when $\Ccal = \Fcal([1])$ and $\beta$ is the walking arrow $\beta_{\univ}$. In this case, the dual morphism to $\iota_{\Fcal([1])}(\beta_{\univ})$ can be computed as the following composition:
\[
\begin{tikzcd}
  & 1 \times 0 \arrow[ld] \arrow[rd, "\id_1 \times \Delta_0"] &                     & 1 \times 0 \times 0 \arrow[ld, "\id"'] \arrow[rd, "\id_1 \times _{\beta_\univ} \times \id_0"] &                     & 1 \times 0 \arrow[ld, "\Delta_1 \times \id_0"] \arrow[rd] &   \\
1 &                                  & 1 \times 0 \times 0 &                                           & 1 \times 1 \times 0 &                                  & 0.
\end{tikzcd}
\]
This recovers the unique span of the form
\[
\begin{tikzcd}
  & 0 \arrow[ld] \arrow[rd] &   \\
1 &                         & 0
\end{tikzcd}
\]
which is $\iota^R_{\Fcal([1])}(\beta_{\univ})$, as desired.
\end{proof}

\subsection{Beck-Chevalley conditions}\label{subsection BC}

We now discuss the universal property of the $2$-category of correspondences.

\begin{definition}\label{definition adjointable}
Let $\Dcal$ be a $2$-category. We say that a commutative diagram
\[
\begin{tikzcd}
d' \arrow{r}{\alpha'} \arrow{d}{\beta'} & d \arrow{d}{\beta} \\e' \arrow{r}{\alpha} & e
\end{tikzcd}
\]
in $\Dcal$ is vertically right adjointable if the following conditions hold:
\begin{itemize}
\item The maps $\beta$ and $\beta'$ admit right adjoints $\beta^R$ and $\beta'^R$.
\item The $2$-cell
\[
\alpha' \beta'^R \rightarrow \beta^R \beta \alpha' \beta'^R = \beta^R \alpha \beta' \beta'^R \rightarrow \beta^R \alpha
\]
built from the unit $\id_d \rightarrow \beta^R \beta$ and the counit $\beta' \beta^R \rightarrow \id_{e'}$, is  an isomorphism.
\end{itemize}

We say that the above diagram is horizontally right adjointable if its transpose is vertically
right adjointable. We say that it is right adjointable if it is both horizontally and vertically
right adjointable. We say that it is (vertically / horizontally) left adjointable if it is (vertically / horizontally) right adjointable as a diagram in the $2$-category $\Dcal^{2\dsh\op}$ obtained from $\Dcal$ by reversing the direction of the $2$-cells.
\end{definition}

\begin{remark}\label{remark adjointability square lax}
Let $\Dcal$ be a $2$-category. A commutative square
\[
\begin{tikzcd}
d' \arrow{r}{\alpha'} \arrow{d}{\beta'} & d \arrow{d}{\beta} \\e' \arrow{r}{\alpha} & e
\end{tikzcd}
\]
in $\Dcal$ defines a morphism $\gamma: \alpha' \rightarrow \alpha$ in the arrow category $\Funct([1],\Dcal)$. It is shown in \cite{Haugsadj} theorem 4.6 that $\beta'$ and $\beta$ admit right adjoints if and only if $\gamma$ admits a right adjoint in the category $\Funct([1], \Dcal)_{\text{lax}}$ of functors $[1] \rightarrow \Dcal$ and lax natural transformations. In that case, the right adjoint to $\gamma$ is given by the lax commutative square 
\[
\begin{tikzcd}
d' \ar[dr,shorten <>=10pt,Rightarrow] \arrow{r}{ \alpha'} & d \\
e'   \arrow{u}{\beta'^R} \arrow{r}{\alpha} & e \arrow{u}{\beta^R}
\end{tikzcd}
\]
where the $2$-cell is the one described in definition \ref{definition adjointable}. It follows that $\gamma$ admits a right adjoint in $\Funct([1], \Dcal) \subset \Funct([1], \Dcal)_{\text{lax}}$ if and only if the original commutative square is vertically right adjointable. In this case, the morphism $\gamma^R$ corresponds to a commutative square in $\Dcal$; we say that this square arises from the original one by passage to right adjoints of vertical arrows.
\end{remark}

\begin{remark}\label{lemma left adj iff right adj}
Let $\Dcal$ be a $2$-category and let
\[
\begin{tikzcd}
d' \arrow{r}{\alpha'} \arrow{d}{\beta'} & d \arrow{d}{\beta} \\e' \arrow{r}{\alpha} & e
\end{tikzcd}
\]
be a commutative diagram in $\Dcal$. Assume that the vertical maps are right adjointable, and that the horizontal maps are left adjointable. We can then construct two lax commutative squares
\[
\begin{tikzcd}
d' \ar[dr,shorten <>=10pt,Rightarrow] \arrow{r}{ \alpha'} & d \\
e'   \arrow{u}{\beta'^R} \arrow{r}{\alpha} & e \arrow{u}{\beta^R}
\end{tikzcd} \hspace{2cm} 
\begin{tikzcd}
d'  \arrow{d}{\beta'}  & d \arrow{l}[swap]{\alpha'^L} \arrow{d}{\beta} \\
e'     &\ar[ul,shorten <>=10pt,Rightarrow]  e \arrow{l}[swap]{\alpha^L}.
\end{tikzcd}
\]
These two are related to each other by passage to left/right adjoints.  In particular we have that our original square is vertically right adjointable if and only if it is horizontally left adjointable.
\end{remark}

\begin{remark}\label{right adjoint in limit}
Denote by $\Adj$ the universal adjunction. This is a $2$-category equipped with an epimorphism $L: [1] \rightarrow \Adj$  such that for every $2$-category $\Dcal$, precomposition with $L$ induces an equivalence between the space of functors $\Adj \rightarrow \Dcal$ and the space of maps $[1] \rightarrow \Dcal$ which pick out a right adjointable arrow in $\Dcal$.  Let $\Ucal_{\text{lax}}$ be the pushout in $\twoCat$ of the following diagram
\[
\begin{tikzcd}
\Adj &                                                                          & {[1] \times [1]} &                                                                          & \Adj \\
     & {[1]} \arrow[lu, "L"] \arrow[ru, "{\lbrace 0 \rbrace \times \id_{[1]}}"] &                  & {[1]} \arrow[lu, swap, "{\lbrace 1 \rbrace \times \id_{[1]}}"] \arrow[ru, "L"] &   .  
\end{tikzcd}
\]
This is the universal lax vertically right adjointable square. It contains a $2$-cell $\mu_{\univ}$ which is the universal instance of the $2$-cell from definition \ref{definition adjointable}. 

Let $\Ucal$ be the $2$-category obtained from $\Ucal_{\text{lax}}$ by inverting the $2$-cell $\mu_{\univ}$. The natural inclusion $i: [1] \times [1] \rightarrow \Ucal$ is an epimorphism in $\twoCat$. For any $2$-category $\Dcal$, precomposition with $i$ induces an equivalence between the space of functors $\Ucal \rightarrow \Dcal$ and the space of functors $[1] \times [1] \rightarrow \Dcal$ which represent a vertically right adjointable square. By virtue of remark \ref{remark adjointability square lax}, we in fact have an equivalence between $i: [1] \times [1] \rightarrow \Ucal$ and $L \times \id_{[1]} : [1] \times [1] \rightarrow  \Adj \times [1]$.

As a consequence of the above we deduce the following fact: if $F: \Ical \rightarrow \twoCat$ is a diagram in $\twoCat$ with limit $\Dcal$, then a commutative square in $\Dcal$ is vertically right adjointable if and only if its image in $F(i)$ is vertically right adjointable, for all $i$ in $\Ical$.
\end{remark}

\begin{definition}
Let $\Ccal$ be a category admitting pullbacks and $\Dcal$ be a $2$-category. We say that a functor $F:\Ccal \rightarrow \Dcal$ satisfies the left Beck-Chevalley condition if for every cospan $x \rightarrow s \leftarrow y$ in $\Ccal$, the induced commutative square in $\Dcal$
\[
\begin{tikzcd}
F(x \times_s y) \arrow{r}{} \arrow{d}{} & F(y) \arrow{d}{} \\ F(x) \arrow{r}{} & F(s)
\end{tikzcd}
\]
is right adjointable.
\end{definition}

\begin{proposition}\label{prop twocorr satisface bc}
Let $\Ccal$ be a category admitting pullbacks. Then the inclusion $\iota_\Ccal: \Ccal \rightarrow \twoCorr(\Ccal)$ satisfies the left Beck-Chevalley condition.
\end{proposition}
\begin{proof}
Let $\Fcal: \Cat \rightarrow \Cat_\pb$ be the left adjoint to the forgetful functor. Let $\Ucal$ be universal cospan: this is the category with three objects $0, 1, 2$ and nontrivial arrows $\alpha: 0 \rightarrow 1 \leftarrow 2: \beta$. Any cospan in $\Ccal$ is the image of a pullback preserving morphism $\Fcal(\Ucal) \rightarrow \Ccal$. It therefore suffices to prove that the image under $\iota_{\Fcal(\Ucal)}$ of the universal cartesian square
\[
\begin{tikzcd}
0 \times_1 2 \arrow{r}{\alpha'} \arrow{d}{\beta'} & 2 \arrow{d}{\beta} \\ 0 \arrow{r}{\alpha} & 1
\end{tikzcd}
\]
is vertically right adjointable. Recall from proposition \ref{prop adjunctos in twocrr} that $\iota_{\Fcal(\Ucal)}(\beta')$ is left adjoint to $\iota^R_{\Fcal(\Ucal)}(\beta')$, and $\iota_{\Fcal(\Ucal)}(\beta)$ is left adjoint to $\iota^R_{\Fcal(\Ucal)}(\beta)$. We have that $\iota_{\Fcal(\Ucal)}(\alpha')\iota^R_{\Fcal(\Ucal)}(\beta')$ and $\iota^R_{\Fcal(\Ucal)}(\beta)\iota_{\Fcal(\Ucal)}(\alpha)$ are both given by the span
\[
\begin{tikzcd}
& 0\times_1 2 \arrow{dl}[swap]{\beta'} \arrow{dr}{\alpha'} & \\ 0  & & 2.
\end{tikzcd}
\]
Our claim now follows from the fact that $\Hom_{\Fcal(\Ucal)}(0 \times_1 2, 0\times_1 2) = [0]$, and thus any $2$-cell 
\[
\iota_{\Fcal(\Ucal)}(\alpha')\iota^R_{\Fcal(\Ucal)}(\beta') \rightarrow  \iota^R_{\Fcal(\Ucal)}(\beta)\iota_{\Fcal(\Ucal)}(\alpha)
\]
is necessarily invertible.
\end{proof}

The rest of this section is devoted to showing that $\twoCorr(\Ccal)$ is the universal $2$-category equipped with a functor from $\Ccal$ which satisfies the left Beck-Chevalley condition (theorem \ref{univ prop twocorr}). A proof was given in \cite{GR} chapter 7 theorem 3.2.2. - here we present an alternative approach using the theory of two-sided fibrations, and in particular the universal property of the span fibration established in proposition \ref{prop univer span}. The proof will need a few preliminary lemmas.

\begin{notation}
Let $\Ccal$ be a category. Recall the universal left adjointable arrow $L: [1] \rightarrow \Adj$ from remark \ref{right adjoint in limit}. We let $\Ccal^R$ be the $2$-category defined by the pushout
\[
\begin{tikzcd}
\Hom_{\Cat}([1], \Ccal) \times [1] \arrow{r}{\ev} \arrow{d}{\id \times L} & \Ccal  \arrow{d}{L_\Ccal}  \\ \Hom_{\Cat}([1],\Ccal) \times \Adj \arrow{r}{} & \Ccal^{R}.
\end{tikzcd}
\]
The functor $L_\Ccal$ is an epimorphism in $\twoCat$. For every $2$-category $\Dcal$, precomposition with $L_\Ccal$ induces an equivalence between the space $\Hom_{\twoCat}(\Ccal^{R} , \Dcal)$ and the space of functors $F: \Ccal \rightarrow \Dcal$ such that for every arrow $\alpha$ in $\Ccal$ the arrow $F(\alpha)$ admits a right adjoint in $\Dcal$.
\end{notation}

\begin{lemma}\label{lemma adjointability in funct}
Let $\Ccal$ be a category and let $\Dcal$ be a $2$-category. Let $\eta: F \rightarrow G$ be a morphism in $\Funct(\Ccal,\Dcal)$. Then the morphism $\eta$ is left adjointable if and only if for every morphism $\alpha: x \rightarrow y$ in $\Ccal$, the commutative square
\[
\begin{tikzcd}
F(x) \arrow{r}{\eta(x)} \arrow{d}{F(\alpha)} & G(x) \arrow{d}{G(\alpha)} \\ F(y) \arrow{r}{\eta(y)} & G(y)
\end{tikzcd}
\]
is horizontally left adjointable.
\end{lemma}
\begin{proof}
Let $\Scal$ be the full subcategory of $\Cat$ on those categories $\Ccal$ for which the lemma holds. As discussed in remark \ref{remark adjointability square lax}, the walking arrow belongs to $\Scal$. To prove this lemma it suffices to show that $\Scal$ is closed under colimits in $\Cat$. Assume given a diagram $C: \Ical \rightarrow \Cat$ with $C(i)$ in $\Scal$ for every $i$ in $\Ical$. Then we have an equivalence
\[
\Funct(\colim_{\Ical} C(i), \Dcal) = \lim_{\Ical^\op} \Funct(C(i), \Dcal).
\]
A morphism $\eta: F \rightarrow G$ in $\Funct(\colim_{\Ical} C(i), \Dcal)$ is left adjointable if and only if its image in $\Funct(C(i),\Dcal)$ is left adjointable for every $i$ in $\Ical$. This happens if and only if the square
\[
\begin{tikzcd}
F(x) \arrow{r}{\eta(x)} \arrow{d}{F(\alpha)} & G(x) \arrow{d}{G(\alpha)} \\ F(y) \arrow{r}{\eta(y)} & G(y)
\end{tikzcd}
\]
is horizontally left adjointable for every arrow $\alpha: x \rightarrow y$ in $\colim_{\Ical} C(i)$ which belongs to the image of the map $C(i) \rightarrow \colim_{\Ical} C(i)$ for some $i$ in $\Ical$. Observe now that the family of arrows $\alpha$ in $\colim_{\Ical}C(i)$ for which the above square is horizontally left adjointable is closed under compositions (this follows from example from the characterization of adjointability of squares from remark \ref{remark adjointability square lax}). We conclude that $\eta$ is left adjointable if and only if the above square is horizontally left adjointable for every arrow, which implies that $\colim_{\Ical}C(i)$ also belongs to $\Scal$, as desired.
\end{proof}

\begin{lemma}\label{lemma identify functs out of Cr}
Let $\Ccal$ be a category. Then the composition of the functor 
\[
(L_\Ccal^{1\dsh\op} \times \id_\Ccal)^* : \Funct((\Ccal^R)^{1\dsh\op} \times  \Ccal, \Catscr)^{\leq 1} \rightarrow  \Funct(\Ccal^\op \times \Ccal, \Cat)
\]
and the two-sided Grothendieck construction \footnote{In this section we use the convention where two-sided fibrations are cartesian over the first coordinate and cocartesian over the second coordinate, unless otherwise stated. Note that this differs from the convention used in section $2$.}
\[
\int_{\Ccal^\op \times \Ccal} : \Funct(\Ccal^\op \times \Ccal, \Cat) \rightarrow \Cat_{/\Ccal \times \Ccal}
\]
is a monomorphism, whose image is the category $\Cat^{\cocart,\bicart}_{/\Ccal \times \Ccal} = \Cat^{\cocart}_{/\Ccal \times \Ccal} \cap \Cat^{\bicart}_{/\Ccal \times \Ccal}$.
\end{lemma}
\begin{proof}
We note that since $L_\Ccal$ is an epimorphism, we have that $L_\Ccal^{1\dsh\op} \times \id_\Ccal$ is also an epimorphism, and hence the induced functor
\[
\Funct((\Ccal^R)^{1\dsh\op} \times \Ccal ,\Catscr)^{\leq 1} \rightarrow \Funct(\Ccal^{\op} \times \Ccal , \Catscr)^{\leq 1} = \Funct(\Ccal^{\op} \times \Ccal , \Cat)
\]
is a indeed a monomorphism. We note that the above map is equivalent to the  map
\[
(L_\Ccal^{1\dsh\op})^*: \Funct((\Ccal^R)^{1\dsh\op}, \Funct(\Ccal, \Catscr))^{\leq 1} \rightarrow \Funct(\Ccal^\op, \Funct(\Ccal, \Cat) ).
\]

It follows from lemma \ref{lemma adjointability in funct} that an object belongs to the image of $(L_\Ccal^{1\dsh\op})^*$ if and only if the associated functor $\Ccal^\op \times \Ccal \rightarrow \Cat$ is left adjointable in the $\Ccal^\op$ coordinate. 
This happens if and only if the associated two-sided fibration is also a  cocartesian fibration. We conclude that the lemma holds at the level of objects.

Consider now the functor
\[
\Funct((\Ccal^R)^{1\dsh\op}, \Funct(\Ccal \times [1],\Catscr))^{\leq 1} \rightarrow \Funct(\Ccal^\op, \Funct(\Ccal \times [1],\Cat))
\]
of precomposition with $L_\Ccal^{1\dsh\op}$. Applying lemma \ref{lemma adjointability in funct} again we conclude that an object belongs to its image if and only if the associated functor $\Ccal^\op \times \Ccal \times [1] \rightarrow \Cat$ is left adjointable in the $\Ccal^\op$ coordinate. Applying remark \ref{remark  adjointability triple product} we see that this happens if and only if the associated morphism in $\Cat^{\bicart}_{/\Ccal \times \Ccal}$ belongs also to $\Cat^{\cocart}_{/\Ccal \times \Ccal}$. This shows that the lemma holds also at the level of morphisms, which finishes the proof.
\end{proof}

\begin{notation}
Let $\Ccal$ be a category admitting pullbacks. Recall the universal vertically right adjointable square $\id_{[1]} \times L : [1] \times [1] \rightarrow [1] \times \Adj$ from remark \ref{right adjoint in limit}. Let $S$ be the space of cartesian commutative squares in $\Ccal$, and let $\ev: S \times ([1] \times [1]) \rightarrow \Ccal$ be the evaluation functor. Let $\twoCorr^{\univ}(\Ccal)$ be the $2$-category defined by the pushout
\[
\begin{tikzcd}
S \times ([1] \times [1]) \arrow{r}{\ev} \arrow{d}{\id_S \times (\id_{[1]} \times L)} & \Ccal  \arrow{d}{\iota_\Ccal^{\univ}}  \\ S \times ([1] \times \Adj) \arrow{r}{} & \twoCorr^{\univ}(\Ccal)
\end{tikzcd}.
\]
The $2$-category $\twoCorr^{\univ}(\Ccal)$ is the universal $2$-category equipped with the a functor from $\Ccal$ which satisfies the left Beck-Chevalley condition. The map $\iota_\Ccal^{\univ}$ is an epimorphism in $\twoCat$: for every $2$-category $\Dcal$, precomposition with $\iota_\Ccal^{\univ}$ induces an equivalence between the space of functors $\twoCorr^{\univ}(\Ccal) \rightarrow \Dcal$ and the space of functors $\Ccal \rightarrow \Dcal$ which satisfy the left Beck-Chevalley condition. It follows in particular that $\iota_\Ccal^\univ$ factors through $\Ccal^R$. We denote by $q: \Ccal^R \rightarrow \twoCorr^\univ(\Ccal)$ the induced functor.
\end{notation}

\begin{lemma}\label{lemma compare funct Cr with funct corr}
Let $\Ccal$ be a category admitting pullbacks and let $\Dcal$ be a $2$-category. The functor
\[
q^*: \Funct(\twoCorr^\univ(\Ccal),\Dcal) \rightarrow \Funct(\Ccal^R, \Dcal)
\]
induces equivalences at the level of Hom categories.
\end{lemma}
\begin{proof}
Observe that $\twoCorr^\univ(\Ccal)$ is obtained out of $\Ccal^R$ by inverting a family of $2$-cells. It therefore suffices to show that if $q': C_2 \rightarrow [1]$ is the projection from the walking $2$-cell to the walking $1$-cell, the induced functor
\[
q'^* : \Funct([1],\Dcal) \rightarrow \Funct(C_2, \Dcal)
\]
induces equivalences at the level of Hom categories. 

Let $\alpha, \beta: [1] \rightarrow \Dcal$ be a pair of objects of $\Funct([1], \Dcal)$. We have a commutative cube of categories as follows:
\[
    \begin{tikzcd}[row sep=1.5em, column sep = 0.25em]
        &  \Hom_{\Funct(C_2, \Dcal)}(q'^*\alpha, q'^*\beta)\arrow[dd, "\ev_1", pos=0.75] \arrow[rr, "\ev_0"] &&
    \Hom_\Dcal(\alpha(0), \beta(0)) \arrow[dd,"q'^*\beta"] \\ \Hom_{\Funct([1],\Dcal)}(\alpha, \beta) 
    \arrow[rr,"\ev_0", pos=0.7] \arrow[ur, "q'^*"] \arrow[dd, "\ev_1"] &&
    \Hom_{\Dcal}(\alpha(0), \beta(0)) \arrow[dd, "\beta" ,pos = 0.25] \arrow[ur,"\id"] \\
 &\Hom_\Dcal(\alpha(1), \beta(1)) \arrow[rr, "q'^*\alpha", pos=0.25] && \Funct([1], \Hom_{\Dcal}(\alpha(0), \beta(1))) \\
  \Hom_\Dcal(\alpha(1), \beta(1))\arrow[rr, "\alpha"] \arrow[ur, "\id"] && \Hom_{\Dcal}(\alpha(0), \beta(1)) \arrow[ur, ""]
    \end{tikzcd}
\]
Here the front and back faces are cartesian. The bottom left and top right diagonal arrows are isomorphisms, and the bottom right diagonal arrow is the degeneracy map, which is fully faithful. It follows that the top left diagonal arrow is an isomorphism, as desired.
\end{proof}

\begin{lemma}\label{lemma check invertibility 2 cell in funct}
Let $\Jcal$ and $\Dcal$ be $2$-categories. Then a $2$-cell in $\Funct(\Jcal, \Dcal)$ is invertible if and only if its image under all evaluation functors is invertible. 
\end{lemma}
\begin{proof}

Let $\Scal$ be the full subcategory of $\twoCat$ on those $2$-categories $\Jcal$ for which the lemma holds. We claim that $\Scal$ is closed under colimits in $\twoCat$. Assume given a  diagram $J: \Ical \rightarrow \twoCat$, with $J(i)$ in $\Scal$ for every $i$ in $\Ical$. Then we have an equivalence 
\[
\Funct(\colim_\Ical J(i) , \Dcal) = \lim_{\Ical^\op}\Funct(J(i), \Dcal).
\]
It follows that a $2$-cell $\gamma: C_2 \rightarrow \Funct(\colim_{\Ical} J(i), \Dcal)$ is invertible if and only if its image in $\Funct(J(i),\Dcal)$ is invertible for all $i$ in $\Ical$. Since $J(i)$ is assumed to belong to $\Scal$, this happens if and only if  the $2$-cell $\ev_j \gamma$ is invertible for all objects $j$ in $\colim_{\Ical} J(i)$, which means that $\colim_{\Ical}J(i)$ also belongs to $\Scal$.

To prove the lemma it then suffices to show that $\Scal$ contains the walking $2$-cell $C_2$. In other words, we have reduced to proving the lemma in the case $\Jcal = C_2$. In this case, for every pair of objects $\mu, \nu$ in $\Funct(C_2, \Dcal)$, we have a pullback diagram
\[
\begin{tikzcd}
\Hom_{\Funct(C_2, \Dcal)}(\mu, \nu) \arrow{r}{} \arrow{d}{} & \Hom_\Dcal(\mu(0), \nu(0)) \arrow{d}{} \\
\Hom_\Dcal(\mu(1), \nu(1)) \arrow{r}{} & \Funct([1], \Hom_\Dcal(\mu(0), \nu(1)))
\end{tikzcd}
\]
where the top and left arrows are the evaluation functors. A $2$-cell in $\Funct(C_2, \Dcal)$ corresponds to an arrow in $\Hom_{\Funct(C_2, \Dcal)}(\mu, \nu) $, and this is invertible if and only if its image in $ \Hom_\Dcal(\mu(0), \nu(0))$ and $\Hom_\Dcal(\mu(1), \nu(1)) $ is invertible, as desired.
\end{proof}

\begin{lemma}\label{lemma check bc a functor cats}
Let $\Jcal, \Dcal$ be $2$-categories, and let $\Ccal$ be a category admitting pullbacks. Then a functor $F: \Ccal^R \rightarrow \Funct(\Jcal, \Dcal)$ extends to $\twoCorr^\univ(\Ccal)$ if and only if for every object $j$ in $\Jcal$ the composite functor $\ev_j F : \Ccal^R \rightarrow \Dcal$ extends to $\twoCorr^\univ(\Ccal)$.
\end{lemma}
\begin{proof}
This follows directly from lemma \ref{lemma check invertibility 2 cell in funct} using the fact that  $\twoCorr^\univ(\Ccal)$ is obtained from $\Ccal^R$ by inverting a family of $2$-cells.
\end{proof}

\begin{lemma}\label{lemma twocorruniv y bc}
Let $\Ccal$ be a category admitting pullbacks. Then the composition of the functor 
\[
((\iota_{\Ccal}^{\univ})^{1\dsh\op} \times \id_\Ccal)^* : \Funct(\twoCorr^{\univ}(\Ccal)^{1\dsh\op} \times  \Ccal, \Catscr)^{\leq 1} \rightarrow  \Funct(\Ccal^\op \times \Ccal, \Cat)
\]
and the two-sided Grothendieck construction
\[
\int_{\Ccal^\op \times \Ccal} : \Funct(\Ccal^\op \times \Ccal, \Cat) \rightarrow \Cat_{/\Ccal \times \Ccal}
\]
is a monomorphism, whose image is the full subcategory of $\Cat^{\cocart, \bicart}_{/\Ccal \times \Ccal}$  on those fibrations which satisfy the Beck-Chevalley condition in the first coordinate.
\end{lemma}
\begin{proof}
The functor $q: \Ccal^R \rightarrow \twoCorr^{\univ}(\Ccal)$ is an epimorphism, so we have a monomorphism
\[
(q^{1\dsh\op} \times \id_\Ccal)^* :  \Funct(\twoCorr^{\univ}(\Ccal)^{1\dsh\op} \times  \Ccal, \Catscr)^{\leq 1}    \rightarrow        \Funct((\Ccal^R)^{1\dsh\op} \times  \Ccal, \Catscr)^{\leq 1} .
\]
The above is equivalent to the functor
\[
\Funct(\twoCorr^\univ(\Ccal)^{1\dsh\op}, \Funct(\Ccal, \Catscr))^{\leq 1} 
\rightarrow \Funct((\Ccal^R)^{1\dsh\op}, \Funct(\Ccal, \Catscr))^{\leq 1}
\]
obtained by precomposition with  $q^{1\dsh\op}$. Applying lemma \ref{lemma compare funct Cr with funct corr} we see that the above functor is fully faithful. Combining this fact with lemma \ref{lemma identify functs out of Cr} we conclude that  precomposition with $(\iota_\Ccal^\univ)^{1\dsh\op} \times \id_\Ccal$ followed by the two-sided Grothendieck construction embeds
\[
\Funct(\twoCorr^{\univ}(\Ccal)^{1\dsh\op} \times \Ccal, \Catscr)^{\leq 1}
\]
as a full subcategory of $\Cat^{\cocart,\text{two-sided}}_{/{\Ccal \times \Ccal}}$. Thanks to lemma \ref{lemma check bc a functor cats}, its image consists of those fibrations which satisfy the Beck-Chevalley condition in the first variable, as desired.
\end{proof}

The next three lemmas use the theory of algebroids and enriched categories as developed in \cite{GH} and \cite{Hinich}, and in particular the approach to the Yoneda embedding via diagonal bimodules from \cite{Hinich}. We refer the reader to \cite{Pres} for our conventions regarding this subject.

\begin{lemma}\label{lemma map of bimodules cocartesian}
Let $\Mcal$ be a presentable monoidal category and let $f:\Acal \rightarrow \Bcal$ be a morphism of $\Mcal$-algebroids which is an equivalence on categories of objects. Let ${}_\Acal \Acal_\Acal$ be the diagonal $\Acal$-bimodule, and let ${}_\Acal \Bcal_\Bcal$ be the restriction of scalars of the diagonal $\Bcal$-bimodule along the map $(f, \id_{\Bcal})$. Then the induced morphism  $f_* : {}_{\Acal} \Acal_\Acal \rightarrow {}_{\Acal} \Bcal_\Bcal$  in $\operatorname{BMod}(\Mcal)$ is a cocartesian lift of the morphism $(\id_\Acal, f): (\Acal, \Acal) \rightarrow (\Acal, \Bcal)$ in $\Algbrd(\Mcal) \times \Algbrd(\Mcal)$. 
\end{lemma}
\begin{proof}
Let $X$ be the category of objects of $\Acal$ and $\Bcal$. By construction, the projection $p:\operatorname{BMod}(\Mcal) \rightarrow \Algbrd(\Mcal) \times \Algbrd(\Mcal)$ is a morphism of cartesian fibrations over $\Cat \times \Cat$. Thanks to \cite{HTT} corollary 4.3.1.15, to show that $f_*$ is $p$-cocartesian it suffices to show that it is cocartesian for the projection 
\[
p_{X, X}: \operatorname{BMod}_{X, X}(\Mcal) \rightarrow \Algbrd_X(\Mcal) \times \Algbrd_X(\Mcal).
\] 

Recall from \cite{Hinich} proposition 3.3.6 that $\Assos_X$ is a flat associative operad. It follows that there is an associative operad $\Mcal_X$ equipped with the universal map of associative operads $\Assos_X \times \Mcal_X \rightarrow \Mcal$. As discussed in \cite{Hinich} corollary 4.4.9, the associative operad $\Mcal_X$ is a presentable monoidal category. The projection $p_{X, X}$ is equivalent to the projection
\[
p'_{X, X}: \Alg_{\BM}(\Mcal_X) \rightarrow \Alg_{\Assos}(\Mcal_X) \times \Alg_{\Assos}(\Mcal_X)
\]
which sends each $\BM$-algebra in $\Mcal_X$ to its underlying associative algebras.

Let $\widetilde{\Acal}$ and $\widetilde{\Bcal}$ be the associative algebras corresponding to the $\Mcal$-algebroids $\Acal$ and $\Bcal$,respectively. Then under the above equivalence, the map $f_*$ corresponds to the map of bimodules $\widetilde{f}_*: {}_{\widetilde{\Acal}} \widetilde{\Acal}_{\widetilde{\Acal}} \rightarrow {}_{\widetilde{\Acal}} \widetilde{\Bcal}_{\widetilde{\Bcal}}$, which is $p'_{X, X}$-cocartesian, as desired. 
\end{proof}

\begin{lemma}\label{lemma homs y epis}
Let $\Mcal$ be a presentable symmetric monoidal category and let $f: \Acal \rightarrow \Bcal$ be a morphism of $\Mcal$-algebroids which is an equivalence on categories of objects. Consider the morphisms 
\[
\Hcal_\Acal = \Acal(-,-): \Acal  \otimes \Acal^\op \rightarrow \overline{\Mcal}
\]
 and 
 \[
 \Hcal_{\Bcal} = \Bcal(-,-): \Bcal \otimes \Bcal^\op \rightarrow \overline{\Mcal}
 \]
 where $\overline{\Mcal}$ denotes the enhancement of $\Mcal$ to an $\Mcal$-enriched category. Let $f_* : \Hcal_{\Acal} \rightarrow \Hcal_{\Bcal}|_{\Acal \otimes \Acal^\op} $
be the induced map. Then for every object $G$ in $\Funct(\Acal \otimes \Bcal^\op, \overline{\Mcal})$, the composition of the restriction map
\[
\tau_\Mcal\Hom_{\Funct(\Acal \otimes \Bcal^\op, \overline{\Mcal})}(\Hcal_\Bcal|_{\Acal \otimes \Bcal^\op}, G) \rightarrow \tau_\Mcal\Hom_{\Funct(\Acal \otimes \Acal^\op, \overline{\Mcal})}(\Hcal_{\Bcal}|_{\Acal \otimes \Acal^\op}, G|_{\Acal \otimes \Acal^\op})
\]
with the precomposition with $f_*$ map
\[
\tau_\Mcal\Hom_{\Funct(\Acal \otimes \Acal^\op, \overline{\Mcal})}(\Hcal_{\Bcal}|_{\Acal \otimes \Acal^\op}, G|_{\Acal \otimes \Acal^\op}) \rightarrow \tau_\Mcal\Hom_{\Funct(\Acal \otimes \Acal^\op, \overline{\Mcal})}(\Hcal_{\Acal}, G|_{\Acal \otimes \Acal^\op})
\]
is an equivalence.
\end{lemma}
\begin{proof}
This is a translation of lemma \ref{lemma map of bimodules cocartesian} under the folding equivalence of \cite{Hinich} section 3.6, and the equivalence of \cite{Hinich} proposition 6.3.7.
\end{proof}

\begin{lemma}\label{lemma identification hom}
Let $\Ccal$ be a category admitting pullbacks. 
\begin{enumerate}[\normalfont (i)]
\item The map $\iota_\Ccal^\univ$ is surjective on objects.
\item The image of the natural transformation
\[
(\iota_\Ccal^{\univ})_*: \Hom_\Ccal(-, -) \rightarrow \Hom_{\twoCorr^\univ(\Ccal)}(-,-)|_{\Ccal^\op \times \Ccal}
\]
under the two-sided Grothendieck construction
\[
\int_{\Ccal^\op \times \Ccal} : \Funct(\Ccal^\op \times \Ccal, \Cat) \rightarrow \Cat_{/\Ccal \times \Ccal}
\]
is equivalent to the morphism of two-sided fibrations over $\Ccal \times \Ccal$
\[
\begin{tikzcd}
 \Funct([1], \Ccal)  \arrow{rr}{\phi'} \arrow{dr}[swap]{(\ev_0, \ev_1)} & & \arrow{dl}{(\ev_1, \ev_2)}  \Funct(\Lambda^2_0, \Ccal) \\ & \Ccal \times \Ccal & 
\end{tikzcd}
\]
where $\Lambda^2_0$ is the category with objects $0, 1, 2$ and nontrivial morphisms $1 \leftarrow 0 \rightarrow 2$, and $\phi'$ is the functor of precomposition with the map $\Lambda^2_0 \rightarrow [1]$ which sends $0, 1, 2$ to $0, 0, 1$, respectively.
\end{enumerate}
\end{lemma}
\begin{proof}
Let $S$ be the space of cartesian commutative squares in $\Ccal$, and let $\ev: S \times ([1] \times [1]) \rightarrow \Ccal$ be the evaluation functor. Let $\twoCorr^{\univ}_{\text{algbrd}}(\Ccal)$ be the $\Cat$-algebroid defined as the pushout
\[
\begin{tikzcd}
S \times ([1] \times [1]) \arrow{r}{\ev} \arrow{d}{\id_S \times (\id_{[1]} \times L)} & \Ccal  \arrow{d}{\iota_{\Ccal, \text{algbrd}}^{\univ}}  \\ S \times ([1] \times \Adj) \arrow{r}{} & \twoCorr^{\univ}_{\text{algbrd}}(\Ccal)
\end{tikzcd}
\]
inside $\Algbrd(\Cat)$. Note that the image of $\iota^\univ_{\Ccal, \text{algbrd}}$ under the localization functor $\Algbrd(\Cat) \rightarrow \twoCat$ recovers  $\iota_\Ccal^\univ$.   It follows from \cite{Pres} proposition 2.3.12 that the map of algebroids $\iota^{\univ}_{\Ccal, \text{algbrd}}$ is an equivalence at the level of objects. Item (i) now follows from this together with the fact that the map from any algebroid to its completion is surjective on objects.

Denote by $j: \twoCorr^\univ_\algbrd(\Ccal) \rightarrow \twoCorr^\univ(\Ccal)$ the canonical map. As shown in \cite{GH} corollary 5.6.3, the map $j$ is fully faithful. In other words, the natural transformation
\[
j_*: \Hom_{\twoCorr^{\univ}_\algbrd(\Ccal)}(- , -) \rightarrow \Hom_{\twoCorr^\univ(\Ccal)}(-, -)|_{ \twoCorr^{\univ}_\algbrd(\Ccal)^{1\dsh\op} \times \twoCorr^{\univ}_\algbrd(\Ccal)}
\]
is an equivalence.  It follows that the image of  the transformation
\[
(\iota_\Ccal^{\univ})_*: \Hom_\Ccal(-, -) \rightarrow \Hom_{\twoCorr^\univ(\Ccal)}(-,-)|_{\Ccal^\op \times \Ccal}
\]
under the two-sided Grothendieck construction is equivalent to the image of 
\[
(\iota_{\Ccal,\algbrd}^{\univ})_*: \Hom_\Ccal(-, -) \rightarrow \Hom_{\twoCorr^\univ_{\algbrd}(\Ccal)}(-,-)|_{\Ccal^\op \times \Ccal}.
\]

  The fact that  $\Catscr$ is complete as a $\Cat$-algebroid implies that the restriction map
\[
\Funct(\twoCorr^{\univ}(\Ccal)^{1\dsh\op} \times \Ccal, \Catscr) \rightarrow  \Funct(\twoCorr^{\univ}_\algbrd(\Ccal)^{1\dsh\op} \times \Ccal, \Catscr)
\]
is an equivalence, and therefore the conclusion of lemma \ref{lemma twocorruniv y bc} remains valid if we replace $\twoCorr^\univ(\Ccal)$ with $\twoCorr^\univ_\algbrd(\Ccal)$. Item (ii) now follows this together with a combination of lemma \ref{lemma homs y epis} and proposition \ref{prop univer span}.
\end{proof} 

\begin{theorem}\label{univ prop twocorr}
Let $\Ccal$ be a category admitting pullbacks, and $\Dcal$ be a $2$-category. Precomposition with the functor $\iota_\Ccal: \Ccal \rightarrow \twoCorr(\Ccal)$ induces an equivalence between the space $\Hom_{\twoCat}(\twoCorr(\Ccal), \Dcal)$ and the subspace of $\Hom_{\twoCat}(\Ccal, \Dcal)$ consisting of those functors $F: \Ccal \rightarrow \Dcal$ which satisfy the left Beck-Chevalley condition.
\end{theorem}
\begin{proof}
 By virtue of proposition \ref{prop twocorr satisface bc} the functor $\iota_\Ccal$ factors through $\twoCorr^{\univ}(\Ccal)$. Our goal is to show that the resulting functor $Q: \twoCorr^{\univ}(\Ccal) \rightarrow \twoCorr(\Ccal)$ is an equivalence. Thanks to item (i) in lemma \ref{lemma identification hom} it  suffices show that for every pair of objects $c, c'$ in $\Ccal$, the induced functor 
\[
Q_*: \Hom_{\twoCorr^{\univ}(\Ccal)}(\iota_\Ccal^{\univ}(c), \iota_\Ccal^{\univ}(c')) \rightarrow \Hom_{\twoCorr(\Ccal)}(\iota_\Ccal(c), \iota_\Ccal(c'))
\]
is an equivalence.

Denote by $\mathcal{R}$ the equivalence  of item (ii) in  lemma \ref{lemma identification hom} between the two-sided fibration associated to the functor $\Hom_{\twoCorr^{\univ}(\Ccal)}(-, -)|_{\Ccal^\op \times \Ccal}$ and $\Funct(\Lambda^2_0, \Ccal)$. In particular $\mathcal{R}$ gives for every pair of objects $c, c'$ in $\Ccal$  an isomorphism
 \[
\Funct(\Lambda^2_0, \Ccal)_{(c, c')}   \xrightarrow{=}  \Hom_{\twoCorr^\univ(\Ccal)}(\iota_\Ccal^\univ(c), \iota_\Ccal^\univ(c')).
\]
Note that the left hand side is also equivalent to $\Hom_{\twoCorr(\Ccal)}(\iota_\Ccal(c), \iota_\Ccal(c'))$, so we have a (a priori non necessarily commutative) diagram of categories
\[
\begin{tikzcd}
\Hom_{\twoCorr^\univ(\Ccal)}(\iota_\Ccal^\univ(c), \iota_\Ccal^\univ(c')) \arrow{dr}[swap]{=} \arrow{rr}{Q_*} & & \arrow{dl}{=} \Hom_{\twoCorr(\Ccal)}(\iota_\Ccal(c), \iota_\Ccal(c')) \\ & \Funct(\Lambda^2_0,\Ccal)_{(c,c')} & .
\end{tikzcd}
\]
Since $\Cat$ is generated by the walking arrow, to show that $Q_*$ is an equivalence, it suffices to show that the diagram of spaces obtained from the above by applying the functor $\Hom_{\Cat}([1], - ): \Cat \rightarrow \Spc$, can be made commutative. In other words, we have to show that $Q_*$ is compatible with both equivalences above at the level of arrows.
 
 Consider a morphism of spans $\eta: T \rightarrow S$ depicted as follows.
\[
\begin{tikzcd}
  & t \arrow[ldd, "\alpha'"'] \arrow[rdd, "\beta'"] \arrow[d, "\mu"]  &    \\
  & s \arrow[ld, "\alpha"] \arrow[rd, "\beta"'] &    \\
c &                                             & c'
\end{tikzcd}
\]
Recall that the projection $(\ev_1, \ev_2): \Funct(\Lambda^2_0, \Ccal) \rightarrow \Ccal \times \Ccal$ is both cartesian and cocartesian. Passing to cartesian lifts of $\alpha, \alpha' , \mu$ and cocartesian lifts of $\mu, \beta, \beta'$ yields the following  diagram of categories:
\[
\begin{tikzcd}
                                                                                                              & { \Funct(\Lambda^2_0, \Ccal)_{(c, t)}} \arrow[rdd, "{(\id_c, \beta')_*}"] \arrow[d, shift left=0.25em, "{(\id_c, \mu)_*}", pos=0.8] &                                         \\
                                                                                                              & { \Funct(\Lambda^2_0, \Ccal)_{(c, s)}} \arrow[rd, "{(\id_c, \beta)_*}"'] \arrow[u, shift left=0.25em, "{(\id_c, \mu)^*}", pos=0.2]  &                                         \\
{ \Funct(\Lambda^2_0, \Ccal)_{(c, c)}} \arrow[ru, "{(\id_c, \alpha)^*}"'] \arrow[ruu, "{(\id_c, \alpha')^*}"] &                                                                                     & { \Funct(\Lambda^2_0, \Ccal)_{(c, c')}}
\end{tikzcd}
\]
Here $(\id_c,\mu)^*$ is right adjoint to $(\id_c, \mu)_*$. The morphism of spans $\eta$ is the image of the identity span of $c$ under the natural transformation
\[
(\id_c,\beta')_* (\id_c, \alpha')^* = (\id_c, \beta)_*(\id_c,\mu)_* (\id_c,\mu)^* (\id_c, \alpha)^* \rightarrow  (\id_c,\beta)_* (\id_c, \alpha)^*
\]
induced by the adjunction $(\id_c, \mu)_* \dashv (\id_c, \mu)^*$. Under the equivalence $\mathcal{R}$, the above diagram becomes
\[
\begin{tikzcd}[column sep = small]
                                                                                                              & { \Hom_{\twoCorr^\univ(\Ccal)}(\iota_\Ccal^\univ(c), \iota_\Ccal^\univ(t))} \arrow[rdd, "{\iota_\Ccal^\univ(\beta')_*}"] \arrow[d, shift left=0.25em, "{\iota_\Ccal^\univ(\mu)_*}", pos=0.8] &                                         \\
                                                                                                              & { \Hom_{\twoCorr^\univ(\Ccal)}(\iota_\Ccal^\univ(c), \iota_\Ccal^\univ(s))} \arrow[rd, "{\iota_\Ccal^\univ(\beta)_*}"'] \arrow[u, shift left=0.25em, "{\iota_\Ccal^\univ(\mu)^R_*}", pos=0.2]  &                                         \\
{ \Hom_{\twoCorr^\univ(\Ccal)}(\iota_\Ccal^\univ(c), \iota_\Ccal^\univ(c))} \arrow[ru, "{\iota_\Ccal^\univ(\alpha)^R_*}"'] \arrow[ruu, "{\iota_\Ccal^\univ(\alpha')^R_*}"] &                                                                                     & { \Hom_{\twoCorr^\univ(\Ccal)}(\iota_\Ccal^\univ(c), \iota_\Ccal^\univ(c'))}.
\end{tikzcd}
\]
Furthermore, the identity span of $c$ becomes $\id_{\iota_\Ccal^\univ(c)}$. We thus see that $\eta$ corresponds, under the equivalence $\mathcal{R}$, to the morphism
\[
\iota_\Ccal^\univ(\beta') \iota_\Ccal^\univ(\alpha')^R = \iota_\Ccal^\univ(\beta) \iota_\Ccal^\univ(\mu) \iota_\Ccal^\univ(\mu)^R \iota_\Ccal^\univ(\alpha)^R \rightarrow \iota_\Ccal^\univ(\beta) \iota_\Ccal^\univ(\alpha)^R
\]
induced by the counit of the adjunction $\iota_\Ccal^\univ(\mu) \dashv \iota_\Ccal^\univ(\mu)^R$. Applying the functor $Q$ recovers the morphism
\[
\iota_\Ccal(\beta')\iota_\Ccal(\alpha')^R = \iota_\Ccal(\beta) \iota_\Ccal(\mu) \iota_\Ccal(\mu)^R \iota_\Ccal(\alpha)^R \rightarrow \iota_\Ccal(\beta) \iota_\Ccal(\alpha)^R
\]
induced by the counit of the adjunction $\iota_\Ccal(\mu) \dashv \iota_\Ccal(\mu)^R$. This agrees with the image of $\eta$ under the usual isomorphism $\Hom_{\twoCorr(\Ccal)}(\iota_\Ccal(c), \iota_\Ccal(c')) = \Funct(\Lambda^2_0, \Ccal)_{(c, c')}$, as we wanted.
\end{proof}


%% file: Higher_correspondences.tex

\tableofcontents

\section{Higher categories of correspondences}

Let $\Ccal$ be a category admitting pullbacks. For each $n \geq 2$ one can construct an $n$-category $\nCorr(\Ccal)$ called the $n$-category of correspondences of $\Ccal$. The case $n = 2$ of this construction was the focus of section \ref{section twocorr}. For $n > 2$ the $n$-category $\nCorr(\Ccal)$ is defined inductively so that its objects agree with the objects of $\Ccal$, and for each pair of objects $c, c'$ in $\Ccal$, the hom $(n-1)$-category between them is $(n-1)\kr\Corr(\Hom_{2\kr\Corr(\Ccal)}(c, c'))$. Our goal in this section is to review the definition and main properties of the $n$-category of correspondences, and to establish two results (theorems \ref{teo univer prop ncorr} and \ref{teo extension}) that provide ways of constructing functors from $\nCorr(\Ccal)$ into a target $n$-category $\Dcal$.

We begin in \ref{subsection ncorr} by reviewing the definition of the $n$-category of correspondences. Here we depart from previous approaches in the literature: rather than defining $\nCorr(\Ccal)$ as an $(n-1)$-fold simplicial category (as in \cite{HaugSpan}), we use the language of enriched category theory to make sense of the fact that $\nCorr(\Ccal)$ is defined by applying the functor $(n-1)\kr\Corr$ at the level of hom categories on $\twoCorr(\Ccal)$. We also show that $n\kr\Corr(\Ccal)$ enjoys strong adjointness properties that extend those of $\twoCorr(\Ccal)$: every $k$-cell in $n\kr\Corr(\Ccal)$ with $k < n-1$ is both left and right adjointable, and its left and right adjoints coincide. In particular, if $\Ccal$ has a symmetric monoidal structure which is compatible with pullbacks, every object of $n\kr\Corr(\Ccal)$ is fully dualizable in the $(n-1)$-category underlying $n\kr\Corr(\Ccal)$.

We are ultimately interested in constructing functors out of $\nCorr(\Ccal)$: these are higher sheaf theories on $\Ccal$. In the case $n = 2$, a way of constructing functors is provided by theorem \ref{univ prop twocorr}: in order to construct a functor out of $\twoCorr(\Ccal)$ it suffices to construct a functor out of $\Ccal$ which satisfies the left Beck-Chevalley condition. In \ref{subsection higher BC} we introduce higher analogs of the Beck-Chevalley condition. While the ordinary Beck-Chevalley condition for a functor $F: \Ccal \rightarrow \Dcal$ involves the adjointability of certain commutative squares in $\Dcal$, the higher Beck-Chevalley condition is an inductively defined criterion that involves the adjointability of certain squares in $\Dcal$, together with the adjointability of certain squares in $\End_{\Dcal}(d)$ for objects $d$ in the image of $F$, together with the adjointability of certain squares in $\End_{\End_\Dcal(d)}(\id_d)$ for all such $d$, and so on. 

The first main result of this section is theorem \ref{teo univer prop ncorr}, which states that $\nCorr(\Ccal)$ is the universal $n$-category equipped with a functor from $\Ccal$ satisfying the left $(n-1)$-fold Beck-Chevalley condition. We also show that in the presence of a symmetric monoidal structure on a functor $F: \Ccal \rightarrow \Dcal$ which satisfies the left $(n-1)$-fold Beck-Chevalley condition, its extension to $\nCorr(\Ccal)$ also comes equipped with a symmetric monoidal structure.

Given a functor $F: \Ccal \rightarrow \Dcal$  into an $n$-category with colimits which is known to satisfy a higher Beck-Chevalley condition, one is sometimes interested in understanding whether the left Kan extension of $F$ to  the presheaf category $\Pcal(\Ccal)$ also satisfies a higher Beck-Chevalley condition. This is instrumental in the study of higher sheaf theories in algebraic geometry, as these start out life as functors on the category of affine schemes which are then extended to prestacks. In \ref{subsection extension} we apply the theory of conical colimits from \cite{Pres} to establish the second main result of this section, theorem \ref{teo extension}: under appropriate conditions on $\Dcal$, if a functor $F: \Ccal \rightarrow \Dcal$ is such that $F$ and $F^{n\dsh\op}$ satisfy the  $(n-1)$-fold Beck-Chevalley condition, then the left Kan extension $F': \Pcal(\Ccal) \rightarrow \Dcal$  satisfies the left $(n-1)$-fold Beck-Chevalley condition. These conditions are for instance verified in the case when $\Dcal$ underlies a presentable $n$-category (see \cite{Pres}).

\subsection{The $n$-category of correspondences}\label{subsection ncorr}

We begin by giving a construction of the $n$-category of correspondences of a category with pullbacks. We will use the language of enriched categories as developed in \cite{GH} and \cite{Hinich} - we refer the reader to \cite{Pres} for our conventions.

\begin{notation}
For each $n \geq 1$ we denote by $\nCat_\pb$ the category $(n-1)\kr\Catenr{\Cat_\pb}$ of $(n-1)$-categories enriched in the cartesian symmetric monoidal category $\Cat_\pb$.
\end{notation}

\begin{remark}\label{remark characterization ncatpb}
The inclusion $\Cat_\pb \rightarrow \Cat$ induces an inclusion $\nCat_\pb \rightarrow \nCat$ for each $n \geq 1$. The fact that the inclusion $\Cat_\pb \rightarrow \Cat$ creates limits implies that  inclusion $\nCat_\pb \rightarrow \Cat$ creates limits as well for all $n \geq 1$. Its image can be characterized inductively for $n \geq 2$ as follows:
\begin{itemize}
\item An $n$-category $\Ccal$ belongs to $\nCat_\pb$ if and only if for every pair of objects $c,c'$ in $\Ccal$ the $(n-1)$-category $\Hom_\Ccal(c, c')$ belongs to $(n-1)\kr\Cat_\pb$, and for every triple of objects $c, c', c''$ the composition map $\Hom_\Ccal(c, c') \times \Hom_\Ccal(c', c'') \rightarrow \Hom_\Ccal(c, c'')$ belongs to $(n-1)\kr\Cat_\pb$.
\item A functor of $n$-categories $F: \Ccal \rightarrow \Dcal$ belongs to $\nCat_\pb$ if and only if for every pair of objects $c, c'$ in $\Ccal$ the induced functor of $(n-1)$-categories $\Hom_\Ccal(c, c') \rightarrow \Hom_{\Dcal}(Fc, Fc')$ belongs to $(n-1)\kr\Cat_\pb$.
\end{itemize}
In particular, the subcategory $\Cat \subset \nCat$ is contained in $\nCat_\pb$ for all $n > 1$. Moreover, if $\Ccal$ is a category and $\Dcal$ is an object of $\nCat_\pb$ for $n > 1$, any functor $F: \Ccal \rightarrow \Dcal$ belongs to $\nCat_\pb$.
\end{remark}

\begin{proposition}
The functor $\twoCorr:\Cat_{\text{\normalfont pb}} \rightarrow \twoCat$ factors through $\twoCat_{\text{\normalfont\pb}}$.
\end{proposition}
\begin{proof}
Let $\Ccal$ be a category admitting pullbacks. We first show that $\twoCorr(\Ccal)$ belongs to $\twoCat_\pb$. Let $c, c'$ be objects in $\Ccal$ and recall from remark \ref{remark hom categories} that we have an equivalence $\Hom_{\twoCorr(\Ccal)}(c, c') = \Ccal_{/c, c'}$. Let $S = (c \leftarrow s \rightarrow c')$ be an object in $\Hom_{\twoCorr(\Ccal)}(c, c')$. The forgetful functor $\Ccal_{/c, c'} \rightarrow \Ccal$ induces an equivalence
\[
\Hom_{\twoCorr(\Ccal)}(c, c')_{/S} = \Ccal_{/s}.
\]
Since $\Ccal$ admits pullbacks, we conclude that the category $\Hom_{\twoCorr(\Ccal)}(c, c')_{/S}$ admits products for all $S$, which means that $\Hom_{\twoCorr(\Ccal)}(c, c')$ has pullbacks.

Let $c, c', c''$ be a triple of objects in $\Ccal$. Recall from remark \ref{remark hom categories} that the composition map for the Segal category $\twoCorr'(\Ccal)$ is equivalent to the composition of the right Kan extension functor 
\[
\Funct(\Tw([2])_\el, \Ccal) \rightarrow \Funct(\Tw([2]), \Ccal)
\]
 with the functor 
 \[
 \Funct(\Tw([2]), \Ccal) \rightarrow \Funct(\Tw([1]), \Ccal)
 \]
  of precomposition with the functor $\Tw([1]) \rightarrow \Tw([2])$ induced by the active arrow $[1] \rightarrow [2]$. Both of these preserve pullbacks, and therefore the composition map for $\twoCorr'(\Ccal)$ preserves pullbacks.  The composition map  
  \[
  \Hom_{\twoCorr(\Ccal)}(c, c') \times \Hom_{\twoCorr(\Ccal)}(c', c'') \rightarrow \Hom_{\twoCorr(\Ccal)}(c, c'')
  \] is obtained from the above by passing to fibers over $(c, c', c'')$ and therefore it also preserves pullbacks. It follows that $\Ccal$ satisfies the criteria of remark \ref{remark characterization ncatpb}, so it belongs to $\twoCat_\pb$, as desired.

Let $F: \Ccal \rightarrow \Dcal$ be a morphism in $\Cat_\pb$ and let $c, c'$ be objects in $\Ccal$. Let $S = (c \leftarrow s \rightarrow c')$ be an object in $\Hom_{\twoCorr(\Ccal)}(c, c')$. We have a commutative diagram of categories
\[
\begin{tikzcd}
\Hom_{\twoCorr(\Ccal)}(c, c')_{/S} \arrow{r}{} \arrow{d}{} & \Ccal_{/s} \arrow{d}{} \\
\Hom_{\twoCorr(\Dcal)}(Fc, Fc')_{/FS} \arrow{r}{}  & \Dcal_{/Fs}
\end{tikzcd}
\]
whose horizontal arrows are isomorphisms. Since $F$ preserves pullbacks, the right vertical arrow preserves products. Hence the left vertical arrow reserves products, and it follows that the morphism $\Hom_{\twoCorr(\Ccal)}(c, c') \rightarrow \Hom_{\twoCorr(\Dcal)}(Fc, Fc')$ induced by $F$ preserves pullbacks. Our result now follows from remark \ref{remark characterization ncatpb}.
\end{proof}

\begin{construction}\label{construction ncorr}
Let $n \geq 3$ and assume given a limit preserving functor 
\[
(n-1)\kr\Corr: \Cat_\pb \rightarrow (n-1)\kr\Cat
\] which factors through $(n-1)\kr\Cat_\pb$, and such that its composition with the truncation functor $(-)^{\leq 0} : (n-1)\kr\Cat \rightarrow \Spc$ is equivalent to the truncation functor $(-)^{\leq 0}: \Cat_\pb \rightarrow \Spc$. Consider the composite functor
\[
\Cat_\pb \xrightarrow{\twoCorr} \twoCat_\pb = \Catenr{\Cat_\pb} \xrightarrow{(n-1)\kr\Corr_!} \Algbrd((n-1)\kr\Cat).
\]
For each object $\Ccal$ in $\Cat_\pb$, the underlying Segal space to $(n-1)\Corr_!\twoCorr(\Ccal)$ is given by \[
(((n-1)\Corr)^{\leq 0})_!\twoCorr(\Ccal) = ((-)^{\leq 0})_!\twoCorr(\Ccal)
\]
which recovers the underlying Segal space of $\twoCorr(\Ccal)$. Since $\twoCorr(\Ccal)$ is a $2$-category, we conclude that $(n-1)\Corr_!\twoCorr(\Ccal)$ is an $n$-category. We denote by $\nCorr$ the resulting functor $\Cat_\pb \rightarrow \nCat$.
\end{construction}
 
 In the setting of construction \ref{construction ncorr}, the resulting functor $\nCorr$ is again limit preserving and factors through $\nCat_\pb$. Moreover, the functor $(\nCorr)^{\leq 0}$ is equivalent to $(\twoCorr)^{\leq 0}$ which is in turn equivalent to the truncation functor $(-)^{\leq 0}: \Cat_\pb \rightarrow \Spc$. We thus see that $\nCorr$ satisfies the hypothesis of construction \ref{construction ncorr} for $n + 1$.  Starting with the functor $\twoCorr$ we can thus produce functors $\nCorr:\Cat_\pb \rightarrow \nCat$ for every $n \geq 3$. In what follows, it will be convenient to allow $n$ to be $1$ as well, by setting $1\kr\Corr: \Cat_\pb \rightarrow \Cat$ to be the forgetful functor.
 
\begin{definition}
Let $\Ccal$ be a category admitting pullbacks and let $n \geq 1$. We call $\nCorr(\Ccal)$ the $n$-category of correspondences of $\Ccal$.
\end{definition}

We now study the relationship of the $n$-category of correspondences for different values of $n$. In what follows, we leave implicit the inclusions $\nCat \rightarrow (n+1)\kr\Cat$. In other words, we work in the category $\omega\kr\Cat = \colim_{n \geq 0} \nCat$, so that all the functors $\nCorr$ for different values of $n$ can be considered to have the same target.

\begin{construction}\label{constr natural transf ncorr}
Let $n \geq 3$ and assume given a natural transformation 
\[
\iota^{n-2, n-1}: (n-2)\kr\Corr \rightarrow (n-1)\kr\Corr.
\]
 Consider the induced natural transformation $((n-2)\kr\Corr)_! \rightarrow ((n-1)\kr\Corr)_!$. Composition with the functor $\twoCorr$ yields a natural transformation 
 \[
\iota^{n-1, n}: (n-1)\kr\Corr \rightarrow n\kr\Corr
 \]
  which we continue denoting by $\iota$. Applying this inductively starting with the natural transformation of construction \ref{constr natural transf} we obtain a sequence of functors $\Cat_\pb \rightarrow \omega\kr\Cat$ and natural transformations between them
\[
1\kr\Corr \rightarrow \twoCorr \rightarrow 3\kr\Corr \rightarrow \ldots.
\]
For each $m \leq n$ we let $\iota^{m, n}: m\kr\Corr \rightarrow n\kr\Corr$ the associated natural transformation. In the case $m = 1$ we will simply write $\iota^n  = \iota^{1, n}$.
\end{construction}

\begin{remark}\label{remark mono y truncation}
Let $n \geq m \geq 1$. Then the natural transformation $\iota^{m, n}$ is a monomorphism. Moreover, its composition with the truncation functor $(-)^{\leq m} : n\kr\Cat \rightarrow m\kr\Cat$ is an equivalence.
\end{remark}

\begin{remark}\label{remark dualizability}
We can think about the sequence of functors and transformations of construction \ref{constr natural transf ncorr} as a functor $-\kr\Corr: \NN \rightarrow \omegaCat$, where $\NN$ denotes the poset of natural numbers. Since $n\kr\Corr$ is limit preserving for all values of $n$, we have that $-\kr\Corr$ is limit preserving. 

In particular, if $\Ccal$ is a symmetric monoidal category admitting pullbacks and such that the functor $\Ccal \times \Ccal \rightarrow \Ccal$ preserves pullbacks, then $n\kr\Corr(\Ccal)$ inherits a symmetric monoidal structure for all $n$, and the functors $\iota_\Ccal^{m, n}$ also inherit symmetric monoidal structures. In the case when $\Ccal$ is a category with finite limits equipped with its cartesian symmetric monoidal structure, it follows from proposition \ref{prop dualizability} that every object of $\nCorr(\Ccal)$ is dualizable.
\end{remark}

\begin{proposition}\label{prop adjointability}
Let $\Ccal$ be a category admitting pullbacks. Let $n \geq 3$ and $1 \leq k < n-1$. Then every $k$-cell $\mu$ in $\nCorr(\Ccal)$ admits both a right adjoint $\mu^R$ and a left adjoint $\mu^L$, and moreover there is an equivalence $\mu^L = \mu^R$. 
\end{proposition}
\begin{proof}
We argue by induction on $k$. Consider first the case $k = 1$. By remark \ref{remark mono y truncation} we have that any arrow in $\nCorr(\Ccal)$ belongs to the image of the functor $\iota_{\Ccal}^{3, n} : 3\kr\Corr(\Ccal) \rightarrow \nCorr(\Ccal)$. It therefore suffices to consider the case $n = 3$. Recall from remark \ref{remark passing op} that there is a natural equivalence $\twoCorr = \twoCorr^{1\dsh\op}$ that restricts to the identity on objects. It follows that there is an equivalence $3\kr\Corr(\Ccal) = 3\kr\Corr(\Ccal)^{2\dsh\op}$ which restricts to the identity on objects and arrows. Therefore for every arrow $\mu$ in $3\kr\Corr(\Ccal)$, we have that $\mu$ admits a right adjoint if and only if it admits a left adjoint, and moreover in that case there is an equivalence $\mu^R = \mu^L$.

It now suffices to show that every arrow in $3\kr\Corr(\Ccal)$ can be written as a composition of arrows which admit either a left or right adjoint. Since the functor $\iota_\Ccal^{2, 3}: \twoCorr(\Ccal) \rightarrow 3\kr\Corr(\Ccal)$ is surjective on arrows, it suffices to show that this is the case in $\twoCorr(\Ccal)$. Indeed, any morphism in $\twoCorr(\Ccal)$ is represented by a span
\[
\begin{tikzcd}
  & s \arrow[ld, "\alpha"'] \arrow[rd, "\beta"] &    \\
c &                                             & c'
\end{tikzcd}
\]
and can thus be written as the composition $\iota_\Ccal(\beta) \iota_\Ccal^R(\alpha)$, which admit adjoints thanks to proposition \ref{prop adjunctos in twocrr}.

Assume now that $k > 1$. Let $c, c'$ be objects in $\Ccal$ such that $\mu$ is a $k$-cell with source and target objects $\iota_\Ccal^n(c)$ and $\iota_\Ccal^n(c')$, respectively. Then $\mu$ can be thought of as a $(k-1)$-cell in 
\[
\Hom_{n\kr\Corr(\Ccal)}(\iota_\Ccal^n(c), \iota_\Ccal^n(c')) = (n-1)\kr\Corr(\Ccal)(\Hom_{\twoCorr(\Ccal)}(\iota_\Ccal(c), \iota_\Ccal(c')).
\]
Our result now follows from the inductive hypothesis.
\end{proof}

The following result appears previously in \cite{HaugSpan} corollary 12.5.
\begin{corollary}\label{coro full dualizability}
Let $\Ccal$ be a category admitting finite limits equipped with its cartesian symmetric monoidal structure, and let $n \geq 2$. Then every object of the $(n-1)$-category $n\kr\Corr(\Ccal)^{\leq n-1}$ is fully dualizable.
\end{corollary}
\begin{proof}
Combine remark \ref{remark dualizability} with proposition \ref{prop adjointability}.
\end{proof}

\subsection{Higher Beck-Chevalley conditions}\label{subsection higher BC}

Our next goal is to generalize the universal property of theorem \ref{univ prop twocorr} to the case of higher categories of correspondences.

\begin{construction}\label{construction comm square}
Let $\Dcal$ be a $2$-category. For each right adjointable arrow $\alpha: d \rightarrow e$ in $\Dcal$ we denote by $\alpha^R$ its right adjoint, and by $\epsilon_\alpha: \alpha \alpha^R \rightarrow \id_d$ the counit of the adjunction.

 Consider a commutative diagram
\[
\begin{tikzcd}
d' \arrow{r}{\alpha'} \arrow{d}{\beta'} & d \arrow{d}{\beta} \\e' \arrow{r}{\alpha} & e
\end{tikzcd}
\]
in $\Dcal$ such that all four arrows admit right adjoints. The diagonal map $\gamma: d' \rightarrow e$ has a right adjoint which can be computed as in two ways as 
\[
\beta'^R \alpha^R = \gamma^R = \alpha'^R \beta^R
\]
and the counit of the adjunction can be described as
\[
(\alpha \epsilon_{\beta'} \alpha^R) \circ \epsilon_\alpha = \epsilon_\gamma  = (\beta \epsilon_{\alpha'} \beta^R) \circ  \epsilon_\beta.
\]
The above equivalence exhibits the following square in  $\End_\Ccal(e)$ as a commutative square:
\[
\begin{tikzcd}[row sep= large, column sep=large]
\gamma \gamma^R \arrow{d}{\alpha \epsilon_{\beta'} \alpha^R} \arrow{r}{\beta \epsilon_{\alpha'} \beta^R} & \beta \beta^R \arrow{d}{\epsilon_\beta}\\
\alpha \alpha^R \arrow{r}{\epsilon_\alpha} & \id_e.
\end{tikzcd}
\]
\end{construction}

\begin{definition}
Let $\Dcal$ be an $n$-category. We say that a commutative diagram
\[
\begin{tikzcd}
d' \arrow{r}{\beta'} \arrow{d}{\alpha'} & d \arrow{d}{\alpha} \\e' \arrow{r}{\beta} & e
\end{tikzcd}
\]
in $\Dcal$ is $1$-fold (vertically / horizontally) right adjointable if it is right adjointable in the  $2$-category underlying $\Dcal$. For each $k \geq 2$ we say that the above diagram is $k$-fold (vertically / horizontally ) right adjointable if it is (vertically / horizontally) right adjointable, all its arrows admit right adjoints, and the commutative diagram square
\[
\begin{tikzcd}[row sep= large, column sep=large]
\gamma \gamma^R \arrow{d}{\alpha \epsilon_{\beta'} \alpha^R} \arrow{r}{\beta \epsilon_{\alpha'} \beta^R} & \beta \beta^R \arrow{d}{\epsilon_\beta}\\
\alpha \alpha^R \arrow{r}{\epsilon_\alpha} & \id_e.
\end{tikzcd}
\]
in $\End_\Dcal(e)$ defined in construction \ref{construction comm square}, is $(k-1)$-fold (vertically/ horizontally) right adjointable.
\end{definition}

\begin{construction}\label{construction nfold right adjoint in limit}
Recall the universal left adjointable arrow $L : [1] \rightarrow \Adj$ and the universal vertically right adjointable square 
\[
\id_{[1]} \times L : [1] \times [1] \rightarrow [1] \times \Adj
\] from remark \ref{right adjoint in limit}. Let $\Ucal^+$ be the colimit in $\twoCat$ of the following diagram:
\[
\begin{tikzcd}
\Adj &                                                                           & {[1]} \times \Adj &                                                                            & \Adj \\
     &                                                                           & {[1] \times [1]} \arrow[u, "\id_{[1]} \times L"] &                                                                            &      \\
     & {[1]} \arrow[ru, "{\id_{[1]} \times \lbrace 0 \rbrace}"] \arrow[luu, "L"] &                            & {[1]} \arrow[lu, "{\lbrace 1 \rbrace \times \id_{[1]}}"'] \arrow[ruu, "L"'] &     
\end{tikzcd}
\]
In other words, $\Ucal^+$ is the universal vertically right adjointable square such that the horizontal arrows are also right adjointable. Construction \ref{construction comm square} provides a functor 
\[
C_{\text{univ}}: \Sigma([1] \times [1]) \rightarrow \Ucal^+
\]
 where $\Sigma$ denotes the functor which associates to each $n$-category $\Tcal$ an $(n+1)$-category $\Sigma(\Tcal)$ with objects $s, t$ and such that its only nontrivial Hom category is $\Hom_{\Sigma(\Tcal)}(s, t) = \Tcal$.

We define for each $n \geq 2$ an $n$-category $\Ucal_n$ equipped with a functor $i_n: [1] \times [1] \rightarrow \Ucal_n$ as follows:
\begin{itemize}
\item When $n = 2$ we set $i_1 = \id_{[1]} \times L : [1] \times [1] \rightarrow [1] \times \Adj$.
\item Assume that $n > 2$. Then we let $\Ucal_n$ be the pushout of the following diagram:
\[
\begin{tikzcd}
\Ucal^+ &                                                                                     & \Sigma(\Ucal_{n-1}) \\
        & {\Sigma([1]\times[1])} \arrow[ru, "\Sigma(i_{n-1})"] \arrow[lu, "C_{\text{univ}}"'] &                    
\end{tikzcd}
\].
\end{itemize}

The functor $i_n$ is the universal vertically $(n-1)$-fold right adjointable commutative square. The map $i_n$ an epimorphism - for each $n$-category $\Dcal$ precomposition with $i_n$ induces an equivalence between the space of functors $\Ucal_n \rightarrow \Dcal$ and the space of functors $[1] \times [1] \rightarrow \Dcal$ which correspond to $(n-1)$-fold right adjointable squares.
\end{construction} 

\begin{remark}\label{remark nfold right adjoint in limit}
As a consequence of the existence of a universal vertically $(n-1)$-fold  we deduce the following fact: if $F: \Ical \rightarrow n\kr\Cat$ is a diagram with limit $\Dcal$, then a commutative square in $\Dcal$ is $(n-1)$-fold vertically right adjointable if and only if its image in $F(i)$ is $(n-1)$-fold vertically right adjointable for all $i$ in $\Ical$.
\end{remark}

\begin{definition}
Let $\Ccal$ be a category admitting pullbacks and $\Dcal$ be an $n$-category. We say that a functor $F: \Ccal \rightarrow \Dcal$ satisfies the left $n$-fold Beck-Chevalley condition if for every cospan $x \rightarrow s \leftarrow y$ in $\Ccal$, the induced commutative square
\[
\begin{tikzcd}
F(x \times_s y) \arrow{r}{} \arrow{d}{} & F(y) \arrow{d}{} \\ F(x) \arrow{r}{} & F(s)
\end{tikzcd}
\]
is $n$-fold right adjointable.
\end{definition}

We are now ready to state the universal property of the $n$-category of correspondences.
\begin{theorem}\label{teo univer prop ncorr}
Let $\Ccal$ be a category admitting pullbacks and $\Dcal$ be an $n$-category. Restriction along the inclusion $\iota^n_\Ccal: \Ccal \rightarrow \nCorr(\Ccal)$ induces an identification of $\Hom_{\nCat}(\nCorr(\Ccal), \Dcal)$ with the subspace of $\Hom_{\nCat}(\Ccal,\Dcal)$ consisting of functors satisfying the left $(n-1)$-fold Beck-Chevalley condition.
\end{theorem}

The proof of theorem \ref{teo univer prop ncorr} needs a few lemmas.

\begin{lemma}\label{lemma epis algebras}
Let $\Ocal$ be an operad \footnote{In this paper we use a language for speaking about operads which is close in spirit to the classical language in terms of objects and operations which satisfy a composition rule. Namely, given an operad $\Ocal$ with associated category of operators $p: \Ocal^\otimes \rightarrow \Fin_*$, we call $p^{-1}(\langle 1 \rangle)$ the category of objects of $\Ocal$, and  arrows in $\Ocal^\otimes$ lying above an active arrow of the form $\langle n \rangle \rightarrow \langle 1 \rangle$ are called operations of $\Ocal$.} and let $\Mcal$ be an $\Ocal$-monoidal category. Let $F: \Acal \rightarrow \Bcal$ be a morphism of $\Ocal$-algebras in $\Mcal$. Assume that for every operation $\mu$ in  $\Ocal$ with source $\lbrace X_s \rbrace_{s \in S}$ and target $X$  the induced map 
\[
F(\mu): \mu(\lbrace\Acal(X_s)\rbrace_{s \in S}) \rightarrow \mu(\lbrace\Bcal(X_s)\rbrace_{s \in S})
\]
 is an epimorphism in $\Mcal(X)$. Then
\begin{enumerate}[\normalfont (i)]
\item The morphism $F$ is an epimorphism in $\Alg_\Ocal(\Mcal)$.
\item Let $\Bcal'$ be another $\Ocal$-algebra in $\Mcal$. Then a morphism $F': \Acal \rightarrow \Bcal'$ factors through $\Bcal$ if and only if for every object $X$ in $\Ocal$ the map $F'(X): \Acal(X) \rightarrow \Bcal'(X)$ factors through $\Bcal(X)$.
\end{enumerate}
\end{lemma}
\begin{proof}
Denote by $p: \Ocal^\otimes \rightarrow \Fin_*$ and $q: \Mcal^\otimes \rightarrow \Ocal^\otimes$ the categories of operators associated to $\Ocal$ and $\Mcal$. Recall that $\Alg_\Ocal(\Mcal)$ is the full subcategory of 
\[
\Funct_{\Ocal^\otimes}(\Ocal^\otimes, \Mcal^\otimes) = \Funct(\Ocal^\otimes , \Mcal^\otimes) \times_{\Funct(\Ocal^\otimes, \Ocal^\otimes)} [0]
\]
on those functors that map inert arrows in $\Ocal^\otimes$ to inert arrows in $\Mcal^\otimes$. The algebras $\Acal, \Bcal$ determine functors $\Acal, \Bcal: \Ocal^\otimes \rightarrow \Mcal^\otimes$, and $F$ is a natural transformation $\Acal \rightarrow \Bcal$.

We will show that $F$ is an epimorphism by showing that the diagram
\[
\begin{tikzcd}
\Acal \arrow{r}{F} \arrow{d}{F} & \Bcal \arrow{d}{\id} \\ \Bcal  \arrow{r}{\id} & \Bcal
\end{tikzcd}
\]
is a pushout diagram in $\Alg_\Ocal(\Mcal)$. This would follow if we are able to show that the above is a pushout in $\Funct_{\Ocal^\otimes}(\Ocal^\otimes, \Mcal^\otimes)$. Note that the square which is constant $\id_{\Ocal^\otimes}$ is a pushout square in $\Funct(\Ocal^\otimes, \Ocal^\otimes)$. Therefore we may in fact restrict to showing that the above square defines a pushout square in $\Funct(\Ocal^\otimes, \Mcal^\otimes)$. This can be done pointwise. Let $\lbrace X_s \rbrace_{s \in S}$ be an object of $\Ocal^\otimes$, corresponding to a finite collection of objects of $\Ocal$. We must show that the square
\[
\begin{tikzcd}[column sep = large, row sep = large]
\lbrace \Acal(X_s) \rbrace_{s \in S} \arrow{r}{\lbrace F(X_s) \rbrace_{s \in S}} \arrow{d}{\lbrace F(X_s) \rbrace_{s \in S}}  & \lbrace \Bcal(X_s) \rbrace_{s \in S} \arrow{d}{\id} \\ \lbrace \Bcal(X_s) \rbrace_{s \in S} \arrow{r}{\id} & \lbrace \Bcal(X_s) \rbrace_{s \in S}
\end{tikzcd}
\]
taking place in  $q^{-1}(\lbrace X_s \rbrace_{s \in S}) = \prod_{s \in S} \Mcal(X_s) \subset \Mcal^\otimes$,  defines a pushout square in $\Mcal^\otimes$. Since the image under $q$ of the above square is constant (and thus a pushout), suffices to show that the above is in fact a $q$-pushout. Using \cite{HTT} proposition 4.3.1.10, we reduce to showing that for every arrow $\alpha: \lbrace X_s \rbrace_{s \in S} \rightarrow \lbrace X_t \rbrace_{t \in T}$ in $\Ocal^\otimes$, the image of the above square under the functor
\[
\alpha_!: q^{-1}(\lbrace X_s \rbrace_{s \in S}) \rightarrow q^{-1}(\lbrace X_t \rbrace_{t \in T})
\]
is a pushout. It suffices to do this in the case when $\alpha$ is either inert or active, and this follows readily from our assumptions on $F$. This proves item (i).

We now prove item (ii). Assume given a morphism $F': \Acal \rightarrow \Bcal'$ such that for every object $X$ in $\Ocal$ the map $F'(X) : \Acal(X) \rightarrow \Bcal'(X)$ factors through $\Bcal(X)$. We have to show that $F'$ factors through $\Bcal$. Let $\lbrace X_s \rbrace_{s \in S}$ be an object in $\Ocal^\otimes$. Note that the map
\[
F'(\lbrace X_s \rbrace_{s \in S}) : \lbrace \Acal(X_s )\rbrace_{s \in S} \rightarrow \lbrace \Bcal'(X_s) \rbrace_{s \in S}
\]
factors through $\lbrace \Bcal(X_s )\rbrace_{s \in S}$. It follows that we have a pushout diagram
\[
\begin{tikzcd}[column sep = large, row sep = large]
\lbrace \Acal(X_s) \rbrace_{s \in S} \arrow{r}{\lbrace F(X_s) \rbrace_{s \in S}} \arrow{d}{\lbrace F'(X_s) \rbrace_{s \in S}}  & \lbrace \Bcal(X_s) \rbrace_{s \in S} \arrow{d}{} \\ \lbrace \Bcal'(X_s) \rbrace_{s \in S} \arrow{r}{\id} & \lbrace \Bcal'(X_s) \rbrace_{s \in S}
\end{tikzcd}
\]
in $q^{-1}(\lbrace X_s \rbrace_{s \in S})$. We claim that the above is a $q$-colimit diagram. As before, by \cite{HTT} proposition 4.3.1.10, we reduce to showing that for every arrow $\alpha: \lbrace X_s \rbrace_{s \in S} \rightarrow \lbrace X_t \rbrace_{t \in T}$ in $\Ocal^\otimes$, the image of the above square under $\alpha_!$ is a pushout. The case when $\alpha$ is inert is clear; the case when $\alpha$ is active follows from our assumption on $F$.

Using \cite{HA} lemma 3.2.2.9 we see that there is a square
\[
\begin{tikzcd}
\Acal \arrow{r}{F} \arrow{d}{F'} & \Bcal \arrow{d}{} \\ \Bcal' \arrow{r}{G} & \Bcal''
\end{tikzcd}
\]
in $\Funct_{\Ocal^\otimes}(\Ocal^\otimes, \Mcal^\otimes)$ with the property that, for every object $\lbrace X_s \rbrace_{s \in S}$ in $\Ocal^\otimes$, the induced square
\[
\begin{tikzcd}[column sep = huge, row sep = large]
\lbrace \Acal(X_s) \rbrace_{s \in S} \arrow{r}{F(\lbrace X_s \rbrace_{s \in S})} \arrow{d}{F'(\lbrace X_s \rbrace_{s \in S})} & \lbrace \Bcal(X_s) \rbrace_{s \in S} \arrow{d}{} \\ \lbrace \Bcal'(X_s) \rbrace_{s \in S} \arrow{r}{G(\lbrace X_s \rbrace_{s \in S})} & \lbrace \Bcal''(X_s) \rbrace_{s \in S}
\end{tikzcd}
\]
is a $q$-colimit diagram. It follows that $G(\lbrace X_s \rbrace_{s \in S})$ is an isomorphism, and hence $G$ is an isomorphism, which implies that $F'$ factors through $\Bcal$, as desired.
 \end{proof}

\begin{lemma}\label{lemma identify square}
Let $\Ccal$ be a category admitting pullbacks and let 
\[
\hspace{0.5cm}
\begin{tikzcd}
x\times_s y \arrow{r}{\alpha'} \arrow{d}{\beta'} & y \arrow{d}{\beta} \\ x \arrow{r}{\alpha} & s
\end{tikzcd} \hspace{0.5cm}(\ast)
\]
be a cartesian square in $\Ccal$. Then the commutative square
\[
\begin{tikzcd}[column sep = 7em]
\iota_\Ccal(\gamma) \iota^R_\Ccal(\gamma) \arrow{r}{\iota_\Ccal(\beta) \epsilon_{\iota_{\Ccal}(\alpha')} \iota_{\Ccal}^R(\beta)} \arrow{d}{\iota_\Ccal(\alpha)\epsilon_{\iota_\Ccal(\beta')} \iota_\Ccal^R(\alpha)} & \iota_\Ccal(\beta) \iota_\Ccal^R(\beta) \arrow{d}{\epsilon_{\iota_\Ccal(\beta)}} \\ \iota_\Ccal(\alpha)\iota_\Ccal(\alpha)^R \arrow{r}{\epsilon_{\iota_\Ccal(\alpha)}} & \id_{\iota_\Ccal(s)}
\end{tikzcd}
\]
in $\End_{\twoCorr(\Ccal)}(\iota_\Ccal(s))$ that results from construction \ref{construction comm square} applied to the image of $(\ast)$ under $\iota_\Ccal$ is equivalent to the image under the canonical map
\[
\Ccal_{/s} \rightarrow \Ccal_{/s, s} = \End_{\twoCorr(\Ccal)}(\iota_\Ccal(s))
\]
of the commutative square
\[
\begin{tikzcd}
\gamma \arrow{r}{\overline{\beta'}} \arrow{d}{\overline{\alpha'}} & \alpha \arrow{d}{\overline{\alpha}} \\ \beta \arrow{r}{\overline{\beta}} & \id_s
\end{tikzcd}
\]
whose image under the forgetful functor $\Ccal_{/s} \rightarrow \Ccal$ recovers $(\ast)$.
\end{lemma}
\begin{proof}
Let $\Ucal$ be the universal cospan: this is the category with objects $0, 1, 2$ and nontrivial arrows $0 \xrightarrow{\alpha_u} 1 \xleftarrow{\beta_u} 2$. Everything in the statement is functorial in $\Ccal$, so we may reduce to the case when $\Ccal = \Fcal(\Ucal)$ is the free category with pullbacks on $\Ucal$, and the cartesian square under consideration is
\[
\begin{tikzcd}
0\times_1 2 \arrow{r}{\alpha'_u} \arrow{d}{\beta'_u} & 2 \arrow{d}{\beta_u} \\ 0 \arrow{r}{\alpha_u} & 1.
\end{tikzcd}
\]
The two squares that we have to show are equivalent can easily be seen to have the same vertices. Our result now follows from the fact that in $\End_{\twoCorr(\Fcal(\Ucal))}(\id_1)$ there is  a unique commutative square with those vertices.
\end{proof}

\begin{lemma}\label{lemma check homs}
Let $\Ccal$ be a category admitting pullbacks and let $\Dcal$ be an $n$-category for $n \geq 3$. Let $F: \Ccal \rightarrow \Dcal$ be a functor satisfying the left Beck-Chevalley condition. Then $F$ satisfies the left $(n-1)$-fold Beck-Chevalley condition if and only if for every pair of objects $z, w$ in $\Ccal$ the induced functor\[
\Hom_{\twoCorr(\Ccal)}(\iota_\Ccal(z), \iota_\Ccal(w)) \rightarrow \Hom_\Dcal(F(z), F(w))\]
 satisfies the left $(n-2)$-fold Beck-Chevalley condition.
\end{lemma}
\begin{proof}
Assume first that for every pair of objects $z, w$ in $\Ccal$ the functor 
\[
\Hom_{\twoCorr(\Ccal)}(\iota_\Ccal(z), \iota_\Ccal(w)) \rightarrow \Hom_\Dcal(F(z), F(w))
\] satisfies the left $(n-2)$-fold Beck-Chevalley condition. Consider a cartesian square
\[
\begin{tikzcd}
x\times_s y \arrow{r}{\alpha'} \arrow{d}{\beta'} & y \arrow{d}{\beta} \\ x \arrow{r}{\alpha} & s
\end{tikzcd}
\]
in $\Ccal$ and denote by $\gamma: x\times_s y \rightarrow s$ the induced map. Since $F$ is assumed to satisfy the left Beck-Chevalley condition, we know that the induced commutative square
\[
\begin{tikzcd}
F(x\times_s y) \arrow{r}{F(\alpha')} \arrow{d}{F(\beta')} & F(y) \arrow{d}{F(\beta)} \\ F(x) \arrow{r}{F(\alpha)} & F(s)
\end{tikzcd}
\]
is right adjointable.  We need to show that the commutative square
\[
\begin{tikzcd}[row sep = large, column sep=7em]
F(\gamma)F(\gamma)^R \arrow{r}{F(\beta) \epsilon_{F(\alpha')} F(\beta)^R} \arrow{d}{F(\alpha)\epsilon_{F(\beta')} F(\alpha)^R} & F(\beta) F(\beta)^R \arrow{d}{\epsilon_{F(\beta)}} \\ F(\alpha)F(\alpha)^R \arrow{r}{\epsilon_{F(\alpha)}} & \id_{F(s)}
\end{tikzcd}
\]
in $\End_\Dcal(F(s))$ arising from construction \ref{construction comm square} is $(n-2)$-fold right adjointable. This square is the image under the functor $\twoCorr(\Ccal) \rightarrow \Dcal$ of the commutative square
\[
\begin{tikzcd}[column sep = 7em]
\iota_\Ccal(\gamma) \iota^R_\Ccal(\gamma) \arrow{r}{\iota_\Ccal(\beta) \epsilon_{\iota_{\Ccal}(\alpha')} \iota_{\Ccal}^R(\beta)} \arrow{d}{\iota_\Ccal(\alpha)\epsilon_{\iota_\Ccal(\beta')} \iota_\Ccal^R(\alpha)} & \iota_\Ccal(\beta) \iota_\Ccal^R(\beta) \arrow{d}{\epsilon_{\iota_\Ccal(\beta)}} \\ \iota_\Ccal(\alpha)\iota_\Ccal(\alpha)^R \arrow{r}{\epsilon_{\iota_\Ccal(\alpha)}} & \id_{\iota_\Ccal(s)}.
\end{tikzcd}
\]
This commutative square is described by lemma \ref{lemma identify square}. In particular, note that the image of it under the canonical map 
\[
\End_{\twoCorr(\Ccal)}(\id_s) = \Ccal_{/s, s} \rightarrow \Ccal
\]
is a cartesian square. Since the projection $\Ccal_{/s, s} \rightarrow \Ccal$ creates pullbacks, we conclude that the above is in fact a cartesian square, and therefore its image inside $\End_\Dcal(F(s))$ is $(n-2)$-fold right adjointable, as we wanted.

Assume now that $F$ satisfies the left $(n-1)$-fold Beck-Chevalley condition, and let $z, w$ be a pair of objects of $\Ccal$. Let
\[\hspace{0.75cm}
\begin{tikzcd}
x^\sharp \times_{s^\sharp} y^{\sharp}  \arrow{r}{\beta'^\sharp} \arrow{d}{\alpha'^\sharp} & y^\sharp \arrow{d}{\alpha^\sharp} \\ x^\sharp \arrow{r}{\beta^\sharp} & s^\sharp
\end{tikzcd}
\hspace{0.5cm} (\star)
\]
be a cartesian square in $\Hom_{\twoCorr(\Ccal)}(\iota_\Ccal(z), \iota_\Ccal(w))$, whose image under the forgetful functor $\Hom_{\twoCorr(\iota_\Ccal(z), \iota_\Ccal(w))} = \Ccal_{/z, w} \rightarrow \Ccal$ is a cartesian square in $\Ccal$ which we denote as follows:
\[
 \begin{tikzcd}
x \times_{s} y  \arrow{r}{\beta'} \arrow{d}{\alpha'} & y \arrow{d}{\alpha} \\ x \arrow{r}{\beta} & s
\end{tikzcd} 
\]
Let $z \leftarrow s : \sigma $ and $\tau: s \xrightarrow{\tau} w$ be the legs of the span $s^\sharp$. Then the square $(\star)$ is equivalent to the image of the square
\[\hspace{2cm}
\begin{tikzcd}
\gamma \arrow{r}{\overline{\beta'}} \arrow{d}{\overline{\alpha'}} & \alpha \arrow{d}{\overline{\alpha}} \\ \beta \arrow{r}{\overline{\beta}} & \id_s
\end{tikzcd}\hspace{0.5cm} (\star \star)
\]
under the composite map
\[
\Ccal_{/s} \rightarrow \Ccal_{/s, s} = \End_{\twoCorr(\Ccal)}(\iota_\Ccal(s)) \xrightarrow{\iota_\Ccal^R(\sigma)^*} \Hom_{\twoCorr(\Ccal)}(\iota_\Ccal(z), \iota_\Ccal(s))\xrightarrow{ \iota_\Ccal(\tau)_*} \Hom_{\twoCorr(\Ccal)}(\iota_\Ccal(z),\iota_\Ccal( w)).
\] 
It therefore suffices to show that the image of the square $(\star \star)$ under the composite map
\[
\Ccal_{/s} \rightarrow \Ccal_{/s, s} = \End_{\twoCorr(\Ccal)}(\iota_\Ccal(s)) \rightarrow \End_{\Dcal}(F(s))
\]
is $(n-2)$-fold right adjointable. This follows from the fact that $F$ satisfies the left $(n-1)$-fold Beck-Chevalley condition, combined with lemma \ref{lemma identify square}.
\end{proof}

\begin{proof}[Proof of theorem \ref{teo univer prop ncorr}]
We argue by induction on $n$. The case $n = 2$ is theorem \ref{univ prop twocorr}, so we assume $n \geq 3$. The inclusion $\iota_\Ccal^{2, n}: \twoCorr(\Ccal) \rightarrow \nCorr(\Ccal) $ defines a morphism in $\Algbrd_{\Ccal^{\leq 0}}((n-1)\kr\Cat)$ which is an epimorphism by item (i) in lemma \ref{lemma epis algebras} combined with our inductive hypothesis. The second part of lemma \ref{lemma epis algebras} guarantees that for any other  $\Assos_{\Ccal^{\leq 0}}$-algebra $\Bcal'$ in $(n-1)\kr\Cat$ precomposition with $\iota_{\Ccal}^{2, n}$ induces an equivalence between the space of maps $\nCorr(\Ccal) \rightarrow \Bcal'$ and the space of maps $\twoCorr(\Ccal) \rightarrow \Bcal'$ such that for every pair of objects $z, w$ in $\Ccal$ the induced map
\[
\Hom_{\twoCorr(\Ccal)}(\iota_\Ccal(z), \iota_\Ccal(w)) \rightarrow \Bcal'(z, w)
\]
satisfies the left $(n-2)$-fold Beck-Chevalley condition. Since the projection 
\[
\Algbrd((n-1)\kr\Cat) \rightarrow \Cat
\]
 is a cartesian fibration, we conclude that $\iota_{\Ccal}^{2, n}$ is in fact an epimorphism in  $\Algbrd((n-1)\kr\Cat)$, and in particular it is an epimorphism in $n\kr\Cat$. Moreover, a morphism $G: \twoCorr(\Ccal) \rightarrow \Dcal$ factors through $\nCorr(\Ccal)$ if and only if for every pair of objects $z, w$ in $\Ccal$ the induced functor
 \[
 \Hom_{\twoCorr(\Ccal)}(\iota_\Ccal(z), \iota_\Ccal(w)) \rightarrow \Hom_{\Dcal}(G(\iota_\Ccal(z)), G(\iota_\Ccal(w)))
 \]
 satisfies the left $(n-2)$-fold Beck-Chevalley condition. Our result now follows from a combination of theorem \ref{univ prop twocorr} and lemma \ref{lemma check homs}. 
\end{proof}

\begin{remark}\label{remark univer prop sym monoidal}
Let $\Ccal$ be a symmetric monoidal category which admits pullbacks and assume such that the tensor product functor $\Ccal \times \Ccal \rightarrow \Ccal$ preserves pullbacks. Combining theorem \ref{teo univer prop ncorr} with lemma \ref{lemma epis algebras} we see that for every $n \geq 2$ the functor $\iota^n_\Ccal: \Ccal \rightarrow \nCorr(\Ccal)$ is an epimorphism of symmetric monoidal $n$-categories. Given a symmetric monoidal $n$-category $\Dcal$, precomposition with $\iota_\Ccal$ yields an equivalence between the space of symmetric monoidal functors $\nCorr(\Ccal) \rightarrow \Dcal$ and the space of symmetric monoidal functors $\Ccal \rightarrow \Dcal$ which satisfy the left $(n-1)$-fold Beck-Chevalley condition.
\end{remark}

\subsection{Extension along the Yoneda embedding} \label{subsection extension}

Our next goal is to show that under suitable conditions, a functor out of $\nCorr(\Ccal)$ can be extended to a functor on $(n-1)\kr\Corr(\Pcal(\Ccal))$. We begin by recalling the passage to adjoints property from \cite{Pres} section 5.5. The definition requires the notion of conical colimits and limits in $n$-categories - we refer the reader to \cite{Pres}  section 4.3 for a general discussion of this concept in the setting of enriched categories.

\begin{notation}
Let $\Dcal$ be an $n$-category. We denote by $\Dcal^{\leq 1}$ the $1$-category underlying $\Dcal$, and by $(\Dcal^{\leq 1})^{\text{radj}}$ (resp. $\Dcal^{\leq 1})^{\text{ladj}}$ the subcategory of $\Dcal^{\leq 1}$ containing all objects, and only those morphisms which are right (resp. left) adjointable in $\Dcal$.
\end{notation}

\begin{definition}\label{def passage right adj prop}
Let $\Dcal$ be an $n$-category. We say that $\Dcal$ satisfies the passage to adjoints property if the following conditions are satisfied:
\begin{itemize}
\item The category $(\Dcal^{\leq 1})^{\text{\normalfont radj}}$ has all colimits, and the inclusion $(\Dcal^{\leq 1})^{\text{\normalfont radj}} \rightarrow \Dcal$ preserves conical colimits.
\item The category $(\Dcal^{\leq 1})^{\text{\normalfont ladj}}$ has all limits, and the inclusion $(\Dcal^{\leq 1})^{\text{\normalfont ladj}} \rightarrow \Dcal$ preserves  conical limits.
\end{itemize}
\end{definition}

\begin{remark}
Let $\Dcal$ be an $n$-category. Then passage to adjoints induces an equivalence between the categories $(\Dcal^{\leq 1})^\text{radj}$ and $(\Dcal^{\leq 1})^{\text{ladj}}$. It follows that if $\Dcal$ satisfies the passage to adjoints property, then a right adjointable diagram $F: \Ical^\rhd \rightarrow \Dcal$ in $\Dcal$ is a conical colimit if and only if the diagram $F^R: (\Ical^\op)^\lhd \rightarrow \Dcal$ is a conical limit.
\end{remark}

\begin{definition}
Let $\Dcal$ be an $n$-category. We say that $\Dcal$ satisfies the $1$-fold passage to adjoints property if it satisfies the passage to  adjoints property. For each $(n-1) \geq k \geq 2$ we say that $\Dcal$ satisfies the $k$-fold passage to  adjoints property if it satisfies the passage to adjoints property and for every pair of objects $d, e$ in $\Dcal$ the $(n-1)$-category $\Hom_\Dcal(d, e)$ satisfies the $(k-1)$-fold passage to adjoints property.
\end{definition}

\begin{definition}
Let $\Dcal$ be an $n$-category. We say that $\Dcal$ is $1$-fold conically cocomplete if it admits all small conical colimits. For each $n \geq k \geq 2$ we say that $\Dcal$ is $k$-fold conically cocomplete if it is $1$-fold conically cocomplete and for every pair of objects $d, e$ in $\Dcal$ the $(n-1)$-category $\Hom_{\Dcal}(d, e)$ is $(k-1)$-fold conically cocomplete.
\end{definition}

\begin{example}
Recall the category $n\kr\Pr^L$ of presentable $n$-categories from \cite{Pres}. Each object $D$ in $n\kr\Pr^L$ has an underlying $n$-category $\psi_n(D)$. It follows from \cite{Pres} theorems 5.4.6 and 5.5.14 and remark 5.3.3 that $\psi_n(D)$ is $n$-fold conically cocomplete and satisfies the $(n-1)$-fold passage to adjoints property. 
\end{example}

\begin{definition}
Let $\Ccal$ be a category admitting pullbacks. We say that a map $f: x \rightarrow y$ in $\Pcal(\Ccal)$ is representable if for every map $c \rightarrow y$ with $c$ in $\Ccal$, the presheaf $x \times_y c$ is representable.
\end{definition}

\begin{remark}
Let $\Ccal$ be a category admitting pullbacks. The class of representable morphisms defines a subcategory $\Pcal(\Ccal)_{\text{rep}}$ of $\Pcal(\Ccal)$. The category $\Pcal(\Ccal)_{\text{rep}}$ contains all pullbacks, and these are preserved by the inclusion into $\Pcal(\Ccal)$.
\end{remark}

We are now ready to state our extension theorem.

\begin{theorem}\label{teo extension}
Let $n \geq 2$. Let $\Ccal$ be a category admitting pullbacks and $\Dcal$ be an $(n-1)$-fold conically cocomplete $n$-category satisfying the $(n-1)$-fold passage to adjoints property. Let $F: \Pcal(\Ccal) \rightarrow \Dcal$ be a conical colimit preserving functor.
\begin{enumerate}[\normalfont (i)]
\item  If $F|_{\Ccal}$ satisfies the left $(n-1)$-fold Beck-Chevalley condition, then for every pair of maps $\alpha: x \rightarrow s$ and $\beta: y \rightarrow s$ in $\Pcal(\Ccal)$ where $\beta$ is representable the commutative square
\[
\begin{tikzcd}
F(x \times_s y) \arrow{r}{} \arrow{d}{} & F(y) \arrow{d}{F(\beta)} \\ F(x) \arrow{r}{F(\alpha)} & F(s)
\end{tikzcd}
\]
is $(n-1)$-fold vertically right adjointable. In particular, $F|_{\Pcal(\Ccal)_{\text{\normalfont rep}}}$ satisfies the left $(n-1)$-fold Beck-Chevalley condition.
\item  If $F|_{\Ccal}: \Ccal \rightarrow \Dcal$ and $(F|_{\Ccal})^{n\dsh\op}: \Ccal \rightarrow \Dcal^{n\dsh\op}$  satisfy the left $(n-1)$-fold Beck-Chevalley conditions, then $F$ satisfies the left $(n-1)$-fold beck-Chevalley condition.
\end{enumerate}
\end{theorem}

Before giving the proof of theorem \ref{teo extension}, we study a few  consequences.

\begin{corollary}\label{coro also have right}
Let $n \geq 3$. Let $\Ccal$ be a category admitting pullbacks and let $\Dcal$ be an $(n-2)$-fold conically cocomplete $n$-category satisfying the $(n-2)$-fold passage to adjoints property. Let $F:\Pcal(\Ccal) \rightarrow \Dcal$ be a conical colimit preserving functor such that $F|_\Ccal$ satisfies the left $(n-1)$-fold Beck-Chevalley condition. Then $F$ satisfies the left $(n-2)$-fold Beck-Chevalley condition.
\end{corollary}
\begin{proof}
By theorem \ref{teo univer prop ncorr}, $F$ admits an extension to $\nCorr(\Ccal)$. Recall from remark \ref{remark passing op} that there is a natural equivalence $\twoCorr = \twoCorr^{1\dsh \op}$ which restricts to the identity on underlying spaces. It follows by induction that for each $m \geq 3$ there is a natural equivalence between the functor $m\kr\Corr$ and the functor $m\kr\Corr^{(m- 1)\dsh\op}$, which is an equivalence upon passage to underlying categories. In particular, we have an equivalence between $\nCorr(\Ccal)^{(n-1) \dsh \op}$ and $\nCorr(\Ccal)$ which is the identity on the underlying category. It follows that $F^{(n-1)\dsh\op}$ admits an extension to a functor from $\nCorr(\Ccal)$ into $\Dcal^{(n-1)\dsh\op}$, and therefore $F^{(n-1)\dsh\op}$ satisfies the $(n-1)$-fold Beck-Chevalley condition. Our result now follows from theorem \ref{teo extension}.
\end{proof}

\begin{remark}
In the setting of theorem \ref{teo extension} item (i), it is not true in general that $F$ will satisfy the full left $(n-1)$-fold  Beck-Chevalley condition, although it does satisfy the left $(n-2)$-fold Beck-Chevalley condition thanks to corollary \ref{coro also have right}. In other words, the statement of part (i) in  theorem \ref{teo extension} does not hold if we remove the representability condition. For example, the functor $\operatorname{mod}: \CAlg(\Sp)^\op \rightarrow (\Prscr^L)^{1\dsh\op}$ that assigns to each commutative ring spectrum its category of modules satisfies the left Beck-Chevalley condition, but its left Kan extension along $\CAlg(\Sp)^\op \rightarrow \Pcal(\CAlg(\Sp))^\op$ is the functor that assigns to each (nonconnective) prestack its category of quasicoherent sheaves, which does not satisfy base change in all generality.
\end{remark}

\begin{definition}
Let $\Dcal$ be an $n$-category. For each $0 \leq k < n$ we denote by $\Dcal^{k\text{\normalfont -adj}}$ the largest $(k+1)$-category contained inside $\Dcal$ and such that all cells of dimension at most $k$ admit both left and right adjoints.  We say that an arrow $\alpha$ in $\Dcal$ is $k$-fold adjointable if it belongs to $\Dcal^{k\text{\normalfont -adj}}$.
\end{definition}

\begin{corollary}
Let $n \geq 3$. Let $\Ccal$ be a category admitting pullbacks  and let $\Dcal$ be an $(n-1)$-fold conically cocomplete $n$-category satisfying the $(n-1)$-fold passage to adjoints property. Let $F:\Pcal(\Ccal) \rightarrow \Dcal$ be a  conical colimit preserving functor such that $F|_\Ccal$ and $(F|_{\Ccal})^{n\dsh\op}$ satisfy the left  $(n-1)$-fold Beck-Chevalley condition.
Let $\beta: x \rightarrow y$ be a morphism in $\Pcal(\Ccal)$. Then  $F(\beta)$ is $(n-2)$-fold adjointable, and its left and right adjoints are equivalent.
\end{corollary}
\begin{proof}
Combine theorems \ref{teo univer prop ncorr} and \ref{teo extension} with proposition \ref{prop adjointability}.
\end{proof}

\begin{corollary}
Let $n \geq 2$. Let $\Ccal$ be a category admitting pullbacks. Let $\Dcal$ be an $(n-1)$-fold conically cocomplete symmetric monoidal $n$-category satisfying the $(n-1)$-fold passage to adjoints property. Equip $\Pcal(\Ccal)$ with its cartesian symmetric monoidal structure. Let $F: \Pcal(\Ccal) \rightarrow \Dcal$ be a symmetric monoidal conical colimit preserving functor such that $F|_\Ccal$ and $(F|_{\Ccal})^{n \dsh \op}$ satisfy the left  $(n-1)$-fold Beck-Chevalley condition. Let $x$ be an object in $\Pcal(\Ccal)$. Then
\begin{enumerate}[\normalfont (i)]
\item The object $F(x)$ in $\Dcal$ is fully dualizable in $\Dcal^{\leq n-1}$.
\item Let $\Delta: x \rightarrow x \times x$ be the diagonal map and $\pi: x \rightarrow \ast$ be the projection to the final object of $\Pcal(\Ccal)$. If $F(\Delta)$ and $F(\pi)$ are $(n-1)$-fold adjointable then $x$ is a fully dualizable object of $\Dcal$.
\end{enumerate}
\end{corollary}
\begin{proof}
Combine theorem \ref{teo extension}, remark \ref{remark univer prop sym monoidal} and corollary \ref{coro full dualizability}, with the description of the self duality of $x$ in $\twoCorr(\Pcal(\Ccal))$ from proposition \ref{prop dualizability}  and remark \ref{remark description dualizable cont}.
\end{proof}

The proof of theorem \ref{teo extension} needs a few lemmas.

\begin{lemma}\label{lemma bc in representables}
Let $\Ccal$ be a category admitting pullbacks and $\Dcal$ be a $2$-category.  Let $F: \Pcal(\Ccal) \rightarrow \Dcal$ be a conical colimit preserving functor.  Let  
\[\hspace{0.5cm}
\begin{tikzcd}
x \times_s y \arrow{d}{\beta'} \arrow{r}{\alpha' } & y \arrow{d}{\beta} \\ x \arrow{r}{\alpha} & s
\end{tikzcd} \hspace{0.5cm} (\star)
\]
be a cartesian square in $\Pcal(\Ccal)$. Assume that for every arrow $c \rightarrow c'$ in $\Ccal_{/s}$ the induced square
\[
\begin{tikzcd}
F(c \times_s y) \arrow{d}{F(\beta_c)} \arrow{r}{} & F(c'\times_s y) \arrow{d}{F(\beta_{c'})} \\ F(c) \arrow{r}{} & F(c')
\end{tikzcd}
\]
is vertically right adjointable. Then the image of $(\star)$ under $F$ is vertically right adjointable.
\end{lemma}
\begin{proof}
Let $X^\rhd: \Ical^\rhd \rightarrow \Pcal(\Ccal)$ be a colimit diagram such that its value on the cone point $X(\ast)$ recovers $x$, and such that $X = X^\rhd|_{\Ical}$ factors through $\Ccal$. Let 
\[
(\Ical^+)^\rhd = (\Ical \sqcup [0])^\rhd = \Ical^\rhd \bigcup_{[0]} [1]
\] be the category obtained by adjoining a final object to $\Ical \sqcup [0]$.  There is an evident inclusion 
\[
i: (\Ical^+)^\rhd \rightarrow \Ical^\rhd \times [1]
\]
defined by the fact that it preserves final objects, and extends the inclusions $\Ical \xrightarrow{\id_{\Ical} \times 1} \Ical \times [1] $ and $[0] \xrightarrow{\ast \times 0}\Ical \times [1]$. Consider the functor $(X^\rhd)^+: (\Ical^+)^\rhd \rightarrow \Pcal(\Ccal)$ which extends $X^\rhd$ and such that $(X^\rhd)^+|_{[0]^\rhd}$ recovers the arrow $\beta'$. 

Let $Y^\rhd : \Ical^\rhd \times [1] \rightarrow \Pcal(\Ccal)$ be the right Kan extension of $(X^\rhd)^+$ along $i$, and let $Y = Y^\rhd |_{\Ical \times [1]}$. In other words, $Y^\rhd$ is such that the restriction $Y^\rhd|_{\ast \times [1]}$ recovers the arrow $\beta'$, the restriction $Y^\rhd|_{\Ical \times \lbrace 1 \rbrace}$ recovers $X^\rhd$, and for every index $i$ in $\Ical$ the induced commutative square
\[
\begin{tikzcd}
(Y(i), 0) \arrow{d}{} \arrow{r}{} & x\times_s y \arrow{d}{\beta'} \\
(Y(i),1) \arrow{r}{} & x
\end{tikzcd}
\]
is a pullback square. Since pullbacks distribute over colimits in presheaf categories, we have that $Y^\rhd|_{\Ical \times \lbrace j \rbrace}$ is a colimit diagram in $\Pcal(\Ccal)$ for $j = 0, 1$.

Let $Z^\rhd: \Ical^\rhd \rightarrow \Funct([1], \Ccal)$ be the functor associated to $Y^\rhd$, and $Z = Z^\rhd|_{\Ical}$. Note that the functor $FZ^\rhd: \Ical^\rhd \rightarrow \Funct([1], \Dcal)$ is a conical colimit diagram for $FZ$.  Our hypothesis on $F$ together with remark \ref{remark adjointability square lax} imply that the functor $FZ$ factors through $\Funct(\Adj, \Dcal)$. Since $\Dcal $ is conically cocomplete, we conclude that $FZ^\rhd$ also factors through $\Funct(\Adj, \Dcal)$, and in fact defines a conical colimit diagram in $\Funct(\Adj,\Dcal)$.  In particular, for all objects $i$ in $\Ical$, the square
\[
\begin{tikzcd}
F(X(i) \times_x x\times_s y) \arrow{d}{F(\beta_{X(i)})} \arrow{r}{} & F(x \times_s y) \arrow{d}{F(\beta')} \\
F(X(i)) \arrow{r}{} & F(x)
\end{tikzcd}
\]
is vertically right adjointable.

Repeating the above argument with the conical colimit diagram $\Ccal_{/s}^\rhd \rightarrow \Pcal(\Ccal)$ shows that for every $i$ in $\Ical$ the square 
\[
\begin{tikzcd}
F(X(i) \times_x x\times_s y) \arrow{d}{F(\beta_{X(i)})} \arrow{r}{} & F(y) \arrow{d}{F(\beta )} \\
F(X(i)) \arrow{r}{} & F(s)
\end{tikzcd}
\]
is vertically right adjointable. In other words, the morphism $F(\beta') \rightarrow F(\beta)$ in $\Funct([1], \Dcal)$ induced by the square $(\star)$ is such that its composition with the map $F(\beta_X(i)) \rightarrow F(\beta')$ factors through the subcategory $\Funct(\Adj, \Dcal)$. Since $FZ^\rhd$ is a conical colimit diagram we conclude that the map $F(\beta') \rightarrow F(\beta)$ factors through $\Funct(\Adj,\Dcal)$, which means that $(\star)$ is vertically right adjointable, as desired.
\end{proof}

\begin{lemma}\label{lemma n3}
Let $\Ccal$ be a category admitting pullbacks, and $\Dcal$ be a $2$-category. Assume that $\Dcal$ is satisfies the passage to adjoints property. Let $F: \Pcal(\Ccal) \rightarrow \Dcal$ be a conical colimit preserving functor and assume that $F|_{\Ccal}$ satisfies the left and right Beck-Chevalley conditions. Then $F$ satisfies the left Beck-Chevalley condition.
\end{lemma}
\begin{proof}
Let
\[
\begin{tikzcd}
x \times_s y \arrow{d}{\beta'} \arrow{r}{\alpha' } & y \arrow{d}{\beta} \\ x \arrow{r}{\alpha} & s
\end{tikzcd}
\]
be a cartesian square in $\Pcal(\Ccal)$. We have to show that its image under $F$ is vertically right adjointable. By virtue of lemma \ref{lemma bc in representables}, we may assume that $x$ and $s$ are representable.

Let $X^\rhd: \Ical^\rhd \rightarrow \Pcal(\Ccal)$ be a colimit diagram such that $X = X^\rhd|_\Ical$ factors through $\Ccal$, and its value on the cone point $X^\rhd(\ast)$ recovers $y$. Since $F$ preserves conical colimits we have that $FX^\rhd$ is a conical colimit diagram.

Note that the composite morphism
\[
F(X(i)) \rightarrow F(y) \xrightarrow{F(\beta)} F(s)
\]
admits a right adjoint for every $i$ in $\Ical$. Since $\Dcal$ satisfies the passage to adjoints property and $FX$ factors through $(\Dcal^{\leq 1})^{\text{radj}}$,  we conclude that $F(\beta)$ admits a right adjoint in $\Dcal$. A similar argument guarantees that $F(\beta')$ admits a right adjoint.

Using remark \ref{lemma left adj iff right adj} we now see that in order to show that the image of our square under $F$ is vertically right adjointable it suffices to show that it is horizontally left adjointable. Applying the dual version of lemma \ref{lemma bc in representables} we may assume that our square takes place in $\Ccal$, and in this case the conclusion follows from the fact that $F|_{\Ccal}$ satisfies the right Beck-Chevalley condition.
\end{proof}

\begin{lemma}\label{lemma cube}
Let $\Dcal$ be an $n$-category for $n \geq 3$ and let 
\[\hspace{0.5cm}
\begin{tikzcd}
\tau_{d'} \arrow{r}{\tau_{\alpha'}} \arrow{d}{\tau_{\beta'}} & \tau_{d} \arrow{d}{\tau_{\beta}} \\ \tau_{e'} \arrow{r}{\tau_{\alpha}} & \tau_e
\end{tikzcd} \hspace{0.5cm}(\star)
\]
be a commutative square in $\Funct([1], \Dcal)$, corresponding to a commutative cube as follows:
\[
    \begin{tikzcd}[row sep=1.5em, column sep = 1.5em]
        & d'_1 \arrow[dd, "\beta'_1", pos=0.75] \arrow[rr, "\alpha'_1"] &&
    d_1 \arrow[dd,"\beta_1"] \\
    d'_0 \arrow[rr,"\alpha'_0", pos=0.7] \arrow[ur, "\tau_{d'}"] \arrow[dd, "\beta'_0"] &&
    d_0 \arrow[dd, "\beta_0" ,pos = 0.25] \arrow[ur,"\tau_d"] \\
 & e'_1 \arrow[rr, "\alpha_1", pos=0.25] && e_1 \\
    e'_0 \arrow[rr, "\alpha_0"] \arrow[ur, "\tau_{e'}"] && e_0 \arrow[ur, "\tau_{e}"]
    \end{tikzcd}
\]
Then the square $(\star)$ is $(n-1)$-fold vertically right adjointable if and only if the front and back faces of the cube are $(n-1)$-fold vertically right adjointable, the left and right faces are vertically right adjointable, and the top and bottom faces are horizontally right adjointable.
\end{lemma}
\begin{proof}
It follows from remark \ref{remark adjointability square lax} that:
\begin{itemize}
\item The morphisms $\tau_{\alpha'}$ and $\tau_\alpha$ are right adjointable if and only if the top and bottom faces of the cube are horizontally right adjointable.
\item The morphisms $\tau_{\beta'}$ and $\tau_{\beta}$ are right adjointable if and only if the left and right faces of the cube are vertically right adjointable.
\end{itemize}

Assume that both both of the items above hold. It follows from lemma \ref{lemma check invertibility 2 cell in funct} that the square $(\star)$ is vertically right adjointable if and only if the front and back faces of the cube are vertically right adjointable. Consider now the pullback square
\[
\begin{tikzcd}
\End_{\Funct([1],\Dcal)}(\tau_e) \arrow{r}{\ev_0}  \arrow{d}{\ev_1} & \End_{\Dcal}(e_0) \arrow{d}{(\tau_e)_*} \\  \arrow{r}{(\tau_e)^*} \End_{\Dcal}(e_1)  & \Hom_{\Dcal}(e_0, e_1).
\end{tikzcd}
\]
Using remark  \ref{remark nfold right adjoint in limit} we see that a commutative square in $\End_{\Funct([1],\Dcal)}(\tau_e)$ is $(n-2)$-fold vertically right adjointable if and only if its image under $\ev_0$ and $\ev_1$ is $(n-2)$-fold vertically right adjointable. The lemma now follows by applying this fact to the square associated to $(\star)$ under construction \ref{construction comm square}.
\end{proof}

\begin{lemma}\label{lemma colim pres}
Let $\Ccal$ be a category admitting pullbacks and $\Dcal$ be a $2$-category. Assume that $\Dcal$  satisfies the passage to adjoints property. Let $F: \Pcal(\Ccal) \rightarrow \Dcal$ be a conical colimit preserving functor satisfying the left Beck-Chevalley condition. Let $c, c'$ be objects of $\Pcal(\Ccal)$. Then the induced functor 
\[
F_* : \Hom_{\twoCorr(\Pcal(\Ccal))}(\iota_{\Pcal(\Ccal)}(c), \iota_{\Pcal(\Ccal)}(c')) \rightarrow \Hom_\Dcal(F(c), F(c'))
\]
is colimit preserving.
\end{lemma}
\begin{proof}
We first show that $F_*$ preserves (infinite) coproducts. Consider a family of spans
\[
\begin{tikzcd}
  & x_i \arrow[ld, "a_i"'] \arrow[rd, "b_i"] &    \\
c &                                          & c'
\end{tikzcd}
\]
in $\Pcal(\Ccal)$, indexed by a set $\Ical$. Let
\[
\begin{tikzcd}
  & x \arrow[ld, "a"'] \arrow[rd, "b"] &    \\
c &                                          & c'
\end{tikzcd}
\]
be their coproduct, and for each $i$ denote by $j_i: x_i \rightarrow x$ the induced map. We must show that the family of morphisms  
\[
F(b_i) F(a_i)^R = F(b) F(j_i) F(j_i)^R F(a)^R \rightarrow   F(b) F(a)^R
\]
 induced from the counit maps $F(j_i) F(j_i)^R \rightarrow \id_{F(x)}$ 
 exhibit $F(b) F(a)^R$ as the coproduct of the family $\lbrace F(b_i)F(a_i)^R \rbrace_{i \in \Ical}$  inside $\Hom_{\Dcal}(F(c), F(c'))$. 
 
 Since $\Dcal$ satisfies the passage to adjoints property, we have that the functors
 \[
 (F(j_i)^R)_*: \Hom_{\Dcal}(F(c), F(x)) \rightarrow \Hom_{\Dcal}(F(c), F(x_i)) 
 \]
 exhibit the left hand side as the product of the categories on the right hand side. It follows that the natural transformations $ (F(j_i) F(j_i)^R)_* \rightarrow \id_{\Hom_{\Dcal}(F(c), F(x))}$ induced from the counits of the adjunctions $F(j_i) \dashv F(j_i)^R$ exhibit  $\id_{\Hom_{\Dcal}(F(c), F(x))}$ as the coproduct of the endofunctors $(F(j_i) F(j_i)^R)_*$. Since $F(b)_*$ is right adjointable we see that the induced natural transformations $(F(b) F(j_i) F(j_i)^R)_* \rightarrow F(b)_*$ exhibit $F(b)_*$ as the coproduct of the functors $F(b) F(j_i) F(j_i)^R = F(b_i) F(j_i)^R$. Our claim now follows from this by evaluation at $F(a)^R$.
 
 We now show that $F_*$ preserves pushouts. Consider a cospan of spans
 \[
 \begin{tikzcd}
                                         & x_0 \arrow[ld, "\mu"'] \arrow[rd, "\nu"] \arrow[ldd, "a_0", pos=0.3] \arrow[rdd, "b_0", swap, pos = 0.3] &                                         \\
x_1 \arrow[d, "a_1"'] \arrow[rrd, "b_1", pos= 0.3] &                                                                                 & x_2 \arrow[lld, "a_2", swap, pos=0.3] \arrow[d, "b_2"] \\
c                                        &                                                                                 & c'                                     
\end{tikzcd}
\]
and let
\[
\begin{tikzcd}
  & x \arrow[ld, "a"'] \arrow[rd, "b"] &    \\
c &                                          & c'
\end{tikzcd}
\]
be its pushout. For each $i$ denote by $j_i : x_i \rightarrow x$ the associated map. We must show that the induced square
\[
\begin{tikzcd}
F(b_0) F(a_0)^R \arrow{r}{} \arrow{d}{} & F(b_2) F(a_2)^R \arrow{d}{} \\
F(b_1) F(a_1)^R \arrow{r}{} & F(b) F(a)^R
\end{tikzcd}
\]
in $\Hom_{\Dcal}(F(c), F(c'))$ is a pushout square. As before, since $\Dcal$ satisfies the passage to adjoints property, we have that the induced commutative square of categories
\[
\begin{tikzcd}
\Hom_\Dcal(F(c), F(x)) \arrow{r}{(F(j_2)^R)_*} \arrow{d}{(F(j_1)^R)_*} & \Hom_\Dcal(F(c), F(x_2)) \arrow{d}{F(\nu)^R_*} \\ \Hom_\Dcal(F(c), F(x_1)) \arrow{r}{F(\mu)^R_*} & \Hom_\Dcal(F(c), F(x_0))
\end{tikzcd}
\]
is a pullback square. It follows that the commutative square
\[
\begin{tikzcd}
(F(j_0)F(j_0)^R)_* \arrow{r}{} \arrow{d}{} & (F(j_2)F(j_2)^R)_* \arrow{d}{} \\ (F(j_1)F(j_1)^R)_* \arrow{r}{} & \id_{\Hom_\Dcal(F(c), F(x))} 
\end{tikzcd}
\]
induced from the counit of the adjunctions $F(j_i) \dashv F(j_i)^R$, is a pushout square. Since $F(b)_*$ is a left adjoint, we obtain that the following commutative square in the category $\Funct(\Hom_\Dcal(F(c), F(x)), \Hom_\Dcal(F(c), F(c')))$ is a pushout square:
\[
\begin{tikzcd}
(F(b)F(j_0)F(j_0)^R)_* \arrow{r}{} \arrow{d}{} & (F(b)F(j_2)F(j_2)^R)_* \arrow{d}{} \\ (F(b)F(j_1)F(j_1)^R)_* \arrow{r}{} & F(b)_*
\end{tikzcd}
\]
Our result now follows from this by evaluation at $F(a)^R$.
\end{proof}

\begin{proof}[Proof of theorem \ref{teo extension}]
We begin by proving item (ii). We argue by induction on $n$. The case $n = 2$ is proven in lemma \ref{lemma n3}, so assume $n > 2$.

It follows from our inductive hypothesis (in the form of corollary \ref{coro also have right}) that $F$ satisfies the left Beck-Chevalley condition. Let
\[
\begin{tikzcd}
x \times_s y \arrow{d}{\beta'} \arrow{r}{\alpha' } & y \arrow{d}{\beta} \\ x \arrow{r}{\alpha} & s
\end{tikzcd}
\]
be a cartesian square in $\Pcal(\Ccal)$. Let $S^\rhd : \Ical^\rhd \rightarrow \Pcal(\Ccal)$ be a colimit diagram such that its value on the cone point $S(\ast)$ recovers $s$, and $S = S^\rhd|_\Ical$ factors through $\Ccal$. Arguing as in the proof of lemma \ref{lemma bc in representables}, we may construct a functor $T^\rhd: \Ical^\rhd \times ([1] \times [1])$ such that $T^\rhd|_{\Ical^\rhd \times \lbrace(1, 1)\rbrace}$ recovers $S^\rhd$, the functor $T^\rhd|_{\ast \times ([1] \times [1])}$ is given by the above square, and moreover for every object $i$ in $\Ical$ the cube
\[
    \begin{tikzcd}[row sep=1.5em, column sep = 1.5em]
        & x\times_s y \arrow[dd, "\beta'", pos=0.75] \arrow[rr, "\alpha'"] &&
    y \arrow[dd,"\beta"] \\
    T^\rhd(i, (0,0)) \arrow[rr] \arrow[ur] \arrow[dd] &&
    T^\rhd(i, (1, 0)) \arrow[dd] \arrow[ur] \\
 & x \arrow[rr, "\alpha", pos=0.25] && s \\
    T^\rhd(i, (0, 1)) \arrow[rr] \arrow[ur] && T^\rhd(i,(1,1)) \arrow[ur]
    \end{tikzcd}
\]
has all its faces cartesian. In particular, for every object $(j, k)$ in $[1] \times [1]$ the restriction $T^\rhd|_{\Ical^\rhd \times \lbrace (j, k) \rbrace}$ is a colimit diagram in $\Pcal(\Ccal)$.

Let $R^\rhd : \Ical^\rhd \rightarrow \Funct([1]\times [1], \Pcal(\Ccal))$ be the induced functor. The composite functor $FR^\rhd : \Ical^\rhd \rightarrow \Funct([1], [1], \Dcal)$ is a conical colimit diagram, and its value on the cone point $\ast$ recovers the image under $F$ of our original square.

Recall from construction \ref{construction nfold right adjoint in limit} the universal $(n-1)$-fold vertically right adjointable square $i_n : [1] \times [1] \rightarrow \Ucal_n$.  Our task is to show that $FR^\rhd(\ast)$ belongs to $\Funct(\Ucal_n, \Dcal)$. This would follow if we are able to show that $FR^\rhd|_{\Ical}$ factors through $\Funct(\Ucal_n, \Dcal)$. Using lemma \ref{lemma cube} together with the fact that $F$ satisfies the left Beck-Chevalley condition, we reduce to showing that $FR^\rhd(i)$ belongs to $\Funct(\Ucal_n, \Dcal)$ for every $i$ in $\Ical$. In other words, we must show that the front face of the above cube is $(n-1)$-fold vertically right adjointable. This face is a cartesian square in $\Pcal(\Ccal)$, and its lower right corner $c = T^\rhd(i,(1,1))$ is representable.

It follows from lemma \ref{lemma identify square} that the commutative square in  $\End_{\twoCorr(\Pcal(\Ccal))}(\iota_{\Pcal(\Ccal)}(c))$ associated to the front face of the cube by construction \ref{construction comm square} is cartesian. It therefore suffices to show that the functor 
\[
\End_{\twoCorr(\Pcal(\Ccal))}(\iota_{\Pcal(\Ccal)}(c)) \rightarrow \End_\Dcal(F(c))
\]
induced from $F$ satisfies the left $(n-2)$-fold Beck-Chevalley condition. Using \cite{GHN} corollary 9.9, we obtain equivalences
\[
\Pcal(\End_{\twoCorr(\Ccal)}(\iota_\Ccal(c))) = \Pcal(\Ccal_{/c, c}) = \Pcal(\Ccal)_{/c, c} = \End_{\twoCorr(\Pcal(\Ccal))}(\iota_{\Pcal(\Ccal)}(c)).
\]
Our result now follows from our inductive hypothesis together with lemmas \ref{lemma colim pres} and \ref{lemma check homs}.

We now give a proof of item (i). We proceed along the same lines as the proof of item (ii) carried out above, with suitable modifications. Again we argue by induction on $n$. The case $n = 2$ follows directly from lemma \ref{lemma bc in representables}. Assume now that $n > 2$. Using our inductive hypothesis (in the form of corollary \ref{coro also have right}) we see  that $F$ satisfies the left Beck-Chevalley condition. Let
\[
\begin{tikzcd}
x \times_s y \arrow{d}{\beta'} \arrow{r}{\alpha' } & y \arrow{d}{\beta} \\ x \arrow{r}{\alpha} & s
\end{tikzcd}
\]
be a cartesian square in $\Pcal(\Ccal)$, where $\beta$ is  representable. We have to show that its image under $F$ is $(n-1)$-fold vertically right adjointable. Arguing in the same way as in the proof of item (ii), we may reduce to the case when $s$ belongs to $\Ccal$. Since $\beta$ is representable, we also have that $y$ belongs to $\Ccal$.

 As before, the square in $\End_{\twoCorr(\Pcal(\Ccal))}(\iota_{\Pcal(\Ccal)}(s))$  associated to the above square via construction \ref{construction comm square}, is cartesian. It follows from lemma \ref{lemma identify square} that the right vertical side of that square corresponds to the morphism of spans
\[
\begin{tikzcd}
  & y \arrow[ldd, "\beta"'] \arrow[rdd, "\beta"] \arrow[d, "\beta"] &   \\
  & s \arrow[ld, "\id"] \arrow[rd, "\id"']                          &   \\
s &                                                                 & s .
\end{tikzcd}
\]
Our result now follows again from our inductive hypothesis combined with  lemmas \ref{lemma colim pres} and \ref{lemma check homs}, using the fact that the above diagram defines a representable morphism in $\Pcal(\Ccal_{/s, s})$.
\end{proof}
